\documentclass[10pt]{amsart}
\usepackage[cp1251]{inputenc}
\usepackage[english]{babel}
\usepackage{amsmath}
\usepackage{amssymb}
\usepackage{amsfonts}

\usepackage[linktocpage=true, colorlinks=true, linkcolor=blue, citecolor=blue, urlcolor=blue]{hyperref}

\begin{document}
\renewcommand{\refname}{References}

\thispagestyle{empty}

\title[Comparative Analysis of the Efficiency of Application]
{Comparative Analysis of the Efficiency of Application of Legendre Polynomials 
and Trigonometric Functions to the Numerical Integration of Ito Stochastic 
Differential Equations}
\author[D.F. Kuznetsov]{Dmitriy F. Kuznetsov}
\address{Dmitriy Feliksovich Kuznetsov
\newline\hphantom{iii} Peter the Great Saint-Petersburg Polytechnic University,
\newline\hphantom{iii} Polytechnicheskaya ul., 29,
\newline\hphantom{iii} 195251, Saint-Petersburg, Russia}%
\email{sde\_kuznetsov@inbox.ru}

\maketitle {\small
\begin{quote}
\noindent{\sc Abstract.} 
The article is devoted to comparative analysis of the efficiency of 
application of Legendre polynomials and trigonometric functions to the 
numerical integration of Ito stochastic differential equations in the 
framework of the method of approximation of iterated Ito and Stratonovich 
stochastic integrals 
based on generalized multiple Fourier series. 
On the example of iterated Ito stochastic integrals of multiplicities 
1 to 3 from the Taylor--Ito expansion it is shown that expansions 
of stochastic integrals based on Legendre polynomials are
essentially simpler
and require significantly less computational costs compared to their 
analogues obtained using the trigonometric system of functions.
The results of the article can be useful for construction
of high-order 
strong numerical methods for Ito stochastic differential equations.
\medskip
\end{quote}
}

\vspace{12mm}


\setlength{\baselineskip}{1.8em}

\tableofcontents

\setlength{\baselineskip}{1.2em}


\vspace{2mm}

\section{Introduction}

\vspace{5mm}

In a lot of author's publications \cite{1}-\cite{new-art-1xxys} 
the mean-square 
approximation method for iterated Ito and Stratonovich stochastic integrals 
based on generalized multiple Fourier series is proposed and developed
(see Theorems 1--8 below).
Further, we will call this method as the method of generalized 
multiple Fourier series.
Under the term "generalized multiple Fourier series" we understand 
the Fourier series constructed using various complete orthonormal systems 
of functions in the space $L_2([t, T])$, and not only using the trigonometric 
system of functions. Here $[t, T]$ is an interval of integration 
of iterated Ito or Stratonovich stochastic integrals.

It is well known the another approach to series expansion of 
stochastic processes  
using eigenfunctions of their covariance operators
(the so-called Karhunen--Loeve 
expansion) \cite{18}. 
If the stochastic process is the Brownian bridge process on the time 
interval $[t, T]$, then the eigenfunctions of its covariance operator
will be trigonometric functions which form a complete orthonormal 
system of functions in the space $L_2([t, T])$ \cite{19}.
This means that the basis functions in the mentioned approach
can only be 
trigonometric functions.
In \cite{19}-\cite{23} the series expansion of the 
Brownian 
bridge process
was used for the expansion and mean-square approximation of iterated Ito 
and Stratonovich stochastic integrals.
Further, we will call this expansion as the Milstein expansion.

As mentioned above, in contrast to the Milstein expansion
the method of generalized 
multiple Fourier 
series \cite{1}-\cite{new-art-1xxys}
(see Theorems 1, 2 below) allows to use different systems of basis functions.
Thus, we can set the problem of choice the optimal system of basis 
functions within the framework of the method of generalized multiple Fourier 
series.
Some ideas on the solution of the
mentioned problem were given in a number of the author's works
\cite{3}-\cite{7}, \cite{9a}-\cite{10aaa}.

For example, in \cite{3}-\cite{7}, \cite{9a}, \cite{10}
it was shown that 
expansions for simplest iterated (double) Stratonovich stochastic integrals 
based on the systems of Haar and Rademacher--Walsh functions
are too complex and ineffective in practice.
In these works, a very brief comparison of the efficiency of 
application of Legendre 
polynomials and trigonometric functions in the framework of the method 
of generalized multiple Fourier series was also carried out.
The subject of this article is the development and refinement of 
the results obtained in \cite{3}-\cite{7}, \cite{9a}, \cite{10} 
in this direction.

\vspace{5mm}

\section{Milstein Approach}

\vspace{5mm}

Let $(\Omega,$ ${\rm F},$ ${\sf P})$ be a complete probability space. 
Let $\{{\rm F}_t, t\in[0,T]\}$ be a nondecreasing right-continous family 
of $\sigma$-algebras of ${\rm F},$
and let ${\bf f}_t$ be a standard $m$-dimensional Wiener 
stochastic process, which is
${\rm F}_t$-measurable for any $t\in[0, T].$ We assume that the components
${\bf f}_{t}^{(i)}$ $(i=1,\ldots,m)$ of this process are independent.
Consider the Brownian bridge process \cite{19}

\begin{equation}
\label{6.5.1}
{\bf f}_t-\frac{t}{\Delta}{\bf f}_{\Delta},\ \ \ t\in [0, \Delta].
\end{equation}

\vspace{4mm}

The componentwise expansion of the stochastic process (\ref{6.5.1}) 
into converging in the mean-square sense trigonometric Fourier series 
(version of the so-called Karhunen--Loeve expansion) has the following form
\cite{19}

\begin{equation}
\label{6.5.2}
{\bf f}_t^{(i)}-\frac{t}{\Delta}{\bf f}_{\Delta}^{(i)}=
\frac{1}{2}a_{i,0}+\sum_{r=1}^{\infty}\left(
a_{i,r}{\rm cos}\frac{2\pi rt}{\Delta} + b_{i,r}{\rm sin}
\frac{2\pi rt}{\Delta}\right),
\end{equation}

\vspace{4mm}
\noindent
where

\vspace{-1mm}
$$
a_{i,r}=\frac{2}{\Delta} \int\limits_0^{\Delta}
\left({\bf f}_s^{(i)}-\frac{s}{\Delta}{\bf f}_{\Delta}^{(i)}\right)
{\rm cos}\frac{2\pi rs}{\Delta}ds,\ \ \
b_{i,r}=\frac{2}{\Delta} \int\limits_0^{\Delta}
\left({\bf f}_s^{(i)}-\frac{s}{\Delta}{\bf f}_{\Delta}^{(i)}\right)
{\rm sin}\frac{2\pi rs}{\Delta}ds,
$$

\vspace{5mm}
\noindent
where $r=0, 1,\ldots;$ $i=1,\ldots,m.$ 

It is easy to demonstrate \cite{19} that the random variables
$a_{i,r}, b_{i,r}$ 
are Gaussian ones and they satisfy the following relations

$$
{\sf M}\left\{a_{i,r}b_{i,r}\right\}=
{\sf M}\left\{a_{i,r}b_{i,k}\right\}=0,\ \ \
{\sf M}\left\{a_{i,r}a_{i,k}\right\}=
{\sf M}\left\{b_{i,r}b_{i,k}\right\}=0,
$$

\vspace{2mm}
$$
{\sf M}\left\{a_{i_1,r}a_{i_2,r}\right\}=
{\sf M}\left\{b_{i_1,r}b_{i_2,r}\right\}=0,\ \ \
{\sf M}\left\{a_{i,r}^2\right\}=
{\sf M}\left\{b_{i,r}^2\right\}=\frac{\Delta}{2\pi^2 r^2},
$$

\vspace{5mm}
\noindent
where $i, i_1, i_2=1,\ldots,m;$ $r\ne k;$ $i_1\ne i_2;$
${\sf M}$ denotes a mathematical expectation.

According to (\ref{6.5.2}), we have

\begin{equation}
\label{6.5.7}
{\bf f}_t^{(i)}={\bf f}_{\Delta}^{(i)}\frac{t}{\Delta}+
\frac{1}{2}a_{i,0}+
\sum_{r=1}^{\infty}\left(
a_{i,r}{\rm cos}\frac{2\pi rt}{\Delta}+b_{i,r}{\rm sin}
\frac{2\pi rt}{\Delta}\right),
\end{equation}

\vspace{4mm}
\noindent
where the series
converges in the mean-square sense.

Denote
\begin{equation}
\label{ito1}
J_{(\lambda_1\ldots \lambda_k)T,t}^{(i_1\ldots i_k)}=
\int\limits_t^T \ldots \int\limits_t^{t_{2}}d{\bf w}_{t_1}^{(i_1)}\ldots
d{\bf w}_{t_k}^{(i_k)},
\end{equation}

\begin{equation}
\label{str1}
J_{(\lambda_1\ldots \lambda_k)T,t}^{*(i_1\ldots i_k)}=
\int\limits_t^{*T}\ldots \int\limits_t^{*t_{2}}
d{\bf w}_{t_1}^{(i_1)}\ldots
d{\bf w}_{t_k}^{(i_k)},
\end{equation}

\begin{equation}
\label{ito}
J[\psi^{(k)}]_{T,t}=\int\limits_t^T\psi_k(t_k) \ldots \int\limits_t^{t_{2}}
\psi_1(t_1) d{\bf w}_{t_1}^{(i_1)}\ldots
d{\bf w}_{t_k}^{(i_k)},
\end{equation}

\begin{equation}
\label{str}
J^{*}[\psi^{(k)}]_{T,t}=
\int\limits_t^{*T}\psi_k(t_k) \ldots \int\limits_t^{*t_{2}}
\psi_1(t_1) d{\bf w}_{t_1}^{(i_1)}\ldots
d{\bf w}_{t_k}^{(i_k)},
\end{equation}

\vspace{3mm}
\noindent
where every $\psi_l(\tau)$ $(l=1,\ldots,k)$ is a 
non-random function 
on $[t,T];$ ${\bf w}_{\tau}^{(i)}={\bf f}_{\tau}^{(i)}$
for $i=1,\ldots,m$ and
${\bf w}_{\tau}^{(0)}=\tau;$\
$i_1,\ldots,i_k = 0, 1,\ldots,m;$ $\lambda_l=0$ for $i_l=0$ and
$\lambda_l=1$ for $i_l=1,\ldots,m$ $(l=1,\ldots,k);$ 

\vspace{-1mm}
$$
\int\limits\ \hbox{and}\ \int\limits^{*}
$$ 

\vspace{4mm}
\noindent
denote Ito and 
Stratonovich stochastic integrals,
respectively.
In this paper we use the definition of the Stratonovich 
stochastic integral from \cite{20}, \cite{21}.

In \cite{19} Milstein G.N. obtained the following expansion of 
$J_{(11)T,t}^{(i_1i_2)}$ using the expansion (\ref{6.5.7})

\begin{equation}
\label{43}
J_{(11)T,t}^{(i_1 i_2)}=\frac{1}{2}(T-t)\Biggl(
\zeta_{0}^{(i_1)}\zeta_{0}^{(i_2)}
+\frac{1}{\pi}
\sum_{r=1}^{\infty}\frac{1}{r}\left(
\zeta_{2r}^{(i_1)}\zeta_{2r-1}^{(i_2)}-
\zeta_{2r-1}^{(i_1)}\zeta_{2r}^{(i_2)}
\right.\biggr.
+\left.\sqrt{2}\left(\zeta_{2r-1}^{(i_1)}\zeta_{0}^{(i_2)}-
\zeta_{0}^{(i_1)}\zeta_{2r-1}^{(i_2)}\right)\right)
\Biggr),
\end{equation}

\vspace{5mm}
\noindent
where the series converges in the mean-square sense;
$i_1\ne i_2;$ $i_1, i_2=1,\ldots,m;$

$$
\zeta_{j}^{(i)}=
\int\limits_t^T \phi_{j}(s) d{\bf f}_s^{(i)}
$$ 

\vspace{3mm}
\noindent
are independent standard Gaussian random variables
for various
$i$ or $j;$ 

\vspace{1mm}
\begin{equation}
\label{666.6}
\phi_j(s)=\frac{1}{\sqrt{T-t}}
\begin{cases}
1\ &\hbox{for}\ j=0\cr\cr
\sqrt{2}{\rm sin}(2\pi r(s-t)/(T-t))\ &\hbox{for}\ j=2r-1\cr\cr
\sqrt{2}{\rm cos}(2\pi r(s-t)/(T-t))\ &\hbox{for}\ j=2r
\end{cases},
\end{equation}

\vspace{5mm}
\noindent
where $r=1, 2,\ldots$
Moreover, \cite{19}

\vspace{-1mm}
\begin{equation}
\label{41}
J_{(1)T,t}^{(i_1)}=\sqrt{T-t}\zeta_0^{(i_1)},
\end{equation}

\vspace{3mm}
\noindent
where $i_1=1,\ldots,m.$

In principle for implementing the strong numerical method with the order $1.0$
of accuracy 
(Milstein method \cite{19}) for Ito stochastic differential equations 
it is sufficient
to take the following approximations

\vspace{-1mm}
\begin{equation}
\label{411}
J_{(1)T,t}^{(i_1)q}\stackrel{\sf def}{=}J_{(1)T,t}^{(i_1)}=\sqrt{T-t}\zeta_0^{(i_1)},
\end{equation}

\vspace{2mm}
\begin{equation}
\label{431}
J_{(11)T,t}^{(i_1 i_2)q}=\frac{1}{2}(T-t)\Biggl(
\zeta_{0}^{(i_1)}\zeta_{0}^{(i_2)}
+\frac{1}{\pi}
\sum_{r=1}^{q}\frac{1}{r}\left(
\zeta_{2r}^{(i_1)}\zeta_{2r-1}^{(i_2)}-
\zeta_{2r-1}^{(i_1)}\zeta_{2r}^{(i_2)}
\right.\biggr.
+\left.\sqrt{2}\left(\zeta_{2r-1}^{(i_1)}\zeta_{0}^{(i_2)}-
\zeta_{0}^{(i_1)}\zeta_{2r-1}^{(i_2)}\right)\right)
\biggl.
\Biggr),
\end{equation}

\vspace{5mm}
\noindent
where $i_1\ne i_2;$ $i_1, i_2=1,\ldots,m.$

It is not difficult to show that

\begin{equation}
\label{801}
{\sf M}\left\{\left(J_{(11)T,t}^{(i_1 i_2)}-
J_{(11)T,t}^{(i_1 i_2)q}
\right)^2\right\}
=\frac{3(T-t)^{2}}{2\pi^2}\Biggl(\frac{\pi^2}{6}-
\sum_{r=1}^q \frac{1}{r^2}\Biggr).
\end{equation}

\vspace{4mm}

However, this approach has an obvious drawback. Indeed, we have too complex 
formulas for the stochastic integrals with Gaussian distribution

\vspace{1mm}
\begin{equation}
\label{100}
J_{(01)T,t}^{(0 i_1)}=
\frac{{(T-t)}^{3/2}}{2}
\biggl(\zeta_0^{(i_1)}-\frac{\sqrt{2}}{\pi}\sum_{r=1}^{\infty}
\frac{1}{r}\zeta_{2r-1}^{(i_1)}\biggr),
\end{equation}

\vspace{2mm}
$$
J_{(001)T,t}^{(00 i_1)}=
(T-t)^{5/2}\Biggl(
\frac{1}{6}\zeta_0^{(i_1)}+\frac{1}{2\sqrt{2}\pi^2}
\sum_{r=1}^{\infty}\frac{1}{r^2}\zeta_{2r}^{(i_1)}-
\frac{1}{2\sqrt{2}\pi}\sum_{r=1}^{\infty}
\frac{1}{r}\zeta_{2r-1}^{(i_1)}\Biggr),
$$

\vspace{4mm}
$$
J_{(01)T,t}^{(0 i_1)q}=
\frac{{(T-t)}^{3/2}}{2}
\biggl(\zeta_0^{(i_1)}-\frac{\sqrt{2}}{\pi}\sum_{r=1}^{q}
\frac{1}{r}\zeta_{2r-1}^{(i_1)}\biggr),
$$

\vspace{2mm}
$$
J_{(001)T,t}^{(00 i_1)q}=
(T-t)^{5/2}\Biggl(
\frac{1}{6}\zeta_0^{(i_1)}+\frac{1}{2\sqrt{2}\pi^2}
\sum_{r=1}^{q}\frac{1}{r^2}\zeta_{2r}^{(i_1)}
-
\frac{1}{2\sqrt{2}\pi}\sum_{r=1}^q
\frac{1}{r}\zeta_{2r-1}^{(i_1)}\Biggr),
$$

\vspace{6mm}
\noindent
where the sense of notations from (\ref{431}) is hold.

In \cite{19} Milstein G.N. proposed the following mean-square 
approximations on the base of the expansions
(\ref{43}), (\ref{100})

\begin{equation}
\label{444}
J_{(01)T,t}^{(0 i_1)q}=\frac{{(T-t)}^{3/2}}{2}
\biggl(\zeta_0^{(i_1)}-\frac{\sqrt{2}}{\pi}\biggl(\sum_{r=1}^{q}
\frac{1}{r}
\zeta_{2r-1}^{(i_1)}+\sqrt{\alpha_q}\xi_q^{(i_1)}\biggr)
\biggr),
\end{equation}

\vspace{2mm}
$$
J_{(11)T,t}^{(i_1 i_2)q}=\frac{1}{2}(T-t)\Biggl(
\zeta_{0}^{(i_1)}\zeta_{0}^{(i_2)}
+\frac{1}{\pi}
\sum_{r=1}^{q}\frac{1}{r}\left(
\zeta_{2r}^{(i_1)}\zeta_{2r-1}^{(i_2)}-
\zeta_{2r-1}^{(i_1)}\zeta_{2r}^{(i_2)}+
\right.\Biggr.
$$

\begin{equation}
\label{555}
+\Biggl.\left.\sqrt{2}\left(\zeta_{2r-1}^{(i_1)}\zeta_{0}^{(i_2)}-
\zeta_{0}^{(i_1)}\zeta_{2r-1}^{(i_2)}\right)\right)
+\frac{\sqrt{2}}{\pi}\sqrt{\alpha_q}\left(
\xi_q^{(i_1)}\zeta_0^{(i_2)}-\zeta_0^{(i_1)}\xi_q^{(i_2)}\right)\Biggr),
\end{equation}

\vspace{5mm}
\noindent
where $i_1\ne i_2$ in (\ref{555}), and

\begin{equation}
\label{333}
\xi_q^{(i)}=\frac{1}{\sqrt{\alpha_q}}\sum_{r=q+1}^{\infty}
\frac{1}{r}~\zeta_{2r-1}^{(i)},\ \ \
\alpha_q=\frac{\pi^2}{6}-\sum_{r=1}^q\frac{1}{r^2},
\end{equation}

\vspace{3mm}
\noindent
where
$\zeta_0^{(i)},$ $\zeta_{2r}^{(i)},$
$\zeta_{2r-1}^{(i)},$ $\xi_q^{(i)};$ $r=1,\ldots,q;$
$i=1,\ldots,m$
are independent standard Gaussian random variables.

Obviously, for the
approximations (\ref{444}) and (\ref{555}) we obtain \cite{19}

\vspace{1mm}
$$
{\sf M}\left\{\left(J_{(01)T,t}^{(0 i_1)}-
J_{(01)T,t}^{(0 i_1)q}
\right)^2\right\}=0,
$$

\vspace{2mm}
\begin{equation}
\label{8010}
{\sf M}\left\{\left(J_{(11)T,t}^{(i_1 i_2)}-
J_{(11)T,t}^{(i_1 i_2)q}
\right)^2\right\}
=\frac{(T-t)^{2}}{2\pi^2}\Biggl(\frac{\pi^2}{6}-
\sum_{r=1}^q \frac{1}{r^2}\Biggr).
\end{equation}

\vspace{5mm}

This idea has been developed in \cite{20}-\cite{22}.
For example, the approximation $J_{(001)T,t}^{(0 0 i_1)q},$ 
which corresponds to (\ref{444}), (\ref{555})
has the form \cite{20}-\cite{22}

\vspace{1mm}
$$
J_{(001)T,t}^{(00 i_1)q}=
(T-t)^{5/2}\Biggl(
\frac{1}{6}\zeta_0^{(i_1)}+\frac{1}{2\sqrt{2}\pi^2}\Biggl(
\sum_{r=1}^{q}\frac{1}{r^2}\zeta_{2r}^{(i_1)}+
\sqrt{\beta_q}\mu_q^{(i_1)}\Biggr)-\Biggl.
$$

\vspace{1mm}
\begin{equation}
\label{1970}
\Biggl.
-\frac{1}{2\sqrt{2}\pi}\Biggl(\sum_{r=1}^q
\frac{1}{r}\zeta_{2r-1}^{(i_1)}+\sqrt{\alpha_q}\xi_q^{(i_1)}\Biggr)\Biggr),
\end{equation}

\vspace{4mm}

$$
{\sf M}\left\{\left(J_{(001)T,t}^{(0 0 i_1)}-
J_{(0 0 1)T,t}^{(0 0 i_1)q}
\right)^2\right\}=0,
$$

\vspace{5mm}
\noindent
where
$\xi_q^{(i)}$ and $\alpha_q$ have the form (\ref{333}) and 

\vspace{1mm}
$$
\mu_q^{(i)}=\frac{1}{\sqrt{\beta_q}}\sum_{r=q+1}^{\infty}
\frac{1}{r^2}~\zeta_{2r}^{(i)},\ \ \
\beta_q=\frac{\pi^4}{90}-\sum_{r=1}^q\frac{1}{r^4},
$$

\vspace{4mm}
\noindent
$\phi_j(s)$\ is defined by (\ref{666.6}); 
$\zeta_0^{(i)},$ $\zeta_{2r}^{(i)},$
$\zeta_{2r-1}^{(i)},$ $\xi_q^{(i)},$ $\mu_q^{(i)}$ ($r=1,\ldots,q;$
$i=1,\ldots,m$) are independent
standard Gaussian random variables.

Nevetheless, the expansions (\ref{444}), (\ref{1970}) are too complex for
approximation of two Gaussian random variables
$J_{(01)T,t}^{(0 i_1)},$ $J_{(001)T,t}^{(0 0 i_1)}.$

Further, we will see that introduction of random 
variables $\xi_q^{(i)}$ and 
$\mu_q^{(i)}$ will sharply 
complicate the approximation of the stochastic 
integral $J_{(111)T,t}^{(i_1 i_2 i_3)};$
$i_1,i_2,i_3=1,\ldots,m$ within the framework of the Milshtein approach.
This is due to the fact that
the number $q$ is fixed for all stochastic integrals 
included into the considered collection. However, it is clear that due 
to the smallness of $T-t$, the number $q$ for $J_{(111)T,t}^{(i_1 i_2 i_3)}$
could be taken significantly 
less than the number $q$ in the formula (\ref{555}). 
This feature is also valid for the formulas (\ref{444}), (\ref{1970}).

On the other hand, the following very simple formulas are well known

\begin{equation}
\label{4}
J_{(1)T,t}^{(i_1)}=\sqrt{T-t}\zeta_0^{(i_1)},
\end{equation}

\vspace{1mm}
\begin{equation}
\label{5}
J_{(01)T,t}^{(0 i_1)}=\frac{(T-t)^{3/2}}{2}\biggl(\zeta_0^{(i_1)}+
\frac{1}{\sqrt{3}}\zeta_1^{(i_1)}\biggr),
\end{equation}

\vspace{2mm}
\begin{equation}
\label{6}
J_{(001)T,t}^{(0 0 i_1)}=\frac{(T-t)^{5/2}}{6}\biggl(\zeta_0^{(i_1)}+
\frac{\sqrt{3}}{2}\zeta_1^{(i_1)}+
\frac{1}{2\sqrt{5}}\zeta_2^{(i_1)}\biggr),
\end{equation}

\vspace{5mm}
\noindent
where
$\zeta_0^{(i)},$ $\zeta_1^{(i)},$ $\zeta_2^{(i)};$ $i=1,\ldots,m$ 
are indepentent standard Gaussian random variables.

Looking ahead, we note that the formulas (\ref{4})-(\ref{6}) 
are part of the method that 
will be discussed in the next section (see Theorems 1, 2 below).

To obtain the Milstein expansion for (\ref{str}) 
the truncated 
expansions (\ref{6.5.7}) of components of the 
Wiener process ${\bf f}_s$ must be
iteratively substituted in the single integrals, and the integrals
must be calculated, starting from the innermost integral.
This is a complicated procedure that obviously does not lead to a general
expansion of (\ref{str}) valid for an arbitrary multiplicity $k.$
For this reason, only expansions of simplest single, double, and triple
integrals (\ref{str}) were obtained \cite{19}-\cite{Zapad-9}.

At that, in \cite{19}, \cite{23} the case 
$\psi_1(s), \psi_2(s)\equiv 1$ and
$i_1, i_2=0, 1,\ldots,m$ is considered. In 
\cite{20}-\cite{22}, \cite{Zapad-9} the attempt to consider the case 
$\psi_1(s), \psi_2(s), \psi_3(s)\equiv 1$ and 
$i_1, i_2, i_3=0, 1,\ldots,m$ is implemented.

Note that generally speaking
the mean-square convergence of the approximation 

\vspace{-1mm}
$$
J_{(111)T,t}^{*(i_1 i_2 i_3)q}
$$ 

\vspace{2mm}
\noindent
(obtained by the Milstein approach)
to the appropriate iterated Stratonovich 
stochastic integral

\vspace{-1mm}
$$
J_{(111)T,t}^{*(i_1 i_2 i_3)}
$$ 

\vspace{2mm}
\noindent
must be proved separately
due to iterated application of passing to the 
limit in the Milstein approach \cite{19}.
However, 
in \cite{20} (pp.~438-439),
\cite{21}
(Sect.~5.8, pp.~202--204), \cite{22} (pp.~82-84),
\cite{Zapad-9} (pp.~263-264) 
the authors use the Wong--Zakai approximation
\cite{W-Z-1}-\cite{Watanabe}
(without rigorous proof)
within the frames of the mentioned approach
based on the Karhunen--Loeve expansion of the Brownian bridge
process \cite{19} (see discussion in Sect.~11 for details).

\vspace{5mm}

\section{Method of Generalized Multiple Fourier Series}

\vspace{5mm}

Let us consider an another approach to the expansion of iterated
Ito stochastic integrals \cite{3}-\cite{new-art-1xxys} (method of generalized
multiple Fourier series).

The idea of this method is as follows: the iterated Ito stochastic 
integral (\ref{ito}) 
of multiplicity $k$ is represented as the multiple stochastic 
integral from the certain non-random discontinuous function of $k$ variables 
defined on the hypercube $[t, T]^k$, where $[t, T]$ is the interval of 
integration of the iterated Ito stochastic integral. Then, the indicated 
non-random function is expanded in the hypercube $[t, T]^k$
into the generalized 
multiple Fourier series converging 
in the sense of norm in Hilbert space
$L_2([t,T]^k)$. After a number of nontrivial transformations we come 
(see Theorems 1, 2 below) to the 
mean-square convergening expansion of the iterated Ito stochastic 
integral (\ref{ito}) into the multiple 
series of products
of standard  Gaussian random 
variables. Coefficients of this 
series are coefficients of 
generalized multiple Fourier series for the mentioned non-random function 
of $k$ variables, which can be calculated using the explicit formula 
regardless of the multiplicity $k$ of the iterated Ito stochastic integral
(\ref{ito}).

Suppose that every $\psi_l(\tau)$ $(l=1,\ldots,k)$ is a 
non-random function from the space $L_2([t, T])$.
Define the following function on the hypercube $[t, T]^k$

\vspace{1mm}
\begin{equation}
\label{ppp}
K(t_1,\ldots,t_k)=
\begin{cases}
\psi_1(t_1)\ldots \psi_k(t_k)\ \ \hbox{for}\ \ t_1<\ldots<t_k\\
~\\
~\\
0\ \ \hbox{otherwise}
\end{cases},\ \ \ t_1,\ldots,t_k\in[t, T],\ \ \ k\ge 2,
\end{equation}

\vspace{6mm}
\noindent
and 
$K(t_1)\equiv\psi_1(t_1)$ for $t_1\in[t, T].$

Suppose that $\{\phi_j(x)\}_{j=0}^{\infty}$
is a complete orthonormal system of functions in the space
$L_2([t, T])$.

The function $K(t_1,\ldots,t_k)$ belongs to the space 
$L_2([t, T]^k).$
At this situation it is well known that the generalized 
multiple Fourier series 
of $K(t_1,\ldots,t_k)\in L_2([t, T]^k)$ is converging 
to $K(t_1,\ldots,t_k)$ in the hypercube $[t, T]^k$ in 
the mean-square sense, i.e.

\begin{equation}
\label{sos1z}
\hbox{\vtop{\offinterlineskip\halign{
\hfil#\hfil\cr
{\rm lim}\cr
$\stackrel{}{{}_{p_1,\ldots,p_k\to \infty}}$\cr
}} }\Biggl\Vert
K(t_1,\ldots,t_k)-
\sum_{j_1=0}^{p_1}\ldots \sum_{j_k=0}^{p_k}
C_{j_k\ldots j_1}\prod_{l=1}^{k} \phi_{j_l}(t_l)\Biggr
\Vert_{L_2([t,T]^k)}=0,
\end{equation}

\vspace{3mm}
\noindent
where

\vspace{-1mm}
\begin{equation}
\label{ppppa}
C_{j_k\ldots j_1}=\int\limits_{[t,T]^k}
K(t_1,\ldots,t_k)\prod_{l=1}^{k}\phi_{j_l}(t_l)dt_1\ldots dt_k
\end{equation}

\vspace{6mm}
\noindent
is the Fourier coefficient,

$$
\left\Vert f\right\Vert_{L_2([t,T]^k)}=\left(\int\limits_{[t,T]^k}
f^2(t_1,\ldots,t_k)dt_1\ldots dt_k\right)^{1/2}.
$$

\vspace{6mm}

Consider the partition $\{\tau_j\}_{j=0}^N$ of $[t,T]$ such that

\vspace{2mm}
\begin{equation}
\label{1111}
t=\tau_0<\ldots <\tau_N=T,\ \ \
\Delta_N=
\hbox{\vtop{\offinterlineskip\halign{
\hfil#\hfil\cr
{\rm max}\cr
$\stackrel{}{{}_{0\le j\le N-1}}$\cr
}} }\Delta\tau_j\to 0\ \ \hbox{if}\ \ N\to \infty,\ \ \
\Delta\tau_j=\tau_{j+1}-\tau_j.
\end{equation}

\vspace{5mm}

{\bf Theorem 1} \cite{3} (2006), \cite{3a}-\cite{17}. {\it Suppose that
every $\psi_l(\tau)$ $(l=1,\ldots, k)$ is a continuous non-random 
function on the interval
$[t, T]$ and
$\{\phi_j(x)\}_{j=0}^{\infty}$ is a complete orthonormal system  
of continuous functions in the space $L_2([t,T]).$ Then

$$
J[\psi^{(k)}]_{T,t}\  =\ 
\hbox{\vtop{\offinterlineskip\halign{
\hfil#\hfil\cr
{\rm l.i.m.}\cr
$\stackrel{}{{}_{p_1,\ldots,p_k\to \infty}}$\cr
}} }\sum_{j_1=0}^{p_1}\ldots\sum_{j_k=0}^{p_k}
C_{j_k\ldots j_1}\Biggl(
\prod_{l=1}^k\zeta_{j_l}^{(i_l)}\ -
\Biggr.
$$

\vspace{3mm}
\begin{equation}
\label{tyyy}
-\ \Biggl.
\hbox{\vtop{\offinterlineskip\halign{
\hfil#\hfil\cr
{\rm l.i.m.}\cr
$\stackrel{}{{}_{N\to \infty}}$\cr
}} }\sum_{(l_1,\ldots,l_k)\in {\rm G}_k}
\phi_{j_{1}}(\tau_{l_1})
\Delta{\bf w}_{\tau_{l_1}}^{(i_1)}\ldots
\phi_{j_{k}}(\tau_{l_k})
\Delta{\bf w}_{\tau_{l_k}}^{(i_k)}\Biggr),
\end{equation}

\vspace{7mm}
\noindent
where $J[\psi^{(k)}]_{T,t}$ is defined by {\rm (\ref{ito}),}

\vspace{2mm}
$$
{\rm G}_k={\rm H}_k\backslash{\rm L}_k,\ \ \
{\rm H}_k=\{(l_1,\ldots,l_k):\ l_1,\ldots,l_k=0,\ 1,\ldots,N-1\},
$$

\vspace{2mm}
$$
{\rm L}_k=\{(l_1,\ldots,l_k):\ l_1,\ldots,l_k=0,\ 1,\ldots,N-1;\
l_g\ne l_r\ (g\ne r);\ g, r=1,\ldots,k\},
$$

\vspace{7mm}
\noindent
${\rm l.i.m.}$ is a limit in the mean-square sense$,$
$i_1,\ldots,i_k=0,1,\ldots,m,$

\vspace{1mm}
\begin{equation}
\label{rr23}
\zeta_{j}^{(i)}=
\int\limits_t^T \phi_{j}(s) d{\bf w}_s^{(i)}
\end{equation} 

\vspace{4mm}
\noindent
are independent standard Gaussian random variables
for various
$i$ or $j$ {\rm(}if $i\ne 0${\rm),}
$C_{j_k\ldots j_1}$ is the Fourier coefficient {\rm(\ref{ppppa}),}
$\Delta{\bf w}_{\tau_{j}}^{(i)}=
{\bf w}_{\tau_{j+1}}^{(i)}-{\bf w}_{\tau_{j}}^{(i)}$
$(i=0, 1,\ldots,m),$
$\left\{\tau_{j}\right\}_{j=0}^{N}$ is a partition of
the interval $[t, T],$ which satisfies the condition {\rm (\ref{1111})}.
}

\vspace{2mm}

Note that the continuity condition of $\phi_j(x)$ can be 
weakened (see \cite{3}-\cite{11}). Moreover, 
Theorem 1 can be generalized to the case
of an arbitrary complete orthonormal systems  
of functions in the space $L_2([t,T])$ (see Theorem 2 below).

In order to evaluate the significance of Theorem 1 for practice we will
demonstrate its transformed particular cases for 
$k=1,\ldots,6$ \cite{3}-\cite{17}

\vspace{2mm}
\begin{equation}
\label{a1}
J[\psi^{(1)}]_{T,t}
=\hbox{\vtop{\offinterlineskip\halign{
\hfil#\hfil\cr
{\rm l.i.m.}\cr
$\stackrel{}{{}_{p_1\to \infty}}$\cr
}} }\sum_{j_1=0}^{p_1}
C_{j_1}\zeta_{j_1}^{(i_1)},
\end{equation}

\vspace{5mm}

\begin{equation}
\label{leto5001}
J[\psi^{(2)}]_{T,t}
=\hbox{\vtop{\offinterlineskip\halign{
\hfil#\hfil\cr
{\rm l.i.m.}\cr
$\stackrel{}{{}_{p_1,p_2\to \infty}}$\cr
}} }\sum_{j_1=0}^{p_1}\sum_{j_2=0}^{p_2}
C_{j_2j_1}\Biggl(\zeta_{j_1}^{(i_1)}\zeta_{j_2}^{(i_2)}
-{\bf 1}_{\{i_1=i_2\ne 0\}}
{\bf 1}_{\{j_1=j_2\}}\Biggr),
\end{equation}

\vspace{7mm}
$$
J[\psi^{(3)}]_{T,t}=
\hbox{\vtop{\offinterlineskip\halign{
\hfil#\hfil\cr
{\rm l.i.m.}\cr
$\stackrel{}{{}_{p_1,\ldots,p_3\to \infty}}$\cr
}} }\sum_{j_1=0}^{p_1}\sum_{j_2=0}^{p_2}\sum_{j_3=0}^{p_3}
C_{j_3j_2j_1}\Biggl(
\zeta_{j_1}^{(i_1)}\zeta_{j_2}^{(i_2)}\zeta_{j_3}^{(i_3)}
-\Biggr.
$$

\vspace{2mm}
\begin{equation}
\label{leto5002}
\Biggl.-{\bf 1}_{\{i_1=i_2\ne 0\}}
{\bf 1}_{\{j_1=j_2\}}
\zeta_{j_3}^{(i_3)}
-{\bf 1}_{\{i_2=i_3\ne 0\}}
{\bf 1}_{\{j_2=j_3\}}
\zeta_{j_1}^{(i_1)}-
{\bf 1}_{\{i_1=i_3\ne 0\}}
{\bf 1}_{\{j_1=j_3\}}
\zeta_{j_2}^{(i_2)}\Biggr),
\end{equation}

\vspace{7mm}
$$
J[\psi^{(4)}]_{T,t}
=
\hbox{\vtop{\offinterlineskip\halign{
\hfil#\hfil\cr
{\rm l.i.m.}\cr
$\stackrel{}{{}_{p_1,\ldots,p_4\to \infty}}$\cr
}} }\sum_{j_1=0}^{p_1}\ldots\sum_{j_4=0}^{p_4}
C_{j_4\ldots j_1}\Biggl(
\prod_{l=1}^4\zeta_{j_l}^{(i_l)}
\Biggr.
-
$$
$$
-
{\bf 1}_{\{i_1=i_2\ne 0\}}
{\bf 1}_{\{j_1=j_2\}}
\zeta_{j_3}^{(i_3)}
\zeta_{j_4}^{(i_4)}
-
{\bf 1}_{\{i_1=i_3\ne 0\}}
{\bf 1}_{\{j_1=j_3\}}
\zeta_{j_2}^{(i_2)}
\zeta_{j_4}^{(i_4)}-
$$
$$
-
{\bf 1}_{\{i_1=i_4\ne 0\}}
{\bf 1}_{\{j_1=j_4\}}
\zeta_{j_2}^{(i_2)}
\zeta_{j_3}^{(i_3)}
-
{\bf 1}_{\{i_2=i_3\ne 0\}}
{\bf 1}_{\{j_2=j_3\}}
\zeta_{j_1}^{(i_1)}
\zeta_{j_4}^{(i_4)}-
$$
$$
-
{\bf 1}_{\{i_2=i_4\ne 0\}}
{\bf 1}_{\{j_2=j_4\}}
\zeta_{j_1}^{(i_1)}
\zeta_{j_3}^{(i_3)}
-
{\bf 1}_{\{i_3=i_4\ne 0\}}
{\bf 1}_{\{j_3=j_4\}}
\zeta_{j_1}^{(i_1)}
\zeta_{j_2}^{(i_2)}+
$$
$$
+
{\bf 1}_{\{i_1=i_2\ne 0\}}
{\bf 1}_{\{j_1=j_2\}}
{\bf 1}_{\{i_3=i_4\ne 0\}}
{\bf 1}_{\{j_3=j_4\}}
+
$$
$$
+
{\bf 1}_{\{i_1=i_3\ne 0\}}
{\bf 1}_{\{j_1=j_3\}}
{\bf 1}_{\{i_2=i_4\ne 0\}}
{\bf 1}_{\{j_2=j_4\}}+
$$
\begin{equation}
\label{leto5003}
+\Biggl.
{\bf 1}_{\{i_1=i_4\ne 0\}}
{\bf 1}_{\{j_1=j_4\}}
{\bf 1}_{\{i_2=i_3\ne 0\}}
{\bf 1}_{\{j_2=j_3\}}\Biggr),
\end{equation}

\vspace{9mm}

$$
J[\psi^{(5)}]_{T,t}
=\hbox{\vtop{\offinterlineskip\halign{
\hfil#\hfil\cr
{\rm l.i.m.}\cr
$\stackrel{}{{}_{p_1,\ldots,p_5\to \infty}}$\cr
}} }\sum_{j_1=0}^{p_1}\ldots\sum_{j_5=0}^{p_5}
C_{j_5\ldots j_1}\Biggl(
\prod_{l=1}^5\zeta_{j_l}^{(i_l)}
-\Biggr.
$$
$$
-
{\bf 1}_{\{i_1=i_2\ne 0\}}
{\bf 1}_{\{j_1=j_2\}}
\zeta_{j_3}^{(i_3)}
\zeta_{j_4}^{(i_4)}
\zeta_{j_5}^{(i_5)}-
{\bf 1}_{\{i_1=i_3\ne 0\}}
{\bf 1}_{\{j_1=j_3\}}
\zeta_{j_2}^{(i_2)}
\zeta_{j_4}^{(i_4)}
\zeta_{j_5}^{(i_5)}-
$$
$$
-
{\bf 1}_{\{i_1=i_4\ne 0\}}
{\bf 1}_{\{j_1=j_4\}}
\zeta_{j_2}^{(i_2)}
\zeta_{j_3}^{(i_3)}
\zeta_{j_5}^{(i_5)}-
{\bf 1}_{\{i_1=i_5\ne 0\}}
{\bf 1}_{\{j_1=j_5\}}
\zeta_{j_2}^{(i_2)}
\zeta_{j_3}^{(i_3)}
\zeta_{j_4}^{(i_4)}-
$$
$$
-
{\bf 1}_{\{i_2=i_3\ne 0\}}
{\bf 1}_{\{j_2=j_3\}}
\zeta_{j_1}^{(i_1)}
\zeta_{j_4}^{(i_4)}
\zeta_{j_5}^{(i_5)}-
{\bf 1}_{\{i_2=i_4\ne 0\}}
{\bf 1}_{\{j_2=j_4\}}
\zeta_{j_1}^{(i_1)}
\zeta_{j_3}^{(i_3)}
\zeta_{j_5}^{(i_5)}-
$$
$$
-
{\bf 1}_{\{i_2=i_5\ne 0\}}
{\bf 1}_{\{j_2=j_5\}}
\zeta_{j_1}^{(i_1)}
\zeta_{j_3}^{(i_3)}
\zeta_{j_4}^{(i_4)}
-{\bf 1}_{\{i_3=i_4\ne 0\}}
{\bf 1}_{\{j_3=j_4\}}
\zeta_{j_1}^{(i_1)}
\zeta_{j_2}^{(i_2)}
\zeta_{j_5}^{(i_5)}-
$$
$$
-
{\bf 1}_{\{i_3=i_5\ne 0\}}
{\bf 1}_{\{j_3=j_5\}}
\zeta_{j_1}^{(i_1)}
\zeta_{j_2}^{(i_2)}
\zeta_{j_4}^{(i_4)}
-{\bf 1}_{\{i_4=i_5\ne 0\}}
{\bf 1}_{\{j_4=j_5\}}
\zeta_{j_1}^{(i_1)}
\zeta_{j_2}^{(i_2)}
\zeta_{j_3}^{(i_3)}+
$$
$$
+
{\bf 1}_{\{i_1=i_2\ne 0\}}
{\bf 1}_{\{j_1=j_2\}}
{\bf 1}_{\{i_3=i_4\ne 0\}}
{\bf 1}_{\{j_3=j_4\}}\zeta_{j_5}^{(i_5)}+
{\bf 1}_{\{i_1=i_2\ne 0\}}
{\bf 1}_{\{j_1=j_2\}}
{\bf 1}_{\{i_3=i_5\ne 0\}}
{\bf 1}_{\{j_3=j_5\}}\zeta_{j_4}^{(i_4)}+
$$
$$
+
{\bf 1}_{\{i_1=i_2\ne 0\}}
{\bf 1}_{\{j_1=j_2\}}
{\bf 1}_{\{i_4=i_5\ne 0\}}
{\bf 1}_{\{j_4=j_5\}}\zeta_{j_3}^{(i_3)}+
{\bf 1}_{\{i_1=i_3\ne 0\}}
{\bf 1}_{\{j_1=j_3\}}
{\bf 1}_{\{i_2=i_4\ne 0\}}
{\bf 1}_{\{j_2=j_4\}}\zeta_{j_5}^{(i_5)}+
$$
$$
+
{\bf 1}_{\{i_1=i_3\ne 0\}}
{\bf 1}_{\{j_1=j_3\}}
{\bf 1}_{\{i_2=i_5\ne 0\}}
{\bf 1}_{\{j_2=j_5\}}\zeta_{j_4}^{(i_4)}+
{\bf 1}_{\{i_1=i_3\ne 0\}}
{\bf 1}_{\{j_1=j_3\}}
{\bf 1}_{\{i_4=i_5\ne 0\}}
{\bf 1}_{\{j_4=j_5\}}\zeta_{j_2}^{(i_2)}+
$$
$$
+
{\bf 1}_{\{i_1=i_4\ne 0\}}
{\bf 1}_{\{j_1=j_4\}}
{\bf 1}_{\{i_2=i_3\ne 0\}}
{\bf 1}_{\{j_2=j_3\}}\zeta_{j_5}^{(i_5)}+
{\bf 1}_{\{i_1=i_4\ne 0\}}
{\bf 1}_{\{j_1=j_4\}}
{\bf 1}_{\{i_2=i_5\ne 0\}}
{\bf 1}_{\{j_2=j_5\}}\zeta_{j_3}^{(i_3)}+
$$
$$
+
{\bf 1}_{\{i_1=i_4\ne 0\}}
{\bf 1}_{\{j_1=j_4\}}
{\bf 1}_{\{i_3=i_5\ne 0\}}
{\bf 1}_{\{j_3=j_5\}}\zeta_{j_2}^{(i_2)}+
{\bf 1}_{\{i_1=i_5\ne 0\}}
{\bf 1}_{\{j_1=j_5\}}
{\bf 1}_{\{i_2=i_3\ne 0\}}
{\bf 1}_{\{j_2=j_3\}}\zeta_{j_4}^{(i_4)}+
$$
$$
+
{\bf 1}_{\{i_1=i_5\ne 0\}}
{\bf 1}_{\{j_1=j_5\}}
{\bf 1}_{\{i_2=i_4\ne 0\}}
{\bf 1}_{\{j_2=j_4\}}\zeta_{j_3}^{(i_3)}+
{\bf 1}_{\{i_1=i_5\ne 0\}}
{\bf 1}_{\{j_1=j_5\}}
{\bf 1}_{\{i_3=i_4\ne 0\}}
{\bf 1}_{\{j_3=j_4\}}\zeta_{j_2}^{(i_2)}+
$$
$$
+
{\bf 1}_{\{i_2=i_3\ne 0\}}
{\bf 1}_{\{j_2=j_3\}}
{\bf 1}_{\{i_4=i_5\ne 0\}}
{\bf 1}_{\{j_4=j_5\}}\zeta_{j_1}^{(i_1)}+
{\bf 1}_{\{i_2=i_4\ne 0\}}
{\bf 1}_{\{j_2=j_4\}}
{\bf 1}_{\{i_3=i_5\ne 0\}}
{\bf 1}_{\{j_3=j_5\}}\zeta_{j_1}^{(i_1)}+
$$
\begin{equation}
\label{a5}
+\Biggl.
{\bf 1}_{\{i_2=i_5\ne 0\}}
{\bf 1}_{\{j_2=j_5\}}
{\bf 1}_{\{i_3=i_4\ne 0\}}
{\bf 1}_{\{j_3=j_4\}}\zeta_{j_1}^{(i_1)}\Biggr),
\end{equation}

\vspace{9mm}

$$
J[\psi^{(6)}]_{T,t}
=\hbox{\vtop{\offinterlineskip\halign{
\hfil#\hfil\cr
{\rm l.i.m.}\cr
$\stackrel{}{{}_{p_1,\ldots,p_6\to \infty}}$\cr
}} }\sum_{j_1=0}^{p_1}\ldots\sum_{j_6=0}^{p_6}
C_{j_6\ldots j_1}\Biggl(
\prod_{l=1}^6
\zeta_{j_l}^{(i_l)}
-\Biggr.
$$
$$
-
{\bf 1}_{\{i_1=i_6\ne 0\}}
{\bf 1}_{\{j_1=j_6\}}
\zeta_{j_2}^{(i_2)}
\zeta_{j_3}^{(i_3)}
\zeta_{j_4}^{(i_4)}
\zeta_{j_5}^{(i_5)}-
{\bf 1}_{\{i_2=i_6\ne 0\}}
{\bf 1}_{\{j_2=j_6\}}
\zeta_{j_1}^{(i_1)}
\zeta_{j_3}^{(i_3)}
\zeta_{j_4}^{(i_4)}
\zeta_{j_5}^{(i_5)}-
$$
$$
-
{\bf 1}_{\{i_3=i_6\ne 0\}}
{\bf 1}_{\{j_3=j_6\}}
\zeta_{j_1}^{(i_1)}
\zeta_{j_2}^{(i_2)}
\zeta_{j_4}^{(i_4)}
\zeta_{j_5}^{(i_5)}-
{\bf 1}_{\{i_4=i_6\ne 0\}}
{\bf 1}_{\{j_4=j_6\}}
\zeta_{j_1}^{(i_1)}
\zeta_{j_2}^{(i_2)}
\zeta_{j_3}^{(i_3)}
\zeta_{j_5}^{(i_5)}-
$$
$$
-
{\bf 1}_{\{i_5=i_6\ne 0\}}
{\bf 1}_{\{j_5=j_6\}}
\zeta_{j_1}^{(i_1)}
\zeta_{j_2}^{(i_2)}
\zeta_{j_3}^{(i_3)}
\zeta_{j_4}^{(i_4)}-
{\bf 1}_{\{i_1=i_2\ne 0\}}
{\bf 1}_{\{j_1=j_2\}}
\zeta_{j_3}^{(i_3)}
\zeta_{j_4}^{(i_4)}
\zeta_{j_5}^{(i_5)}
\zeta_{j_6}^{(i_6)}-
$$
$$
-
{\bf 1}_{\{i_1=i_3\ne 0\}}
{\bf 1}_{\{j_1=j_3\}}
\zeta_{j_2}^{(i_2)}
\zeta_{j_4}^{(i_4)}
\zeta_{j_5}^{(i_5)}
\zeta_{j_6}^{(i_6)}-
{\bf 1}_{\{i_1=i_4\ne 0\}}
{\bf 1}_{\{j_1=j_4\}}
\zeta_{j_2}^{(i_2)}
\zeta_{j_3}^{(i_3)}
\zeta_{j_5}^{(i_5)}
\zeta_{j_6}^{(i_6)}-
$$
$$
-
{\bf 1}_{\{i_1=i_5\ne 0\}}
{\bf 1}_{\{j_1=j_5\}}
\zeta_{j_2}^{(i_2)}
\zeta_{j_3}^{(i_3)}
\zeta_{j_4}^{(i_4)}
\zeta_{j_6}^{(i_6)}-
{\bf 1}_{\{i_2=i_3\ne 0\}}
{\bf 1}_{\{j_2=j_3\}}
\zeta_{j_1}^{(i_1)}
\zeta_{j_4}^{(i_4)}
\zeta_{j_5}^{(i_5)}
\zeta_{j_6}^{(i_6)}-
$$
$$
-
{\bf 1}_{\{i_2=i_4\ne 0\}}
{\bf 1}_{\{j_2=j_4\}}
\zeta_{j_1}^{(i_1)}
\zeta_{j_3}^{(i_3)}
\zeta_{j_5}^{(i_5)}
\zeta_{j_6}^{(i_6)}-
{\bf 1}_{\{i_2=i_5\ne 0\}}
{\bf 1}_{\{j_2=j_5\}}
\zeta_{j_1}^{(i_1)}
\zeta_{j_3}^{(i_3)}
\zeta_{j_4}^{(i_4)}
\zeta_{j_6}^{(i_6)}-
$$
$$
-
{\bf 1}_{\{i_3=i_4\ne 0\}}
{\bf 1}_{\{j_3=j_4\}}
\zeta_{j_1}^{(i_1)}
\zeta_{j_2}^{(i_2)}
\zeta_{j_5}^{(i_5)}
\zeta_{j_6}^{(i_6)}-
{\bf 1}_{\{i_3=i_5\ne 0\}}
{\bf 1}_{\{j_3=j_5\}}
\zeta_{j_1}^{(i_1)}
\zeta_{j_2}^{(i_2)}
\zeta_{j_4}^{(i_4)}
\zeta_{j_6}^{(i_6)}-
$$
$$
-
{\bf 1}_{\{i_4=i_5\ne 0\}}
{\bf 1}_{\{j_4=j_5\}}
\zeta_{j_1}^{(i_1)}
\zeta_{j_2}^{(i_2)}
\zeta_{j_3}^{(i_3)}
\zeta_{j_6}^{(i_6)}+
$$
$$
+
{\bf 1}_{\{i_1=i_2\ne 0\}}
{\bf 1}_{\{j_1=j_2\}}
{\bf 1}_{\{i_3=i_4\ne 0\}}
{\bf 1}_{\{j_3=j_4\}}
\zeta_{j_5}^{(i_5)}
\zeta_{j_6}^{(i_6)}+
{\bf 1}_{\{i_1=i_2\ne 0\}}
{\bf 1}_{\{j_1=j_2\}}
{\bf 1}_{\{i_3=i_5\ne 0\}}
{\bf 1}_{\{j_3=j_5\}}
\zeta_{j_4}^{(i_4)}
\zeta_{j_6}^{(i_6)}+
$$
$$
+
{\bf 1}_{\{i_1=i_2\ne 0\}}
{\bf 1}_{\{j_1=j_2\}}
{\bf 1}_{\{i_4=i_5\ne 0\}}
{\bf 1}_{\{j_4=j_5\}}
\zeta_{j_3}^{(i_3)}
\zeta_{j_6}^{(i_6)}
+
{\bf 1}_{\{i_1=i_3\ne 0\}}
{\bf 1}_{\{j_1=j_3\}}
{\bf 1}_{\{i_2=i_4\ne 0\}}
{\bf 1}_{\{j_2=j_4\}}
\zeta_{j_5}^{(i_5)}
\zeta_{j_6}^{(i_6)}+
$$
$$
+
{\bf 1}_{\{i_1=i_3\ne 0\}}
{\bf 1}_{\{j_1=j_3\}}
{\bf 1}_{\{i_2=i_5\ne 0\}}
{\bf 1}_{\{j_2=j_5\}}
\zeta_{j_4}^{(i_4)}
\zeta_{j_6}^{(i_6)}
+{\bf 1}_{\{i_1=i_3\ne 0\}}
{\bf 1}_{\{j_1=j_3\}}
{\bf 1}_{\{i_4=i_5\ne 0\}}
{\bf 1}_{\{j_4=j_5\}}
\zeta_{j_2}^{(i_2)}
\zeta_{j_6}^{(i_6)}+
$$
$$
+
{\bf 1}_{\{i_1=i_4\ne 0\}}
{\bf 1}_{\{j_1=j_4\}}
{\bf 1}_{\{i_2=i_3\ne 0\}}
{\bf 1}_{\{j_2=j_3\}}
\zeta_{j_5}^{(i_5)}
\zeta_{j_6}^{(i_6)}
+
{\bf 1}_{\{i_1=i_4\ne 0\}}
{\bf 1}_{\{j_1=j_4\}}
{\bf 1}_{\{i_2=i_5\ne 0\}}
{\bf 1}_{\{j_2=j_5\}}
\zeta_{j_3}^{(i_3)}
\zeta_{j_6}^{(i_6)}+
$$
$$
+
{\bf 1}_{\{i_1=i_4\ne 0\}}
{\bf 1}_{\{j_1=j_4\}}
{\bf 1}_{\{i_3=i_5\ne 0\}}
{\bf 1}_{\{j_3=j_5\}}
\zeta_{j_2}^{(i_2)}
\zeta_{j_6}^{(i_6)}
+
{\bf 1}_{\{i_1=i_5\ne 0\}}
{\bf 1}_{\{j_1=j_5\}}
{\bf 1}_{\{i_2=i_3\ne 0\}}
{\bf 1}_{\{j_2=j_3\}}
\zeta_{j_4}^{(i_4)}
\zeta_{j_6}^{(i_6)}+
$$
$$
+
{\bf 1}_{\{i_1=i_5\ne 0\}}
{\bf 1}_{\{j_1=j_5\}}
{\bf 1}_{\{i_2=i_4\ne 0\}}
{\bf 1}_{\{j_2=j_4\}}
\zeta_{j_3}^{(i_3)}
\zeta_{j_6}^{(i_6)}
+
{\bf 1}_{\{i_1=i_5\ne 0\}}
{\bf 1}_{\{j_1=j_5\}}
{\bf 1}_{\{i_3=i_4\ne 0\}}
{\bf 1}_{\{j_3=j_4\}}
\zeta_{j_2}^{(i_2)}
\zeta_{j_6}^{(i_6)}+
$$
$$
+
{\bf 1}_{\{i_2=i_3\ne 0\}}
{\bf 1}_{\{j_2=j_3\}}
{\bf 1}_{\{i_4=i_5\ne 0\}}
{\bf 1}_{\{j_4=j_5\}}
\zeta_{j_1}^{(i_1)}
\zeta_{j_6}^{(i_6)}
+
{\bf 1}_{\{i_2=i_4\ne 0\}}
{\bf 1}_{\{j_2=j_4\}}
{\bf 1}_{\{i_3=i_5\ne 0\}}
{\bf 1}_{\{j_3=j_5\}}
\zeta_{j_1}^{(i_1)}
\zeta_{j_6}^{(i_6)}+
$$
$$
+
{\bf 1}_{\{i_2=i_5\ne 0\}}
{\bf 1}_{\{j_2=j_5\}}
{\bf 1}_{\{i_3=i_4\ne 0\}}
{\bf 1}_{\{j_3=j_4\}}
\zeta_{j_1}^{(i_1)}
\zeta_{j_6}^{(i_6)}
+
{\bf 1}_{\{i_6=i_1\ne 0\}}
{\bf 1}_{\{j_6=j_1\}}
{\bf 1}_{\{i_3=i_4\ne 0\}}
{\bf 1}_{\{j_3=j_4\}}
\zeta_{j_2}^{(i_2)}
\zeta_{j_5}^{(i_5)}+
$$
$$
+
{\bf 1}_{\{i_6=i_1\ne 0\}}
{\bf 1}_{\{j_6=j_1\}}
{\bf 1}_{\{i_3=i_5\ne 0\}}
{\bf 1}_{\{j_3=j_5\}}
\zeta_{j_2}^{(i_2)}
\zeta_{j_4}^{(i_4)}
+
{\bf 1}_{\{i_6=i_1\ne 0\}}
{\bf 1}_{\{j_6=j_1\}}
{\bf 1}_{\{i_2=i_5\ne 0\}}
{\bf 1}_{\{j_2=j_5\}}
\zeta_{j_3}^{(i_3)}
\zeta_{j_4}^{(i_4)}+
$$
$$
+
{\bf 1}_{\{i_6=i_1\ne 0\}}
{\bf 1}_{\{j_6=j_1\}}
{\bf 1}_{\{i_2=i_4\ne 0\}}
{\bf 1}_{\{j_2=j_4\}}
\zeta_{j_3}^{(i_3)}
\zeta_{j_5}^{(i_5)}
+
{\bf 1}_{\{i_6=i_1\ne 0\}}
{\bf 1}_{\{j_6=j_1\}}
{\bf 1}_{\{i_4=i_5\ne 0\}}
{\bf 1}_{\{j_4=j_5\}}
\zeta_{j_2}^{(i_2)}
\zeta_{j_3}^{(i_3)}+
$$
$$
+
{\bf 1}_{\{i_6=i_1\ne 0\}}
{\bf 1}_{\{j_6=j_1\}}
{\bf 1}_{\{i_2=i_3\ne 0\}}
{\bf 1}_{\{j_2=j_3\}}
\zeta_{j_4}^{(i_4)}
\zeta_{j_5}^{(i_5)}
+
{\bf 1}_{\{i_6=i_2\ne 0\}}
{\bf 1}_{\{j_6=j_2\}}
{\bf 1}_{\{i_3=i_5\ne 0\}}
{\bf 1}_{\{j_3=j_5\}}
\zeta_{j_1}^{(i_1)}
\zeta_{j_4}^{(i_4)}+
$$
$$
+
{\bf 1}_{\{i_6=i_2\ne 0\}}
{\bf 1}_{\{j_6=j_2\}}
{\bf 1}_{\{i_4=i_5\ne 0\}}
{\bf 1}_{\{j_4=j_5\}}
\zeta_{j_1}^{(i_1)}
\zeta_{j_3}^{(i_3)}
+
{\bf 1}_{\{i_6=i_2\ne 0\}}
{\bf 1}_{\{j_6=j_2\}}
{\bf 1}_{\{i_3=i_4\ne 0\}}
{\bf 1}_{\{j_3=j_4\}}
\zeta_{j_1}^{(i_1)}
\zeta_{j_5}^{(i_5)}+
$$
$$
+
{\bf 1}_{\{i_6=i_2\ne 0\}}
{\bf 1}_{\{j_6=j_2\}}
{\bf 1}_{\{i_1=i_5\ne 0\}}
{\bf 1}_{\{j_1=j_5\}}
\zeta_{j_3}^{(i_3)}
\zeta_{j_4}^{(i_4)}
+
{\bf 1}_{\{i_6=i_2\ne 0\}}
{\bf 1}_{\{j_6=j_2\}}
{\bf 1}_{\{i_1=i_4\ne 0\}}
{\bf 1}_{\{j_1=j_4\}}
\zeta_{j_3}^{(i_3)}
\zeta_{j_5}^{(i_5)}+
$$
$$
+
{\bf 1}_{\{i_6=i_2\ne 0\}}
{\bf 1}_{\{j_6=j_2\}}
{\bf 1}_{\{i_1=i_3\ne 0\}}
{\bf 1}_{\{j_1=j_3\}}
\zeta_{j_4}^{(i_4)}
\zeta_{j_5}^{(i_5)}
+
{\bf 1}_{\{i_6=i_3\ne 0\}}
{\bf 1}_{\{j_6=j_3\}}
{\bf 1}_{\{i_2=i_5\ne 0\}}
{\bf 1}_{\{j_2=j_5\}}
\zeta_{j_1}^{(i_1)}
\zeta_{j_4}^{(i_4)}+
$$
$$
+
{\bf 1}_{\{i_6=i_3\ne 0\}}
{\bf 1}_{\{j_6=j_3\}}
{\bf 1}_{\{i_4=i_5\ne 0\}}
{\bf 1}_{\{j_4=j_5\}}
\zeta_{j_1}^{(i_1)}
\zeta_{j_2}^{(i_2)}
+
{\bf 1}_{\{i_6=i_3\ne 0\}}
{\bf 1}_{\{j_6=j_3\}}
{\bf 1}_{\{i_2=i_4\ne 0\}}
{\bf 1}_{\{j_2=j_4\}}
\zeta_{j_1}^{(i_1)}
\zeta_{j_5}^{(i_5)}+
$$
$$
+
{\bf 1}_{\{i_6=i_3\ne 0\}}
{\bf 1}_{\{j_6=j_3\}}
{\bf 1}_{\{i_1=i_5\ne 0\}}
{\bf 1}_{\{j_1=j_5\}}
\zeta_{j_2}^{(i_2)}
\zeta_{j_4}^{(i_4)}
+
{\bf 1}_{\{i_6=i_3\ne 0\}}
{\bf 1}_{\{j_6=j_3\}}
{\bf 1}_{\{i_1=i_4\ne 0\}}
{\bf 1}_{\{j_1=j_4\}}
\zeta_{j_2}^{(i_2)}
\zeta_{j_5}^{(i_5)}+
$$
$$
+
{\bf 1}_{\{i_6=i_3\ne 0\}}
{\bf 1}_{\{j_6=j_3\}}
{\bf 1}_{\{i_1=i_2\ne 0\}}
{\bf 1}_{\{j_1=j_2\}}
\zeta_{j_4}^{(i_4)}
\zeta_{j_5}^{(i_5)}
+
{\bf 1}_{\{i_6=i_4\ne 0\}}
{\bf 1}_{\{j_6=j_4\}}
{\bf 1}_{\{i_3=i_5\ne 0\}}
{\bf 1}_{\{j_3=j_5\}}
\zeta_{j_1}^{(i_1)}
\zeta_{j_2}^{(i_2)}+
$$
$$
+
{\bf 1}_{\{i_6=i_4\ne 0\}}
{\bf 1}_{\{j_6=j_4\}}
{\bf 1}_{\{i_2=i_5\ne 0\}}
{\bf 1}_{\{j_2=j_5\}}
\zeta_{j_1}^{(i_1)}
\zeta_{j_3}^{(i_3)}
+
{\bf 1}_{\{i_6=i_4\ne 0\}}
{\bf 1}_{\{j_6=j_4\}}
{\bf 1}_{\{i_2=i_3\ne 0\}}
{\bf 1}_{\{j_2=j_3\}}
\zeta_{j_1}^{(i_1)}
\zeta_{j_5}^{(i_5)}+
$$
$$
+
{\bf 1}_{\{i_6=i_4\ne 0\}}
{\bf 1}_{\{j_6=j_4\}}
{\bf 1}_{\{i_1=i_5\ne 0\}}
{\bf 1}_{\{j_1=j_5\}}
\zeta_{j_2}^{(i_2)}
\zeta_{j_3}^{(i_3)}
+
{\bf 1}_{\{i_6=i_4\ne 0\}}
{\bf 1}_{\{j_6=j_4\}}
{\bf 1}_{\{i_1=i_3\ne 0\}}
{\bf 1}_{\{j_1=j_3\}}
\zeta_{j_2}^{(i_2)}
\zeta_{j_5}^{(i_5)}+
$$
$$
+
{\bf 1}_{\{i_6=i_4\ne 0\}}
{\bf 1}_{\{j_6=j_4\}}
{\bf 1}_{\{i_1=i_2\ne 0\}}
{\bf 1}_{\{j_1=j_2\}}
\zeta_{j_3}^{(i_3)}
\zeta_{j_5}^{(i_5)}
+
{\bf 1}_{\{i_6=i_5\ne 0\}}
{\bf 1}_{\{j_6=j_5\}}
{\bf 1}_{\{i_3=i_4\ne 0\}}
{\bf 1}_{\{j_3=j_4\}}
\zeta_{j_1}^{(i_1)}
\zeta_{j_2}^{(i_2)}+
$$
$$
+
{\bf 1}_{\{i_6=i_5\ne 0\}}
{\bf 1}_{\{j_6=j_5\}}
{\bf 1}_{\{i_2=i_4\ne 0\}}
{\bf 1}_{\{j_2=j_4\}}
\zeta_{j_1}^{(i_1)}
\zeta_{j_3}^{(i_3)}
+
{\bf 1}_{\{i_6=i_5\ne 0\}}
{\bf 1}_{\{j_6=j_5\}}
{\bf 1}_{\{i_2=i_3\ne 0\}}
{\bf 1}_{\{j_2=j_3\}}
\zeta_{j_1}^{(i_1)}
\zeta_{j_4}^{(i_4)}+
$$
$$
+
{\bf 1}_{\{i_6=i_5\ne 0\}}
{\bf 1}_{\{j_6=j_5\}}
{\bf 1}_{\{i_1=i_4\ne 0\}}
{\bf 1}_{\{j_1=j_4\}}
\zeta_{j_2}^{(i_2)}
\zeta_{j_3}^{(i_3)}
+
{\bf 1}_{\{i_6=i_5\ne 0\}}
{\bf 1}_{\{j_6=j_5\}}
{\bf 1}_{\{i_1=i_3\ne 0\}}
{\bf 1}_{\{j_1=j_3\}}
\zeta_{j_2}^{(i_2)}
\zeta_{j_4}^{(i_4)}+
$$
$$
+
{\bf 1}_{\{i_6=i_5\ne 0\}}
{\bf 1}_{\{j_6=j_5\}}
{\bf 1}_{\{i_1=i_2\ne 0\}}
{\bf 1}_{\{j_1=j_2\}}
\zeta_{j_3}^{(i_3)}
\zeta_{j_4}^{(i_4)}-
$$
$$
-
{\bf 1}_{\{i_6=i_1\ne 0\}}
{\bf 1}_{\{j_6=j_1\}}
{\bf 1}_{\{i_2=i_5\ne 0\}}
{\bf 1}_{\{j_2=j_5\}}
{\bf 1}_{\{i_3=i_4\ne 0\}}
{\bf 1}_{\{j_3=j_4\}}-
$$
$$
-
{\bf 1}_{\{i_6=i_1\ne 0\}}
{\bf 1}_{\{j_6=j_1\}}
{\bf 1}_{\{i_2=i_4\ne 0\}}
{\bf 1}_{\{j_2=j_4\}}
{\bf 1}_{\{i_3=i_5\ne 0\}}
{\bf 1}_{\{j_3=j_5\}}-
$$
$$
-
{\bf 1}_{\{i_6=i_1\ne 0\}}
{\bf 1}_{\{j_6=j_1\}}
{\bf 1}_{\{i_2=i_3\ne 0\}}
{\bf 1}_{\{j_2=j_3\}}
{\bf 1}_{\{i_4=i_5\ne 0\}}
{\bf 1}_{\{j_4=j_5\}}-
$$
$$
-
{\bf 1}_{\{i_6=i_2\ne 0\}}
{\bf 1}_{\{j_6=j_2\}}
{\bf 1}_{\{i_1=i_5\ne 0\}}
{\bf 1}_{\{j_1=j_5\}}
{\bf 1}_{\{i_3=i_4\ne 0\}}
{\bf 1}_{\{j_3=j_4\}}-
$$
$$
-
{\bf 1}_{\{i_6=i_2\ne 0\}}
{\bf 1}_{\{j_6=j_2\}}
{\bf 1}_{\{i_1=i_4\ne 0\}}
{\bf 1}_{\{j_1=j_4\}}
{\bf 1}_{\{i_3=i_5\ne 0\}}
{\bf 1}_{\{j_3=j_5\}}-
$$
$$
-
{\bf 1}_{\{i_6=i_2\ne 0\}}
{\bf 1}_{\{j_6=j_2\}}
{\bf 1}_{\{i_1=i_3\ne 0\}}
{\bf 1}_{\{j_1=j_3\}}
{\bf 1}_{\{i_4=i_5\ne 0\}}
{\bf 1}_{\{j_4=j_5\}}-
$$
$$
-
{\bf 1}_{\{i_6=i_3\ne 0\}}
{\bf 1}_{\{j_6=j_3\}}
{\bf 1}_{\{i_1=i_5\ne 0\}}
{\bf 1}_{\{j_1=j_5\}}
{\bf 1}_{\{i_2=i_4\ne 0\}}
{\bf 1}_{\{j_2=j_4\}}-
$$
$$
-
{\bf 1}_{\{i_6=i_3\ne 0\}}
{\bf 1}_{\{j_6=j_3\}}
{\bf 1}_{\{i_1=i_4\ne 0\}}
{\bf 1}_{\{j_1=j_4\}}
{\bf 1}_{\{i_2=i_5\ne 0\}}
{\bf 1}_{\{j_2=j_5\}}-
$$
$$
-
{\bf 1}_{\{i_3=i_6\ne 0\}}
{\bf 1}_{\{j_3=j_6\}}
{\bf 1}_{\{i_1=i_2\ne 0\}}
{\bf 1}_{\{j_1=j_2\}}
{\bf 1}_{\{i_4=i_5\ne 0\}}
{\bf 1}_{\{j_4=j_5\}}-
$$
$$
-
{\bf 1}_{\{i_6=i_4\ne 0\}}
{\bf 1}_{\{j_6=j_4\}}
{\bf 1}_{\{i_1=i_5\ne 0\}}
{\bf 1}_{\{j_1=j_5\}}
{\bf 1}_{\{i_2=i_3\ne 0\}}
{\bf 1}_{\{j_2=j_3\}}-
$$
$$
-
{\bf 1}_{\{i_6=i_4\ne 0\}}
{\bf 1}_{\{j_6=j_4\}}
{\bf 1}_{\{i_1=i_3\ne 0\}}
{\bf 1}_{\{j_1=j_3\}}
{\bf 1}_{\{i_2=i_5\ne 0\}}
{\bf 1}_{\{j_2=j_5\}}-
$$
$$
-
{\bf 1}_{\{i_6=i_4\ne 0\}}
{\bf 1}_{\{j_6=j_4\}}
{\bf 1}_{\{i_1=i_2\ne 0\}}
{\bf 1}_{\{j_1=j_2\}}
{\bf 1}_{\{i_3=i_5\ne 0\}}
{\bf 1}_{\{j_3=j_5\}}-
$$
$$
-
{\bf 1}_{\{i_6=i_5\ne 0\}}
{\bf 1}_{\{j_6=j_5\}}
{\bf 1}_{\{i_1=i_4\ne 0\}}
{\bf 1}_{\{j_1=j_4\}}
{\bf 1}_{\{i_2=i_3\ne 0\}}
{\bf 1}_{\{j_2=j_3\}}-
$$
$$
-
{\bf 1}_{\{i_6=i_5\ne 0\}}
{\bf 1}_{\{j_6=j_5\}}
{\bf 1}_{\{i_1=i_2\ne 0\}}
{\bf 1}_{\{j_1=j_2\}}
{\bf 1}_{\{i_3=i_4\ne 0\}}
{\bf 1}_{\{j_3=j_4\}}-
$$
\begin{equation}
\label{a6}
\Biggl.-
{\bf 1}_{\{i_6=i_5\ne 0\}}
{\bf 1}_{\{j_6=j_5\}}
{\bf 1}_{\{i_1=i_3\ne 0\}}
{\bf 1}_{\{j_1=j_3\}}
{\bf 1}_{\{i_2=i_4\ne 0\}}
{\bf 1}_{\{j_2=j_4\}}\Biggr),
\end{equation}

\vspace{6mm}
\noindent
where ${\bf 1}_A$ is the indicator of the set $A$.

As a result we obtain the following new possibilities and advantages 
compared with the method based on the Milstein approach \cite{19}.

There is the explicit formula (see (\ref{ppppa})) for calculation 
of expansion coefficients 
of the iterated Ito stochastic integral (\ref{ito}) with any
fixed multiplicity $k$. 

We have new possibilities for exact calculation of the mean-square 
error of approximation of the iterated Ito stochastic integrals (\ref{ito})
(see Theorem 8 below).

Since the used multiple Fourier series is a generalized in the sense
that it is built using various complete orthonormal
systems of functions in the space $L_2([t, T])$, we have new possibilities 
for approximation --- we can
use not only trigonometric functions as in the Milstein approach \cite{19}
but Legendre polynomials.
As it turned out (see below), it is more convenient to work 
with Legendre polynomials for constructing approximations of the iterated
Ito stochastic integrals (\ref{ito}). We can choose different numbers $q$ 
(see Sect.~7) for
approximations of different iterated Ito stochastic integrals
from the family (\ref{ito}).
This is impossible for approximations based on the Milstein approach
\cite{19}.
Approximations based on Legendre polynomials essentially simpler 
than approximations based on trigonometric functions
(see (\ref{444}), (\ref{1970}), (\ref{5}), (\ref{6})).

As we mentioned before,
the Milstein approach \cite{19} based on the 
Karhunen--Loeve expansion of the Brownian bridge process
leads to iterated series (in contrast with multiple 
series from Theorems 1--7) starting at least from the second or third 
multiplicity of iterated stochastic integrals.
Multiple series are more convenient for approximation than the iterated ones, 
since partial sums of multiple series converge for any possible case of  
convergence to infinity of their upper limits of summation 
(let us denote them as $p_1,\ldots, p_k$). 
For example,
when $p_1=\ldots=p_k=p\to\infty$. 
For iterated series, the condition $p_1=\ldots=p_k=p\to\infty$ obviously 
does not guarantee the convergence of this series.

However, 
in \cite{20} (pp.~438-439),
\cite{21}
(Sect.~5.8, pp.~202--204), \cite{22} (pp.~82-84),
\cite{Zapad-9} (pp.~263-264) 
the authors use (without rigorous proof)
the condition $p_1=p_2=p_3=q\to\infty$
within the frames of the mentioned approach
based on the Karhunen--Loeve expansion of the Brownian bridge
process \cite{19} together with
the Wong--Zakai approximation
\cite{W-Z-1}-\cite{Watanabe}
(see discussion in Sect.~11 for details).

For further consideration, let us 
consider the generalization of formulas (\ref{a1})--(\ref{a6})                 
for the case of an arbitrary multiplicity $k$ $(k\in\mathbb{N})$ of 
the iterated Ito stochastic integral $J[\psi^{(k)}]_{T,t}$ defined by (\ref{ito}).
In order to do this, let us
introduce some notations. 
Consider the unordered
set $\{1, 2, \ldots, k\}$ 
and separate it into two parts:
the first part consists of $r$ unordered 
pairs (sequence order of these pairs is also unimportant) and the 
second one consists of the 
remaining $k-2r$ numbers.
So, we have

\begin{equation}
\label{leto5007}
(\{
\underbrace{\{g_1, g_2\}, \ldots, 
\{g_{2r-1}, g_{2r}\}}_{\small{\hbox{part 1}}}
\},
\{\underbrace{q_1, \ldots, q_{k-2r}}_{\small{\hbox{part 2}}}
\}),
\end{equation}

\vspace{4mm}
\noindent
where 

\vspace{-2mm}
$$
\{g_1, g_2, \ldots, 
g_{2r-1}, g_{2r}, q_1, \ldots, q_{k-2r}\}=\{1, 2, \ldots, k\},
$$

\vspace{4mm}
\noindent
braces   
mean an unordered 
set, and pa\-ren\-the\-ses mean an ordered set.

We will say that (\ref{leto5007}) is a partition 
and consider the sum with respect to all possible
partitions

\begin{equation}
\label{leto5008}
\sum_{\stackrel{(\{\{g_1, g_2\}, \ldots, 
\{g_{2r-1}, g_{2r}\}\}, \{q_1, \ldots, q_{k-2r}\})}
{{}_{\{g_1, g_2, \ldots, 
g_{2r-1}, g_{2r}, q_1, \ldots, q_{k-2r}\}=\{1, 2, \ldots, k\}}}}
a_{g_1 g_2, \ldots, 
g_{2r-1} g_{2r}, q_1 \ldots q_{k-2r}}.
\end{equation}

\vspace{4mm}

Below there are several examples of sums in the form (\ref{leto5008})

\vspace{2mm}
$$
\sum_{\stackrel{(\{g_1, g_2\})}{{}_{\{g_1, g_2\}=\{1, 2\}}}}
a_{g_1 g_2}=a_{12},
$$

\vspace{3mm}
$$
\sum_{\stackrel{(\{\{g_1, g_2\}, \{g_3, g_4\}\})}
{{}_{\{g_1, g_2, g_3, g_4\}=\{1, 2, 3, 4\}}}}
a_{g_1 g_2, g_3 g_4}=a_{12,34} + a_{13,24} + a_{23,14},
$$

\vspace{3mm}
$$
\sum_{\stackrel{(\{g_1, g_2\}, \{q_1, q_{2}\})}
{{}_{\{g_1, g_2, q_1, q_{2}\}=\{1, 2, 3, 4\}}}}
a_{g_1 g_2, q_1 q_{2}}=
$$

$$
=a_{12,34}+a_{13,24}+a_{14,23}
+a_{23,14}+a_{24,13}+a_{34,12},
$$

\vspace{3mm}
$$
\sum_{\stackrel{(\{g_1, g_2\}, \{q_1, q_{2}, q_3\})}
{{}_{\{g_1, g_2, q_1, q_{2}, q_3\}=\{1, 2, 3, 4, 5\}}}}
a_{g_1 g_2, q_1 q_{2}q_3}
=
$$

$$
=a_{12,345}+a_{13,245}+a_{14,235}
+a_{15,234}+a_{23,145}+a_{24,135}+
$$
$$
+a_{25,134}+a_{34,125}+a_{35,124}+a_{45,123},
$$

\vspace{4mm}
$$
\sum_{\stackrel{(\{\{g_1, g_2\}, \{g_3, g_{4}\}\}, \{q_1\})}
{{}_{\{g_1, g_2, g_3, g_{4}, q_1\}=\{1, 2, 3, 4, 5\}}}}
a_{g_1 g_2, g_3 g_{4},q_1}
=
$$

$$
=
a_{12,34,5}+a_{13,24,5}+a_{14,23,5}+
a_{12,35,4}+a_{13,25,4}+a_{15,23,4}+
$$
$$
+a_{12,54,3}+a_{15,24,3}+a_{14,25,3}+a_{15,34,2}+a_{13,54,2}+a_{14,53,2}+
$$
$$
+
a_{52,34,1}+a_{53,24,1}+a_{54,23,1}.
$$

\vspace{7mm}

Now we can write (\ref{tyyy}) as

\vspace{1mm}

$$
J[\psi^{(k)}]_{T,t}=
\hbox{\vtop{\offinterlineskip\halign{
\hfil#\hfil\cr
{\rm l.i.m.}\cr
$\stackrel{}{{}_{p_1,\ldots,p_k\to \infty}}$\cr
}} }
\sum\limits_{j_1=0}^{p_1}\ldots
\sum\limits_{j_k=0}^{p_k}
C_{j_k\ldots j_1}\Biggl(
\prod_{l=1}^k\zeta_{j_l}^{(i_l)}+\sum\limits_{r=1}^{[k/2]}
(-1)^r \times
\Biggr.
$$

\vspace{3mm}
\begin{equation}
\label{leto6000hh}
\times
\sum_{\stackrel{(\{\{g_1, g_2\}, \ldots, 
\{g_{2r-1}, g_{2r}\}\}, \{q_1, \ldots, q_{k-2r}\})}
{{}_{\{g_1, g_2, \ldots, 
g_{2r-1}, g_{2r}, q_1, \ldots, q_{k-2r}\}=\{1, 2, \ldots, k\}}}}
\prod\limits_{s=1}^r
{\bf 1}_{\{i_{g_{{}_{2s-1}}}=~i_{g_{{}_{2s}}}\ne 0\}}
\Biggl.{\bf 1}_{\{j_{g_{{}_{2s-1}}}=~j_{g_{{}_{2s}}}\}}
\prod_{l=1}^{k-2r}\zeta_{j_{q_l}}^{(i_{q_l})}\Biggr),
\end{equation}

\vspace{5mm}
\noindent
where $[x]$ is an integer part of a real number $x;$
another notations are the same as in Theorem {\bf 1}.

\vspace{2mm}

In particular, from (\ref{leto6000hh}) for $k=5$ we obtain

\vspace{3mm}

$$
J[\psi^{(5)}]_{T,t}=
\hbox{\vtop{\offinterlineskip\halign{
\hfil#\hfil\cr
{\rm l.i.m.}\cr
$\stackrel{}{{}_{p_1,\ldots,p_5\to \infty}}$\cr
}} }\sum_{j_1=0}^{p_1}\ldots\sum_{j_5=0}^{p_5}
C_{j_5\ldots j_1}\Biggl(
\prod_{l=1}^5\zeta_{j_l}^{(i_l)}-\Biggr.
$$

\vspace{2mm}
$$
-
\sum\limits_{\stackrel{(\{g_1, g_2\}, \{q_1, q_{2}, q_3\})}
{{}_{\{g_1, g_2, q_{1}, q_{2}, q_3\}=\{1, 2, 3, 4, 5\}}}}
{\bf 1}_{\{i_{g_{{}_{1}}}=~i_{g_{{}_{2}}}\ne 0\}}
{\bf 1}_{\{j_{g_{{}_{1}}}=~j_{g_{{}_{2}}}\}}
\prod_{l=1}^{3}\zeta_{j_{q_l}}^{(i_{q_l})}+
$$

\vspace{2mm}
$$
+
\sum_{\stackrel{(\{\{g_1, g_2\}, 
\{g_{3}, g_{4}\}\}, \{q_1\})}
{{}_{\{g_1, g_2, g_{3}, g_{4}, q_1\}=\{1, 2, 3, 4, 5\}}}}
{\bf 1}_{\{i_{g_{{}_{1}}}=~i_{g_{{}_{2}}}\ne 0\}}
{\bf 1}_{\{j_{g_{{}_{1}}}=~j_{g_{{}_{2}}}\}}
\Biggl.{\bf 1}_{\{i_{g_{{}_{3}}}=~i_{g_{{}_{4}}}\ne 0\}}
{\bf 1}_{\{j_{g_{{}_{3}}}=~j_{g_{{}_{4}}}\}}
\zeta_{j_{q_1}}^{(i_{q_1})}\Biggr).
$$

\vspace{7mm}
\noindent
The last equality obviously agrees with
(\ref{a5}).

Let us consider the generalization of Theorem 1 for the case
of an arbitrary complete orthonormal systems  
of functions in the space $L_2([t,T])$ 
and $\psi_1(\tau),\ldots,\psi_k(\tau)\in L_2([t, T]).$

\vspace{2mm}

{\bf Theorem~2}\ \cite{10a} (Sect.~1.11), \cite{11} (Sect.~15).
{\it Suppose that
$\psi_1(\tau),\ldots,\psi_k(\tau)\in L_2([t, T])$ and
$\{\phi_j(x)\}_{j=0}^{\infty}$ is an arbitrary complete orthonormal system  
of functions in the space $L_2([t,T]).$
Then the following expansion

\vspace{1mm}
$$
J[\psi^{(k)}]_{T,t}=
\hbox{\vtop{\offinterlineskip\halign{
\hfil#\hfil\cr
{\rm l.i.m.}\cr
$\stackrel{}{{}_{p_1,\ldots,p_k\to \infty}}$\cr
}} }
\sum\limits_{j_1=0}^{p_1}\ldots
\sum\limits_{j_k=0}^{p_k}
C_{j_k\ldots j_1}\Biggl(
\prod_{l=1}^k\zeta_{j_l}^{(i_l)}+\sum\limits_{r=1}^{[k/2]}
(-1)^r \times
\Biggr.
$$

\vspace{2mm}
\begin{equation}
\label{leto6000}
\times
\sum_{\stackrel{(\{\{g_1, g_2\}, \ldots, 
\{g_{2r-1}, g_{2r}\}\}, \{q_1, \ldots, q_{k-2r}\})}
{{}_{\{g_1, g_2, \ldots, 
g_{2r-1}, g_{2r}, q_1, \ldots, q_{k-2r}\}=\{1, 2, \ldots, k\}}}}
\prod\limits_{s=1}^r
{\bf 1}_{\{i_{g_{{}_{2s-1}}}=~i_{g_{{}_{2s}}}\ne 0\}}
\Biggl.{\bf 1}_{\{j_{g_{{}_{2s-1}}}=~j_{g_{{}_{2s}}}\}}
\prod_{l=1}^{k-2r}\zeta_{j_{q_l}}^{(i_{q_l})}\Biggr)
\end{equation}

\vspace{6mm}
\noindent
con\-verg\-ing in the mean-square sense is valid,
where $[x]$ is an integer part of a real number $x;$
another notations are the same as in Theorem~{\rm 1}.}

\vspace{2mm}

It should be noted that an analogue of Theorem 2 was considered 
in \cite{Rybakov1000}. 
Note that we use another notations 
\cite{10a} (Sect.~1.11), \cite{11} (Sect.~15)
in comparison with \cite{Rybakov1000}.
Moreover, the proof of an analogue of Theorem 2
from \cite{Rybakov1000} is somewhat different from the proof given in 
\cite{10a} (Sect.~1.11), \cite{11} (Sect.~15).

\vspace{5mm}

\section{Expansions of Iterated Stratonovich Stochastic Integrals of Multiplicities 1 to 6}

\vspace{5mm}

In a number the author's works \cite{5}-\cite{10aaa}, \cite{12}, \cite{15a}
Theorems 1, 2 have been adapted for the integrals
(\ref{str}) of multiplicities 2 to 4.
Let us collect some old results in the following theorem.

\vspace{2mm}

{\bf Theorem 3} \cite{5}-\cite{10aaa}, \cite{12}, \cite{15a}. 
{\it Suppose that 
$\{\phi_j(x)\}_{j=0}^{\infty}$ is a complete orthonormal system of 
Legendre polynomials or trigonometric functions in the space $L_2([t, T]).$
At the same time $\psi_2(\tau)$ is a continuously differentiable 
function on $[t, T]$ and $\psi_1(\tau), \psi_3(\tau)$ are twice
continuously differentiable functions on $[t, T]$. Then

\begin{equation}
\label{a}
J^{*}[\psi^{(2)}]_{T,t}=
\hbox{\vtop{\offinterlineskip\halign{
\hfil#\hfil\cr
{\rm l.i.m.}\cr
$\stackrel{}{{}_{p_1,p_2\to \infty}}$\cr
}} }\sum_{j_1=0}^{p_1}\sum_{j_2=0}^{p_2}
C_{j_2j_1}\zeta_{j_1}^{(i_1)}\zeta_{j_2}^{(i_2)}\ \ \ (i_1,i_2=1,\ldots,m),
\end{equation}

\vspace{1mm}
\begin{equation}
\label{feto19000ab}
J^{*}[\psi^{(3)}]_{T,t}=
\hbox{\vtop{\offinterlineskip\halign{
\hfil#\hfil\cr
{\rm l.i.m.}\cr
$\stackrel{}{{}_{p_1,p_2,p_3\to \infty}}$\cr
}} }\sum_{j_1=0}^{p_1}\sum_{j_2=0}^{p_2}\sum_{j_3=0}^{p_3}
C_{j_3 j_2 j_1}\zeta_{j_1}^{(i_1)}\zeta_{j_2}^{(i_2)}\zeta_{j_3}^{(i_3)}\ \ \
(i_1,i_2,i_3=0, 1,\ldots,m),
\end{equation}

\vspace{1mm}
\begin{equation}
\label{feto19000a}
J^{*}[\psi^{(3)}]_{T,t}=
\hbox{\vtop{\offinterlineskip\halign{
\hfil#\hfil\cr
{\rm l.i.m.}\cr
$\stackrel{}{{}_{p\to \infty}}$\cr
}} }
\sum\limits_{j_1,j_2,j_3=0}^{p}
C_{j_3 j_2 j_1}\zeta_{j_1}^{(i_1)}\zeta_{j_2}^{(i_2)}\zeta_{j_3}^{(i_3)}\ \ \
(i_1,i_2,i_3=1,\ldots,m),
\end{equation}

\vspace{1mm}
\begin{equation}
\label{uu}
J^{*}[\psi^{(4)}]_{T,t}=
\hbox{\vtop{\offinterlineskip\halign{
\hfil#\hfil\cr
{\rm l.i.m.}\cr
$\stackrel{}{{}_{p\to \infty}}$\cr
}} }
\sum\limits_{j_1, \ldots, j_4=0}^{p}
C_{j_4 j_3 j_2 j_1}\zeta_{j_1}^{(i_1)}
\zeta_{j_2}^{(i_2)}\zeta_{j_3}^{(i_3)}\zeta_{j_4}^{(i_4)}\ \ \
(i_1,\ldots,i_4=0, 1,\ldots,m),
\end{equation}

\vspace{6mm}
\noindent
where $J^{*}[\psi^{(k)}]_{T,t}$ is defined by {\rm (\ref{str})}, and
$\psi_l(\tau)\equiv 1$ $(l=1,\ldots,4)$ in {\rm (\ref{feto19000ab})}, 
{\rm (\ref{uu});} another notations are the same as in Theorems {\rm 1, 2.}
}

\vspace{2mm}

Recently, a new approach to the expansion and mean-square 
approximation of iterated Stratonovich stochastic integrals has been obtained
\cite{10a} (Sect.~2.10--2.16), \cite{12} (Sect.~13--19), 
\cite{arxiv-11} (Sect.~7--13), \cite{15a} (Sect.~5--11),
\cite{new-art-1-xxy} (Sect.~4--9).
Let us formulate four theorems that were proved using this approach.

\vspace{2mm}

{\bf Theorem 4}\ \cite{10a}, \cite{12}, \cite{arxiv-11}, \cite{15a}, \cite{new-art-1-xxy}.\
{\it Suppose 
that $\{\phi_j(x)\}_{j=0}^{\infty}$ is a complete orthonormal system of 
Legendre polynomials or trigonometric functions in the space $L_2([t, T]).$
Furthermore, let $\psi_1(\tau), \psi_2(\tau),$ $\psi_3(\tau)$ are continuously dif\-ferentiable 
nonrandom functions on $[t, T].$ 
Then, for the 
iterated Stra\-to\-no\-vich stochastic integral of third multiplicity

\vspace{-1mm}
$$
J^{*}[\psi^{(3)}]_{T,t}={\int\limits_t^{*}}^T\psi_3(t_3)
{\int\limits_t^{*}}^{t_3}\psi_2(t_2)
{\int\limits_t^{*}}^{t_2}\psi_1(t_1)
d{\bf w}_{t_1}^{(i_1)}
d{\bf w}_{t_2}^{(i_2)}d{\bf w}_{t_3}^{(i_3)}\ \ \ (i_1,i_2,i_3=0,1,\ldots,m)
$$

\vspace{3mm}
\noindent
the following 
relations

\vspace{-2mm}
\begin{equation}
\label{fin1}
J^{*}[\psi^{(3)}]_{T,t}
=\hbox{\vtop{\offinterlineskip\halign{
\hfil#\hfil\cr
{\rm l.i.m.}\cr
$\stackrel{}{{}_{p\to \infty}}$\cr
}} }
\sum\limits_{j_1, j_2, j_3=0}^{p}
C_{j_3 j_2 j_1}\zeta_{j_1}^{(i_1)}\zeta_{j_2}^{(i_2)}\zeta_{j_3}^{(i_3)},
\end{equation}

\vspace{2mm}
\begin{equation}
\label{fin2}
{\sf M}\left\{\left(
J^{*}[\psi^{(3)}]_{T,t}-
\sum\limits_{j_1, j_2, j_3=0}^{p}
C_{j_3 j_2 j_1}\zeta_{j_1}^{(i_1)}\zeta_{j_2}^{(i_2)}\zeta_{j_3}^{(i_3)}\right)^2\right\}
\le \frac{C}{p}
\end{equation}

\vspace{4mm}
\noindent
are fulfilled, where $i_1, i_2, i_3=0,1,\ldots,m$ in {\rm (\ref{fin1})} and 
$i_1, i_2, i_3=1,\ldots,m$ in {\rm (\ref{fin2})},
constant $C$ is independent of $p,$

\vspace{-1mm}
$$
C_{j_3 j_2 j_1}=\int\limits_t^T\psi_3(t_3)\phi_{j_3}(t_3)
\int\limits_t^{t_3}\psi_2(t_2)\phi_{j_2}(t_2)
\int\limits_t^{t_2}\psi_1(t_1)\phi_{j_1}(t_1)dt_1dt_2dt_3
$$

\vspace{3mm}
\noindent
and
$$
\zeta_{j}^{(i)}=
\int\limits_t^T \phi_{j}(\tau) d{\bf f}_{\tau}^{(i)}
$$ 

\vspace{2mm}
\noindent
are independent standard Gaussian random variables for various 
$i$ or $j$ {\rm (}in the case when $i\ne 0${\rm );} 
another notations are the same as in Theorems~{\rm 1, 2}.}

\vspace{2mm}

{\bf Theorem 5}\ \cite{10a}, \cite{12}, \cite{arxiv-11}, \cite{15a}, \cite{new-art-1-xxy}.\
{\it Let
$\{\phi_j(x)\}_{j=0}^{\infty}$ be a complete orthonormal system of 
Legendre polynomials or trigonometric functions in the space $L_2([t, T]).$
Furthermore, let $\psi_1(\tau), \ldots,$ $\psi_4(\tau)$ be continuously dif\-ferentiable 
nonrandom functions on $[t, T].$ 
Then, for the 
iterated Stra\-to\-no\-vich stochastic integral of fourth multiplicity

\vspace{-1mm}
\begin{equation}
\label{fin0}
J^{*}[\psi^{(4)}]_{T,t}={\int\limits_t^{*}}^T\psi_4(t_4)
{\int\limits_t^{*}}^{t_4}\psi_3(t_3)
{\int\limits_t^{*}}^{t_3}\psi_2(t_2)
{\int\limits_t^{*}}^{t_2}\psi_1(t_1)
d{\bf w}_{t_1}^{(i_1)}
d{\bf w}_{t_2}^{(i_2)}d{\bf w}_{t_3}^{(i_3)}d{\bf w}_{t_4}^{(i_4)}
\end{equation}

\vspace{3mm}
\noindent
the following 
relations

\vspace{-1mm}
\begin{equation}
\label{fin3}
J^{*}[\psi^{(4)}]_{T,t}
=\hbox{\vtop{\offinterlineskip\halign{
\hfil#\hfil\cr
{\rm l.i.m.}\cr
$\stackrel{}{{}_{p\to \infty}}$\cr
}} }
\sum\limits_{j_1, j_2, j_3,j_4=0}^{p}
C_{j_4j_3 j_2 j_1}\zeta_{j_1}^{(i_1)}\zeta_{j_2}^{(i_2)}\zeta_{j_3}^{(i_3)}\zeta_{j_4}^{(i_4)},
\end{equation}

\vspace{2mm}

\begin{equation}
\label{fin4}
{\sf M}\left\{\left(
J^{*}[\psi^{(4)}]_{T,t}-
\sum\limits_{j_1, j_2, j_3, j_4=0}^{p}
C_{j_4 j_3 j_2 j_1}\zeta_{j_1}^{(i_1)}\zeta_{j_2}^{(i_2)}\zeta_{j_3}^{(i_3)}
\zeta_{j_4}^{(i_4)}
\right)^2\right\}
\le \frac{C}{p^{1-\varepsilon}}
\end{equation}

\vspace{4mm}
\noindent
are fulfilled, where $i_1, \ldots , i_4=0,1,\ldots,m$ in {\rm (\ref{fin0}),} {\rm (\ref{fin3})} 
and $i_1, \ldots, i_4=1,\ldots,m$ in {\rm (\ref{fin4}),}
constant $C$ does not depend on $p,$
$\varepsilon$ is an arbitrary
small positive real number 
for the case of complete orthonormal system of 
Legendre polynomials in the space $L_2([t, T])$
and $\varepsilon=0$ for the case of
complete orthonormal system of 
trigonometric functions in the space $L_2([t, T]),$

\vspace{-1mm}
$$
C_{j_4 j_3 j_2 j_1}=
\int\limits_t^T\psi_4(t_4)\phi_{j_4}(t_4)
\int\limits_t^{t_4}\psi_3(t_3)\phi_{j_3}(t_3)
\int\limits_t^{t_3}\psi_2(t_2)\phi_{j_2}(t_2)
\int\limits_t^{t_2}\psi_1(t_1)\phi_{j_1}(t_1)dt_1dt_2dt_3dt_4;
$$

\vspace{3mm}
\noindent
another notations are the same as in Theorem~{\rm 4}.}

\vspace{2mm}

{\bf Theorem 6}\ \cite{10a}, \cite{12}, \cite{arxiv-11}, \cite{15a}, \cite{new-art-1-xxy}.\
{\it Assume 
that $\{\phi_j(x)\}_{j=0}^{\infty}$ is a complete orthonormal system of 
Legendre polynomials or trigonometric functions in the space $L_2([t, T])$
and $\psi_1(\tau), \ldots,$ $\psi_5(\tau)$ are continuously dif\-ferentiable 
nonrandom functions on $[t, T].$ 
Then, for the 
iterated Stra\-to\-no\-vich stochastic integral of fifth multiplicity

\vspace{-1mm}
\begin{equation}
\label{fin7}
J^{*}[\psi^{(5)}]_{T,t}={\int\limits_t^{*}}^T\psi_5(t_5)
\ldots
{\int\limits_t^{*}}^{t_2}\psi_1(t_1)
d{\bf w}_{t_1}^{(i_1)}
\ldots d{\bf w}_{t_5}^{(i_5)}
\end{equation}

\vspace{3mm}
\noindent
the following 
relations

\vspace{-1mm}
\begin{equation}
\label{fin8}
J^{*}[\psi^{(5)}]_{T,t}
=\hbox{\vtop{\offinterlineskip\halign{
\hfil#\hfil\cr
{\rm l.i.m.}\cr
$\stackrel{}{{}_{p\to \infty}}$\cr
}} }
\sum\limits_{j_1,\ldots,j_5=0}^{p}
C_{j_5 \ldots j_1}\zeta_{j_1}^{(i_1)}\ldots \zeta_{j_5}^{(i_5)},
\end{equation}

\vspace{2mm}

\begin{equation}
\label{fin9}
{\sf M}\left\{\left(
J^{*}[\psi^{(5)}]_{T,t}-
\sum\limits_{j_1, \ldots, j_5=0}^{p}
C_{j_5 \ldots j_1}\zeta_{j_1}^{(i_1)}\ldots
\zeta_{j_5}^{(i_5)}
\right)^2\right\}
\le \frac{C}{p^{1-\varepsilon}}
\end{equation}

\vspace{4mm}
\noindent
are fulfilled, where $i_1, \ldots , i_5=0,1,\ldots,m$ in {\rm (\ref{fin7}),} {\rm (\ref{fin8})} 
and $i_1, \ldots, i_5=1,\ldots,m$ in {\rm (\ref{fin9}),}
constant $C$ is independent of $p,$
$\varepsilon$ is an arbitrary
small positive real number 
for the case of complete orthonormal system of 
Legendre polynomials in the space $L_2([t, T])$
and $\varepsilon=0$ for the case of
complete orthonormal system of 
trigonometric functions in the space $L_2([t, T]),$

\vspace{-1mm}
$$
C_{j_5 \ldots j_1}=
\int\limits_t^T\psi_5(t_5)\phi_{j_5}(t_5)\ldots
\int\limits_t^{t_2}\psi_1(t_1)\phi_{j_1}(t_1)dt_1\ldots dt_5;
$$

\vspace{3mm}
\noindent
another notations are the same as in Theorems~{\rm 4, 5}.}

\vspace{2mm} 

{\bf Theorem 7}\ \cite{10a}, \cite{12}, \cite{arxiv-11}, \cite{15a}, \cite{new-art-1xxys}.\
{\it Suppose that 
$\{\phi_j(x)\}_{j=0}^{\infty}$ is a complete orthonormal system of 
Legendre polynomials or trigonometric functions in the space $L_2([t, T]).$
Then, for the 
iterated Stratonovich stochastic integral of sixth multiplicity

\vspace{-1mm}
\begin{equation}
\label{after10001qu1}
J_{T,t}^{*(i_1\ldots i_6)}={\int\limits_t^{*}}^T
\ldots
{\int\limits_t^{*}}^{t_2}
d{\bf w}_{t_1}^{(i_1)}
\ldots d{\bf w}_{t_6}^{(i_6)}
\end{equation}

\vspace{3mm}
\noindent
the following 
expansion 

\vspace{-2mm}
$$
J_{T,t}^{*(i_1\ldots i_6)}
=\hbox{\vtop{\offinterlineskip\halign{
\hfil#\hfil\cr
{\rm l.i.m.}\cr
$\stackrel{}{{}_{p\to \infty}}$\cr
}} }
\sum\limits_{j_1, \ldots, j_6=0}^{p}
C_{j_6 \ldots j_1}\zeta_{j_1}^{(i_1)}\ldots
\zeta_{j_6}^{(i_6)}
$$

\vspace{3mm}
\noindent
that converges in the mean-square sense is valid, where
$i_1, \ldots, i_6=0, 1,\ldots,m,$

\vspace{-1mm}
$$
C_{j_6 \ldots j_1}=
\int\limits_t^T\phi_{j_6}(t_6)\ldots
\int\limits_t^{t_2}\phi_{j_1}(t_1)dt_1\ldots dt_6
$$

\vspace{3mm}
\noindent
another notations are the same as in Theorems~{\rm 4--6}.}

\vspace{2mm}

Note that an analogue of Theorem 3 for the case of iterated Stratonovich stochastic
integrals of multiplicity 1 follows from (\ref{a1}).

The results of Theorems~4--7 were developed in \cite{10a} (Chapter~2), 
\cite{12}, \cite{arxiv-11}, \cite{15a}, \cite{2024xx}-\cite{2025xxxaaa}.
In particular, analogues of Theorem~7 for iterated Stratonovich stochastic
integrals of multiplicities 7 and 8 were obtained in \cite{10a} (Sect.~2.36, 2.37).
In addition, the variants of Thorems 4--7 
were obtained
for the case when $\{\phi_j(x)\}_{j=0}^{\infty}$ is an arbitrary complete orthonormal system
of functions in $L_2([t, T])$ \cite{10a} (Sect.~2.1.4, 2.23, 2.24, 2.31--2.34),
\cite{12}, \cite{arxiv-11}, \cite{15a}, \cite{2024xx}-\cite{2025xxxaaa}.

\vspace{5mm}

\section{Exact Calculation of the Mean-Square Error in Theorems~1, 2}

\vspace{5mm}

As we mentioned above,
Theorems 1, 2 give new possibilities for exact calculation of 
the mean-square 
error of approximation of iterated Ito stochastic integrals
(see Theorem 8 below).

Assume that $J[\psi^{(k)}]_{T,t}^{p_1 \ldots p_k}$ is the approximation 
of (\ref{ito}), which is
the expression before passing to the limit 
$\hbox{\vtop{\offinterlineskip\halign{
\hfil#\hfil\cr
{\rm l.i.m.}\cr
$\stackrel{}{{}_{p_1,\ldots,p_k\to \infty}}$\cr
}} }
$ on the right-hand side of (\ref{leto6000})
                                                                     
\vspace{1mm}
$$
J[\psi^{(k)}]_{T,t}^{p_1 \ldots p_k}=
\sum\limits_{j_1=0}^{p_1}\ldots
\sum\limits_{j_k=0}^{p_k}
C_{j_k\ldots j_1}\Biggl(
\prod_{l=1}^k\zeta_{j_l}^{(i_l)}+\sum\limits_{r=1}^{[k/2]}
(-1)^r \times
\Biggr.
$$

\vspace{4mm}
$$
\times
\sum_{\stackrel{(\{\{g_1, g_2\}, \ldots, 
\{g_{2r-1}, g_{2r}\}\}, \{q_1, \ldots, q_{k-2r}\})}
{{}_{\{g_1, g_2, \ldots, 
g_{2r-1}, g_{2r}, q_1, \ldots, q_{k-2r}\}=\{1, 2, \ldots, k\}}}}
\prod\limits_{s=1}^r
{\bf 1}_{\{i_{g_{{}_{2s-1}}}=~i_{g_{{}_{2s}}}\ne 0\}}
\Biggl.{\bf 1}_{\{j_{g_{{}_{2s-1}}}=~j_{g_{{}_{2s}}}\}}
\prod_{l=1}^{k-2r}\zeta_{j_{q_l}}^{(i_{q_l})}\Biggr),
$$

\vspace{6mm}
\noindent
where $[x]$ is an integer part of a real number $x;$
another notations are the same as in Theorems~{\rm 1, 2}.

\vspace{4mm} 

Let us denote

\vspace{1mm}

$$
{\sf M}\left\{\left(J[\psi^{(k)}]_{T,t}-
J[\psi^{(k)}]_{T,t}^{p_1,\ldots,p_k}\right)^2\right\}\stackrel{{\rm def}}
{=}E_k^{p_1,\ldots,p_k},
$$

\vspace{2mm}
$$
E_k^{p_1,\ldots,p_k}\stackrel{{\rm def}}{=}E_k^p\ \ \hbox{if}\ \ 
p_1=\ldots=p_k=p,
$$

\vspace{2mm}
$$
\left\Vert K\right\Vert_{L_2([t,T]^k)}^2=\int\limits_{[t,T]^k}
K^2(t_1,\ldots,t_k)dt_1\ldots dt_k\stackrel{{\rm def}}{=}I_k.
$$

\vspace{5mm}

In \cite{7}-\cite{11}, \cite{15b} it was shown that 

\vspace{1mm}

\begin{equation}
\label{star00011}
E_k^{p_1,\ldots,p_k}\le k!\left(I_k-\sum_{j_1=0}^{p_1}\ldots
\sum_{j_k=0}^{p_k}C^2_{j_k\ldots j_1}\right)
\end{equation}

\vspace{5mm}
\noindent
if $i_1,\ldots,i_k=1,\ldots,m$ $(0<T-t<\infty$) or 
$i_1,\ldots,i_k=0, 1,\ldots,m$ $(0<T-t<1).$

Moreover \cite{10a} (Sect.~1.1.9, 1.11, 1.12), \cite{11} (Sect.~6, 15, 16)

\vspace{1mm}
$$
{\sf M}\left\{\left(J[\psi^{(k)}]_{T,t}-
J[\psi^{(k)}]_{T,t}^{p_1,\ldots,p_k}\right)^{2n}\right\}\le C_{n,k}
\left(I_k-\sum_{j_1=0}^{p_1}\ldots
\sum_{j_k=0}^{p_k}C^2_{j_k\ldots j_1}\right)^n,
$$

\vspace{5mm}
\noindent
where $C_{n,k}=(k!)^{n}(2n-1)^{nk}$.

The value $E_k^{p}$
can be calculated exactly.

\vspace{2mm}

{\bf Theorem 8} \cite{10a} (Sect.~1.12), \cite{15b} (Sect.~6).\
{\it Suppose that $\{\phi_j(x)\}_{j=0}^{\infty}$ 
is an arbitrary complete orthonormal system  
of functions in the space $L_2([t,T])$ and
$\psi_1(\tau),\ldots,\psi_k(\tau)\in L_2([t, T]),$  $i_1,\ldots, i_k=1,\ldots,m$.
Then

\begin{equation}
\label{tttr11}
E_k^p=I_k- \sum_{j_1,\ldots, j_k=0}^{p}
C_{j_k\ldots j_1}
{\sf M}\left\{J[\psi^{(k)}]_{T,t}
\sum\limits_{(j_1,\ldots,j_k)}
\int\limits_t^T \phi_{j_k}(t_k)
\ldots
\int\limits_t^{t_{2}}\phi_{j_{1}}(t_{1})
d{\bf f}_{t_1}^{(i_1)}\ldots
d{\bf f}_{t_k}^{(i_k)}\right\},
\end{equation}

\vspace{5mm}
\noindent
where $i_1,\ldots,i_k = 1,\ldots,m;$
the expression 

\vspace{-1mm}
$$
\sum\limits_{(j_1,\ldots,j_k)}
$$ 

\vspace{3mm}
\noindent
means the sum with respect to all
possible permutations 
$(j_1,\ldots,j_k)$. At the same time if 
$j_r$ swapped with $j_q$ in the permutation $(j_1,\ldots,j_k),$
then $i_r$ swapped with $i_q$ in the permutation
$(i_1,\ldots,i_k);$
another notations are the same as in Theorems {\rm 1, 2.}
}

\vspace{2mm}

Note that 

\vspace{-1mm}
$$
{\sf M}\left\{J[\psi^{(k)}]_{T,t}
\int\limits_t^T \phi_{j_k}(t_k)
\ldots
\int\limits_t^{t_{2}}\phi_{j_{1}}(t_{1})
d{\bf f}_{t_1}^{(i_1)}\ldots
d{\bf f}_{t_k}^{(i_k)}\right\}=C_{j_k\ldots j_1}.
$$

\vspace{5mm}

Then from Theorem 8 for pairwise different $i_1,\ldots,i_k$ 
and for $i_1=\ldots=i_k$
we obtain

\vspace{-1mm}
$$
E_k^p= I_k- \sum_{j_1,\ldots,j_k=0}^{p}
C_{j_k\ldots j_1}^2,
$$

\vspace{2mm}
$$ 
E_k^p= I_k - \sum_{j_1,\ldots,j_k=0}^{p}
C_{j_k\ldots j_1}\left(\sum\limits_{(j_1,\ldots,j_k)}
C_{j_k\ldots j_1}\right).
$$

\vspace{6mm}

Consider some examples of application of Theorem 8
$(i_1, i_2 ,i_3=1,\ldots,m)$

\vspace{1mm}
$$
E_2^p
=I_2
-\sum_{j_1,j_2=0}^p
C_{j_2j_1}^2-
\sum_{j_1,j_2=0}^p
C_{j_2j_1}C_{j_1j_2}\ \ \ (i_1=i_2),
$$

\vspace{3mm}
\begin{equation}
\label{881}
E_3^p=I_3
-\sum_{j_3,j_2,j_1=0}^p C_{j_3j_2j_1}^2-
\sum_{j_3,j_2,j_1=0}^p C_{j_3j_1j_2}C_{j_3j_2j_1}\ \ \ (i_1=i_2\ne i_3),
\end{equation}

\vspace{3mm}
\begin{equation}
\label{882}
E_3^p=I_3-
\sum_{j_3,j_2,j_1=0}^p C_{j_3j_2j_1}^2-
\sum_{j_3,j_2,j_1=0}^p C_{j_2j_3j_1}C_{j_3j_2j_1}\ \ \ (i_1\ne i_2=i_3),
\end{equation}

\vspace{3mm}
\begin{equation}
\label{883}
E_3^p=I_3
-\sum_{j_3,j_2,j_1=0}^p C_{j_3j_2j_1}^2-
\sum_{j_3,j_2,j_1=0}^p C_{j_3j_2j_1}C_{j_1j_2j_3}\ \ \ (i_1=i_3\ne i_2).
\end{equation}

\vspace{5mm}

\section{Comparative Analysis of the Efficiency of Application 
of Legendre Polynomials 
and Trigonometric Functions for the Integral 
$J_{(11)T,t}^{(i_1 i_2)}$}

\vspace{5mm}

Using Theorems 1, 2 and complete orthonormal system of 
Legendre polynomials
in the space $L_2([t, T])$ it is shown \cite{3}-\cite{new-art-1xxys} 
(also see
\cite{1}-\cite{2}) that 

\vspace{1mm}
\begin{equation}
\label{4004}
J_{(11)T,t}^{(i_1 i_2)}=
\frac{T-t}{2}\Biggl(\zeta_0^{(i_1)}\zeta_0^{(i_2)}+\sum_{i=1}^{\infty}
\frac{1}{\sqrt{4i^2-1}}\left(
\zeta_{i-1}^{(i_1)}\zeta_{i}^{(i_2)}-
\zeta_i^{(i_1)}\zeta_{i-1}^{(i_2)}\right)-{\bf 1}_{\{i_1=i_2\}}\Biggr),
\end{equation}

\vspace{4mm}
\noindent
where the series converges in the mean-square sense; $i_1, i_2=1,\ldots,m,$

$$
\zeta_{j}^{(i)}=
\int\limits_t^T \phi_{j}(s) d{\bf f}_s^{(i)}
$$ 

\vspace{3mm}
\noindent
are independent standard Gaussian random variables
for various
$i$ or $j$,

\vspace{1mm}
\begin{equation}
\label{4009}
\phi_j(x)=\sqrt{\frac{2j+1}{T-t}}P_j\left(\left(
x-\frac{T+t}{2}\right)\frac{2}{T-t}\right);\ j=0, 1, 2,\ldots,
\end{equation}

\vspace{4mm}
\noindent
where $P_j(x)$ is the Legendre polynomial.

The formula (\ref{4004}) can also be found in \cite{1}-\cite{2}.
It is not difficult to show that \cite{1}-\cite{17}

\vspace{1mm}
\begin{equation}
\label{400}
{\sf M}\left\{\left(J_{(11)T,t}^{(i_1 i_2)}-
J_{(11)T,t}^{(i_1 i_2)q}
\right)^2\right\}
=\frac{(T-t)^2}{2}\Biggl(\frac{1}{2}-\sum_{i=1}^q
\frac{1}{4i^2-1}\Biggr),
\end{equation}

\vspace{3mm}
\noindent
where

\begin{equation}
\label{401}
J_{(11)T,t}^{(i_1 i_2)q}=
\frac{T-t}{2}\Biggl(\zeta_0^{(i_1)}\zeta_0^{(i_2)}+\sum_{i=1}^{q}
\frac{1}{\sqrt{4i^2-1}}\left(
\zeta_{i-1}^{(i_1)}\zeta_{i}^{(i_2)}-
\zeta_i^{(i_1)}\zeta_{i-1}^{(i_2)}\right)-{\bf 1}_{\{i_1=i_2\}}\Biggr).
\end{equation}

\vspace{5mm}

Let us compare (\ref{401}) with (\ref{555}) and (\ref{400}) with (\ref{8010}).
Consider minimal natural numbers $q_{\rm trig}$ and 
$q_{\rm pol},$ which satisfy to (see Table 1)

\vspace{1mm}
$$
\frac{(T-t)^2}{2}\Biggl(\frac{1}{2}-\sum_{i=1}^{q_{\rm pol}}
\frac{1}{4i^2-1}\Biggr)\le (T-t)^3,
$$

\vspace{2mm}
$$
\frac{(T-t)^{2}}{2\pi^2}\Biggl(\frac{\pi^2}{6}-
\sum_{r=1}^{q_{\rm trig}}\frac{1}{r^2}\Biggr)\le (T-t)^3.
$$

\vspace{5mm}

Thus, we have

$$
\frac{q_{\rm pol}}{q_{\rm trig}}\ \ \approx\ \ 1.67,\ \ 2.22,\ \ 2.43,\ \ 2.36,\ \ 2.41,\ 
\ 2.43,\ \ 2.45,\ \ 2.45.
$$

\vspace{4mm}

The formula (\ref{555}) includes $(4q+4)m$ independent
standard Gaussian random variables. At the same time the folmula
(\ref{401}) includes only $(2q+2)m$ independent
standard Gaussian random variables. Moreover, the formula
(\ref{401}) is simpler than the formula (\ref{555}).
Thus, in this case we can talk about approximately equal computational costs
for the formulas (\ref{555}) and (\ref{401}).

\begin{table}
\centering
\caption{Numbers $q_{\rm trig},$ $q_{\rm trig}^{*}$, $q_{\rm pol}$}
\label{tab:1}      
\begin{tabular}{p{1.1cm}p{1.1cm}p{1.1cm}p{1.1cm}p{1.1cm}p{1.1cm}p{1.1cm}p{1.1cm}p{1.1cm}}
\hline\noalign{\smallskip}
$T-t$&$2^{-5}$&$2^{-6}$&$2^{-7}$&$2^{-8}$&$2^{-9}$&$2^{-10}$&$2^{-11}$&$2^{-12}$\\
\noalign{\smallskip}\hline\noalign{\smallskip}
$q_{\rm trig}$&3&4&7&14&27&53&105&209\\
\noalign{\smallskip}\hline\noalign{\smallskip}
$q_{\rm trig}^{*}$&6&11&20&40&79&157&312&624\\
\noalign{\smallskip}\hline\noalign{\smallskip}
$q_{\rm pol}$&5&9&17&33&65&129&257&513\\
\noalign{\smallskip}\hline\noalign{\smallskip}
\end{tabular}
\end{table}

There is one important feature. 
As we mentioned above, further we will see that introduction of random 
variables $\xi_q^{(i)}$ and 
$\mu_q^{(i)}$ will sharply 
complicate the approximation of 
the iterated stochastic integral $J_{(111)T,t}^{(i_1 i_2 i_3)};$
$i_1,i_2,i_3=1,\ldots,m.$ 
This is due to the fact that
the number $q$ is fixed for all stochastic integrals, which 
included into the considered collection. However, it is clear that due 
to the smallness of $T-t$, the number $q$ for $J_{(111)T,t}^{(i_1 i_2 i_3)}$
could be chosen significantly 
less than in the formula (\ref{555}). 
This feature is also valid for the formulas (\ref{444}), (\ref{1970}).
However, for the case of Legendre polynomials we can choose different 
numbers $q$ for different iterated stochastic integrals.

From the other hand, if we will not introduce the random 
variables $\xi_q^{(i)}$ and 
$\mu_q^{(i)},$ then the mean-square error of approximation of the
iterated 
stochastic integral $J_{(11)T,t}^{(i_1 i_2)}$ will be three times larger
(see (\ref{801})). 
Moreover, in this case the stochastic integrals 
$J_{(01)T,t}^{(0 i_1)}$, $J_{(001)T,t}^{(00i_1)}$ 
(with Gaussian distribution)
will be approximated worse.

Consider minimal natural numbers $q_{\rm trig}^{*}$, 
which satisfy to (see Table 1)

\vspace{1mm}
$$
\frac{3(T-t)^{2}}{2\pi^2}\Biggl(\frac{\pi^2}{6}-
\sum_{r=1}^{q_{\rm trig}^{*}}\frac{1}{r^2}\Biggr)\le (T-t)^3.
$$

\vspace{4mm}

In this situation we can talk about
the advantage of Legendre polynomials ($q_{\rm trig}^{*} > q_{pol}$ and
(\ref{555}) is more complex than (\ref{401})).

\vspace{5mm}

\section{Comparative Analysis of the Efficiency 
of Application of Legendre Polynomials 
and Trigonometric Functions for the Integrals $J_{(1)T,t}^{(i_1)},$
$J_{(11)T,t}^{(i_1 i_2)},$ $J_{(01)T,t}^{(0 i_1)},$ $J_{(10)T,t}^{(i_1 0)},$
$J_{(111)T,t}^{(i_1 i_2 i_3)}$}

\vspace{5mm}

It is well known \cite{19}-\cite{23} that for implementation
of
strong Taylor--Ito numerical methods with the 
order 1.5 of accuracy for Ito stochastic differential
equations we need to
approximate the following collection of iterated Ito stochastic integrals

\vspace{1mm}
$$
J_{(1)T,t}^{(i_1)},\ \ \
J_{(11)T,t}^{(i_1 i_2)},\ \ \ 
J_{(01)T,t}^{(0 i_1)},\ \ \ J_{(10)T,t}^{(i_1 0)},\ \ \ 
J_{(111)T,t}^{(i_1 i_2 i_3)}.
$$

\vspace{4mm}

Using Theorems 1, 2 for the system of trigonometric 
functions, we have \cite{3}-\cite{17} (also see \cite{1}-\cite{2})

\vspace{1mm}
\begin{equation}
\label{x1}
J_{(1)T,t}^{(i_1)}=\sqrt{T-t}\zeta_0^{(i_1)},
\end{equation}

\vspace{2mm}
$$
J_{(11)T,t}^{(i_1 i_2)q}=\frac{1}{2}(T-t)\Biggl(
\zeta_{0}^{(i_1)}\zeta_{0}^{(i_2)}
+\frac{1}{\pi}
\sum_{r=1}^{q}\frac{1}{r}\left(
\zeta_{2r}^{(i_1)}\zeta_{2r-1}^{(i_2)}-
\zeta_{2r-1}^{(i_1)}\zeta_{2r}^{(i_2)}+
\right.\Biggr.
$$

\vspace{1mm}
$$
+\left.\sqrt{2}\left(\zeta_{2r-1}^{(i_1)}\zeta_{0}^{(i_2)}-
\zeta_{0}^{(i_1)}\zeta_{2r-1}^{(i_2)}\right)\right)
+
$$
\begin{equation}
\label{x2}
\Biggl.+\frac{\sqrt{2}}{\pi}\sqrt{\alpha_q}\left(
\xi_q^{(i_1)}\zeta_0^{(i_2)}-\zeta_0^{(i_1)}\xi_q^{(i_2)}\right)-
{\bf 1}_{\{i_1=i_2\}}\Biggr),
\end{equation}

\vspace{5mm}

\begin{equation}
\label{x3}
J_{(01)T,t}^{(0 i_1)q}=\frac{{(T-t)}^{3/2}}{2}
\Biggl(\zeta_0^{(i_1)}-\frac{\sqrt{2}}{\pi}\left(\sum_{r=1}^{q}
\frac{1}{r}
\zeta_{2r-1}^{(i_1)}+\sqrt{\alpha_q}\xi_q^{(i_1)}\right)
\Biggr),
\end{equation}

\vspace{5mm}

\begin{equation}
\label{x4}
J_{(10)T,t}^{(i_1 0)q}=\frac{{(T-t)}^{3/2}}{2}
\Biggl(\zeta_0^{(i_1)}+\frac{\sqrt{2}}{\pi}\left(\sum_{r=1}^{q}
\frac{1}{r}
\zeta_{2r-1}^{(i_1)}+\sqrt{\alpha_q}\xi_q^{(i_1)}\right)
\Biggr),
\end{equation}

\vspace{7mm}

$$
J_{(111)T,t}^{(i_1 i_2 i_3)q}=(T-t)^{3/2}\Biggl(\frac{1}{6}
\zeta_{0}^{(i_1)}\zeta_{0}^{(i_2)}\zeta_{0}^{(i_3)}+\Biggr.
\frac{\sqrt{\alpha_q}}{2\sqrt{2}\pi}\left(
\xi_q^{(i_1)}\zeta_0^{(i_2)}\zeta_0^{(i_3)}-\xi_q^{(i_3)}\zeta_0^{(i_2)}
\zeta_0^{(i_1)}\right)+
$$
$$
+\frac{1}{2\sqrt{2}\pi^2}\sqrt{\beta_q}\left(
\mu_q^{(i_1)}\zeta_0^{(i_2)}\zeta_0^{(i_3)}-2\mu_q^{(i_2)}\zeta_0^{(i_1)}
\zeta_0^{(i_3)}+\mu_q^{(i_3)}\zeta_0^{(i_1)}\zeta_0^{(i_2)}\right)+
$$
$$
+
\frac{1}{2\sqrt{2}}\sum_{r=1}^{q}
\Biggl(\frac{1}{\pi r}\left(
\zeta_{2r-1}^{(i_1)}
\zeta_{0}^{(i_2)}\zeta_{0}^{(i_3)}-
\zeta_{2r-1}^{(i_3)}
\zeta_{0}^{(i_2)}\zeta_{0}^{(i_1)}\right)+\Biggr.
$$
$$
\Biggl.+
\frac{1}{\pi^2 r^2}\left(
\zeta_{2r}^{(i_1)}
\zeta_{0}^{(i_2)}\zeta_{0}^{(i_3)}-
2\zeta_{2r}^{(i_2)}
\zeta_{0}^{(i_3)}\zeta_{0}^{(i_1)}+
\zeta_{2r}^{(i_3)}
\zeta_{0}^{(i_2)}\zeta_{0}^{(i_1)}\right)\Biggr)+
$$
$$
+
\sum_{r=1}^{q}
\Biggl(\frac{1}{4\pi r}\left(
\zeta_{2r}^{(i_1)}
\zeta_{2r-1}^{(i_2)}\zeta_{0}^{(i_3)}-
\zeta_{2r-1}^{(i_1)}
\zeta_{2r}^{(i_2)}\zeta_{0}^{(i_3)}-
\zeta_{2r-1}^{(i_2)}
\zeta_{2r}^{(i_3)}\zeta_{0}^{(i_1)}+
\zeta_{2r-1}^{(i_3)}
\zeta_{2r}^{(i_2)}\zeta_{0}^{(i_1)}\right)+\Biggr.
$$
$$
+
\frac{1}{8\pi^2 r^2}\left(
3\zeta_{2r-1}^{(i_1)}
\zeta_{2r-1}^{(i_2)}\zeta_{0}^{(i_3)}+
\zeta_{2r}^{(i_1)}
\zeta_{2r}^{(i_2)}\zeta_{0}^{(i_3)}-
6\zeta_{2r-1}^{(i_1)}
\zeta_{2r-1}^{(i_3)}\zeta_{0}^{(i_2)}+\right.
$$
\begin{equation}
\label{x5}
\Biggl.\left.
+
3\zeta_{2r-1}^{(i_2)}
\zeta_{2r-1}^{(i_3)}\zeta_{0}^{(i_1)}-
2\zeta_{2r}^{(i_1)}
\zeta_{2r}^{(i_3)}\zeta_{0}^{(i_2)}+
\zeta_{2r}^{(i_3)}
\zeta_{2r}^{(i_2)}\zeta_{0}^{(i_1)}\right)\Biggr)
\Biggl.+D_{T,t}^{(i_1i_2i_3)q}\Biggr),
\end{equation}

\vspace{6mm}
\noindent
where

$$
D_{T,t}^{(i_1i_2i_3)q}=
\frac{1}{2\pi^2}\sum_{\stackrel{r,l=1}{{}_{r\ne l}}}^{q}
\Biggl(\frac{1}{r^2-l^2}\biggl(
\zeta_{2r}^{(i_1)}
\zeta_{2l}^{(i_2)}\zeta_{0}^{(i_3)}-
\zeta_{2r}^{(i_2)}
\zeta_{0}^{(i_1)}\zeta_{2l}^{(i_3)}+\biggr.\Biggr.
$$
$$
\Biggl.+\biggl.
\frac{r}{l}
\zeta_{2r-1}^{(i_1)}
\zeta_{2l-1}^{(i_2)}\zeta_{0}^{(i_3)}-\frac{l}{r}
\zeta_{0}^{(i_1)}
\zeta_{2r-1}^{(i_2)}\zeta_{2l-1}^{(i_3)}\biggr)-
\frac{1}{rl}\zeta_{2r-1}^{(i_1)}
\zeta_{0}^{(i_2)}\zeta_{2l-1}^{(i_3)}\Biggr)+
$$
$$
+
\frac{1}{4\sqrt{2}\pi^2}\Biggl(
\sum_{r,m=1}^{q}\Biggl(\frac{2}{rm}
\left(-\zeta_{2r-1}^{(i_1)}
\zeta_{2m-1}^{(i_2)}\zeta_{2m}^{(i_3)}+
\zeta_{2r-1}^{(i_1)}
\zeta_{2r}^{(i_2)}\zeta_{2m-1}^{(i_3)}+
\right.\Biggr.\Biggr.
$$
$$
\left.+
\zeta_{2r-1}^{(i_1)}
\zeta_{2m}^{(i_2)}\zeta_{2m-1}^{(i_3)}-
\zeta_{2r}^{(i_1)}
\zeta_{2r-1}^{(i_2)}\zeta_{2m-1}^{(i_3)}\right)+
$$
$$
+\frac{1}{m(r+m)}
\left(-\zeta_{2(m+r)}^{(i_1)}
\zeta_{2r}^{(i_2)}\zeta_{2m}^{(i_3)}-
\zeta_{2(m+r)-1}^{(i_1)}
\zeta_{2r-1}^{(i_2)}\zeta_{2m}^{(i_3)}-
\right.
$$
$$
\Biggl.\left.
-\zeta_{2(m+r)-1}^{(i_1)}
\zeta_{2r}^{(i_2)}\zeta_{2m-1}^{(i_3)}+
\zeta_{2(m+r)}^{(i_1)}
\zeta_{2r-1}^{(i_2)}\zeta_{2m-1}^{(i_3)}\right)\Biggr)+
$$
$$
+
\sum_{m=1}^{q}\sum_{l=m+1}^{q}\Biggl(\frac{1}{m(l-m)}
\left(\zeta_{2(l-m)}^{(i_1)}
\zeta_{2l}^{(i_2)}\zeta_{2m}^{(i_3)}+
\zeta_{2(l-m)-1}^{(i_1)}
\zeta_{2l-1}^{(i_2)}\zeta_{2m}^{(i_3)}-
\right.\Biggr.
$$
$$
\left.
-\zeta_{2(l-m)-1}^{(i_1)}
\zeta_{2l}^{(i_2)}\zeta_{2m-1}^{(i_3)}+
\zeta_{2(l-m)}^{(i_1)}
\zeta_{2l-1}^{(i_2)}\zeta_{2m-1}^{(i_3)}\right)+
$$
$$
+
\frac{1}{l(l-m)}
\left(-\zeta_{2(l-m)}^{(i_1)}
\zeta_{2m}^{(i_2)}\zeta_{2l}^{(i_3)}+
\zeta_{2(l-m)-1}^{(i_1)}
\zeta_{2m-1}^{(i_2)}\zeta_{2l}^{(i_3)}-
\right.
$$
$$
\Biggl.
\Biggl.
\Biggl.
\left.
-\zeta_{2(l-m)-1}^{(i_1)}
\zeta_{2m}^{(i_2)}\zeta_{2l-1}^{(i_3)}-
\zeta_{2(l-m)}^{(i_1)}
\zeta_{2m-1}^{(i_2)}\zeta_{2l-1}^{(i_3)}\right)\Biggr)\Biggr),
$$

\vspace{5mm}
\noindent
where

\vspace{1mm}
$$
\xi_q^{(i)}=\frac{1}{\sqrt{\alpha_q}}\sum_{r=q+1}^{\infty}
\frac{1}{r}~\zeta_{2r-1}^{(i)},\ \ \
\alpha_q=\frac{\pi^2}{6}-\sum_{r=1}^q\frac{1}{r^2},
$$

\vspace{1mm}
$$
\mu_q^{(i)}=\frac{1}{\sqrt{\beta_q}}\sum_{r=q+1}^{\infty}
\frac{1}{r^2}~\zeta_{2r}^{(i)},\ \ \
\beta_q=\frac{\pi^4}{90}-\sum_{r=1}^q\frac{1}{r^4},\
$$

\vspace{6mm}
\noindent
and $\zeta_0^{(i)},$ $\zeta_{2r}^{(i)},$
$\zeta_{2r-1}^{(i)},$ $\xi_q^{(i)},$ $\mu_q^{(i)}$ 
($r=1,\ldots,q;$
$i=1,\ldots,m$) are independent
standard Gaussian random variables. Moreover, in (\ref{x5})
we suppose that $i_1\ne i_2,$ $i_1\ne i_3,$ $i_2\ne i_3$.

Mean-square errors for the approximations (\ref{x2})--(\ref{x5})
are represented by the formulas

\vspace{1mm}
$$
{\sf M}\left\{\left(J_{(01)T,t}^{(0 i_1)}-
J_{(01)T,t}^{(0 i_1)q}
\right)^2\right\}=0,
$$

\vspace{2mm}
$$
{\sf M}\left\{\left(J_{(10)T,t}^{(i_1 0)}-
J_{(10)T,t}^{(i_1 0)q}
\right)^2\right\}=0,
$$

\vspace{2mm}
\begin{equation}
\label{y1}
{\sf M}\left\{\left(J_{(11)T,t}^{(i_1 i_2)}-
J_{(11)T,t}^{(i_1 i_2)q}
\right)^2\right\}
=\frac{(T-t)^{2}}{2\pi^2}\Biggl(\frac{\pi^2}{6}-
\sum_{r=1}^q \frac{1}{r^2}\Biggr),
\end{equation}

\vspace{4mm}
$$
{\sf M}\left\{\left(J_{(111)T,t}^{(i_1 i_2 i_3)}-
J_{(111)T,t}^{(i_1 i_2 i_3)q}\right)^2\right\}=
(T-t)^3\Biggl(\frac{4}{45}-\frac{1}{4\pi^2}\sum_{r=1}^q\frac{1}{r^2}-
\Biggl.
$$

\vspace{1mm}
\begin{equation}
\label{x7}
\Biggl.-\frac{55}{32\pi^4}\sum_{r=1}^q\frac{1}{r^4}-
\frac{1}{4\pi^4}\sum_{\stackrel{r,l=1}{{}_{r\ne l}}}^q
\frac{5l^4+4r^4-3r^2l^2}{r^2 l^2 \left(r^2-l^2\right)^2}\Biggr),
\end{equation}

\vspace{6mm}
\noindent
where $i_1\ne i_2,\ i_1\ne i_3,\ i_2\ne i_3$.

\begin{table}
\centering
\caption{Confirmation of the formula (\ref{x7})}
\label{tab:2}      
\begin{tabular}{p{2.1cm}p{1.7cm}p{1.7cm}p{2.1cm}p{2.3cm}p{2.3cm}p{2.3cm}}
\hline\noalign{\smallskip}
$\varepsilon/(T-t)^3$&0.0459&0.0072&$7.5722\cdot 10^{-4}$
&$7.5973\cdot 10^{-5}$&
$7.5990\cdot 10^{-6}$\\
\noalign{\smallskip}\hline\noalign{\smallskip}
$q$&1&10&100&1000&10000\\
\noalign{\smallskip}\hline\noalign{\smallskip}
\end{tabular}
\end{table}

In Table 2, we can see the numerical confirmation of 
the formula (\ref{x7}) ($\varepsilon$ is a right-hand side of (\ref{x7})).

Note that the formulas (\ref{x1}), (\ref{x2}) have been obtained 
for the first time in \cite{19}. Using 
(\ref{x1}), (\ref{x2}), we can realize numerically 
the explicit one-step strong Taylor--Ito numerical 
method with the order 1.0 of accuracy 
(Milstein scheme \cite{19}).
The analogue of the formula (\ref{x5}) has been obtained 
for the first time in 
\cite{20}-\cite{22}.                         

As we mentioned above, the Milstein approach (see Sect.~2) leads to iterated 
application of the operation
of limit transition. The analogue of (\ref{x5}) has been derived in 
\cite{20}-\cite{22}, \cite{Zapad-9}
on the base of the Milstein approach \cite{19}.
It means that the authors of the works \cite{20}-\cite{22},
\cite{Zapad-9} 
could not formally use
the double sum with the upper limit $q$ in the analogue of
(\ref{x5}) in 
\cite{20} (pp.~438-439),
\cite{21}
(Sect.~5.8, pp.~202--204), \cite{22} (pp.~82-84),
\cite{Zapad-9} (pp.~263-264)
on the base of the Wong--Zakai approximation
\cite{W-Z-1}-\cite{Watanabe}
(see discussion in Sect.~11 for details).
From the other hand, the correctness of (\ref{x5}) follows 
directly from Theorems 1, 2. Note that (\ref{x5}) has been
obtained  reasonably
for the first time in \cite{3}. The version of (\ref{x5}) but without
using the random variables $\xi_q^{(i)}$ and $\mu_q^{(i)}$ 
can be found in \cite{1}-\cite{2}.

The formula (\ref{y1}) appears for the first time in \cite{19}.
The mean-square error (\ref{x7}) has been obtained for the first time
in \cite{3} on the base of the simplified variant of 
Theorem 8 (the case of pairwise different $i_1,\ldots,i_k$).

As we noted above, the number $q$  must be the same
in (\ref{x2})--(\ref{x5}). This is the main drawback of this approach
because really the number $q$ in (\ref{x5})
can be chosen essentially smaller than in (\ref{x2}).

Note that in (\ref{x5}) we can replace 
$J_{(111)T,t}^{(i_1 i_2 i_3)q}$ with $J_{(111)T,t}^{*(i_1 i_2 i_3)q}$
and (\ref{x5}) will when be valid for any $i_1, i_2, i_3 = 0, 1,\ldots,m$
(see Theorem 3).

Consider approximations of the iterated Ito stochastic integrals

$$
J_{(1)T,t}^{(i_1)},\ \ \ 
J_{(11)T,t}^{(i_1 i_2)},\ \ \ 
J_{(01)T,t}^{(0 i_1)},\ \ \ J_{(10)T,t}^{(i_1 0)},\ \ \ 
J_{(111)T,t}^{(i_1 i_2 i_3)}\ \ \ (i_1,i_2,i_3=1,\ldots,m)
$$

\vspace{4mm}
\noindent
on the base of Theorems 1, 2 (the case of Legendre 
polynomials) \cite{1}-\cite{17}

\vspace{1mm}
\begin{equation}
\label{y2}
J_{(1)T,t}^{(i_1)}=\sqrt{T-t}\zeta_0^{(i_1)},
\end{equation}

\vspace{1mm}
\begin{equation}
\label{y3}
J_{(11)T,t}^{(i_1 i_2)q}=
\frac{T-t}{2}\Biggl(\zeta_0^{(i_1)}\zeta_0^{(i_2)}+\sum_{i=1}^{q}
\frac{1}{\sqrt{4i^2-1}}\left(
\zeta_{i-1}^{(i_1)}\zeta_{i}^{(i_2)}-
\zeta_i^{(i_1)}\zeta_{i-1}^{(i_2)}\right)-{\bf 1}_{\{i_1=i_2\}}\Biggr),
\end{equation}

\vspace{2mm}
\begin{equation}
\label{y3a}
J_{(01)T,t}^{(0 i_1)}=\frac{(T-t)^{3/2}}{2}\Biggl(\zeta_0^{(i_1)}+
\frac{1}{\sqrt{3}}\zeta_1^{(i_1)}\Biggr),
\end{equation}

\vspace{2mm}
\begin{equation}
\label{y4a}
J_{(10)T,t}^{(i_1 0)}=\frac{(T-t)^{3/2}}{2}\Biggl(\zeta_0^{(i_1)}-
\frac{1}{\sqrt{3}}\zeta_1^{(i_1)}\Biggr),
\end{equation}

\vspace{4mm}
$$
J_{(111)T,t}^{(i_1 i_2 i_3)q_1}=
\sum_{j_1,j_2,j_3=0}^{q_1}
C_{j_3j_2j_1}
\Biggl(
\zeta_{j_1}^{(i_1)}\zeta_{j_2}^{(i_2)}\zeta_{j_3}^{(i_3)}
-{\bf 1}_{\{i_1=i_2\}}
{\bf 1}_{\{j_1=j_2\}}
\zeta_{j_3}^{(i_3)}-\Biggr.
$$

\vspace{1mm}
\begin{equation}
\label{y4}
\Biggl.
-{\bf 1}_{\{i_2=i_3\}}
{\bf 1}_{\{j_2=j_3\}}
\zeta_{j_1}^{(i_1)}-
{\bf 1}_{\{i_1=i_3\}}
{\bf 1}_{\{j_1=j_3\}}
\zeta_{j_2}^{(i_2)}\Biggr),\ \ \ q_1\ll q,
\end{equation}

\vspace{5mm}

$$
J_{(111)T,t}^{(i_1 i_1 i_1)}
=\frac{1}{6}(T-t)^{3/2}\left(
\left(\zeta_0^{(i_1)}\right)^3-3
\zeta_0^{(i_1)}\right),
$$

\vspace{5mm}

$$
C_{j_3j_2j_1}=\int\limits_{t}^{T}\phi_{j_3}(z)
\int\limits_{t}^{z} \phi_{j_2}(y)
\int\limits_{t}^{y}
\phi_{j_1}(x) dx dy dz=
$$

\vspace{1mm}
$$
=\frac{\sqrt{(2j_1+1)(2j_2+1)(2j_3+1)}}{8}(T-t)^{3/2}\bar
C_{j_3j_2j_1},
$$

\vspace{4mm}

\begin{equation}
\label{ww}
\bar C_{j_3j_2j_1}=\int\limits_{-1}^{1}P_{j_3}(z)
\int\limits_{-1}^{z}P_{j_2}(y)
\int\limits_{-1}^{y}
P_{j_1}(x)dx dy dz,
\end{equation}

\vspace{6mm}
\noindent
where $\phi_{j}(x)$ has the form (\ref{4009}) 
and $P_i(x)$ is the Legendre polynomial
$(i=0, 1, 2,\ldots)$.

Mean-square errors for the approximations (\ref{y3}), (\ref{y4})
are represented by the 
formulas (see Theo\-rem 8 and (\ref{star00011})) \cite{1}-\cite{17}

\vspace{1mm}
\begin{equation}
\label{sad004}
{\sf M}\left\{\left(J_{(11)T,t}^{(i_1 i_2)}-
J_{(11)T,t}^{(i_1 i_2)q}
\right)^2\right\}
=\frac{(T-t)^2}{2}\Biggl(\frac{1}{2}-\sum_{i=1}^q
\frac{1}{4i^2-1}\Biggr)\ \ \ (i_1\ne i_2),
\end{equation}

\vspace{3mm}

\begin{equation}
\label{sad005}
{\sf M}\left\{\left(
J_{(111)T,t}^{(i_1 i_2 i_3)}-
J_{(111)T,t}^{(i_1 i_2 i_3)q_1}\right)^2\right\}=
\frac{(T-t)^{3}}{6}-
\sum_{j_3,j_2,j_1=0}^{q_1}
C_{j_3j_2j_1}^2\ \ \ (i_1\ne i_2, i_1\ne i_3, i_2\ne i_3),
\end{equation}

\vspace{7mm}
$$
{\sf M}\left\{\left(
J_{(111)T,t}^{(i_1 i_2 i_3)}-
J_{(111)T,t}^{(i_1 i_2 i_3)q_1}\right)^2\right\}=
\frac{(T-t)^{3}}{6}-\sum_{j_3,j_2,j_1=0}^{q_1}
C_{j_3j_2j_1}^2-
$$

\vspace{1mm}
\begin{equation}
\label{sad006}
-\sum_{j_3,j_2,j_1=0}^{q_1}
C_{j_2j_3j_1}C_{j_3j_2j_1}\ \ \ (i_1\ne i_2=i_3),
\end{equation}

\vspace{6mm}
$$
{\sf M}\left\{\left(
J_{(111)T,t}^{(i_1 i_2 i_3)}-
J_{(111)T,t}^{(i_1 i_2 i_3)q_1}\right)^2\right\}=
\frac{(T-t)^{3}}{6}
-\sum_{j_3,j_2,j_1=0}^{q_1}
C_{j_3j_2j_1}^2-
$$

\vspace{1mm}
\begin{equation}
\label{sad007}
-\sum_{j_3,j_2,j_1=0}^{q_1}
C_{j_3j_2j_1}C_{j_1j_2j_3}\ \ \ (i_1=i_3\ne i_2),
\end{equation}

\vspace{6mm}
$$
{\sf M}\left\{\left(
J_{(111)T,t}^{(i_1 i_2 i_3)}-
J_{(111)T,t}^{(i_1 i_2 i_3)q_1}\right)^2\right\}=
\frac{(T-t)^{3}}{6}-\sum_{j_3,j_2,j_1=0}^{q_1}
C_{j_3j_2j_1}^2-
$$

\vspace{1mm}
\begin{equation}
\label{sad008}
-\sum_{j_3,j_2,j_1=0}^{q_1}
C_{j_3j_1j_2}C_{j_3j_2j_1}\ \ \ (i_1=i_2\ne i_3),
\end{equation}

\vspace{4mm}

\begin{equation}
\label{sad009}
{\sf M}\left\{\left(
J_{(111)T,t}^{(i_1 i_2 i_3)}-
J_{(111)T,t}^{(i_1 i_2 i_3)q_1}\right)^2\right\}\le
6\left(\frac{(T-t)^{3}}{6}-\sum_{j_3,j_2,j_1=0}^{q_1}
C_{j_3j_2j_1}^2\right)\ \ \ (i_1,i_2,i_3=1,\ldots,m).
\end{equation}

\vspace{6mm}

Let us compare 
the efficiency of application of Legendre polynomials 
and trigonometric functions for the iterated stochastic integrals
$J_{(11)T,t}^{(i_1 i_2)},$ $J_{(111)T,t}^{(i_1 i_2 i_3)}$.

Consider the following conditions $(i_1\ne i_2,\ 
i_1\ne i_3,\ i_2\ne i_3)$

\vspace{1mm}
\begin{equation}
\label{z1}
\frac{(T-t)^2}{2}\Biggl(\frac{1}{2}-\sum_{i=1}^{q}
\frac{1}{4i^2-1}\Biggr)\le (T-t)^4,
\end{equation}

\vspace{2mm}
\begin{equation}
\label{z2}
(T-t)^{3}\Biggl(\frac{1}{6}-\sum_{j_1,j_2,j_3=0}^{q_1}
\frac{\left(C_{j_3j_2j_1}\right)^2}{(T-t)^3}\Biggr)\le (T-t)^4,
\end{equation}

\vspace{2mm}
\begin{equation}
\label{z3}
\frac{(T-t)^{2}}{2\pi^2}\Biggl(\frac{\pi^2}{6}-
\sum_{r=1}^{p}\frac{1}{r^2}\Biggr)\le (T-t)^4,
\end{equation}

\vspace{2mm}
\begin{equation}
\label{z4}
(T-t)^3\Biggl(\frac{4}{45}-\frac{1}{4\pi^2}\sum_{r=1}^{p_1}\frac{1}{r^2}
\Biggl.
\Biggl.-\frac{55}{32\pi^4}\sum_{r=1}^{p_1}\frac{1}{r^4}-
\frac{1}{4\pi^4}\sum_{\stackrel{r,l=1}{{}_{r\ne l}}}^{p_1}
\frac{5l^4+4r^4-3r^2l^2}{r^2 l^2 \left(r^2-l^2\right)^2}\Biggr)\le (T-t)^4,
\end{equation}

\vspace{4mm}
\noindent
where

\vspace{-1mm}
$$
C_{j_3j_2j_1}=\frac{\sqrt{(2j_1+1)(2j_2+1)(2j_3+1)}}{8}(T-t)^{3/2}\bar
C_{j_3j_2j_1},
$$

\vspace{2mm}
$$
\bar C_{j_3j_2j_1}=\int\limits_{-1}^{1}P_{j_3}(z)
\int\limits_{-1}^{z}P_{j_2}(y)
\int\limits_{-1}^{y}
P_{j_1}(x)dx dy dz,
$$

\vspace{5mm}
\noindent
where $P_i(x)$ is the Legendre polynomial.

In Tables 3 and 4, we can see minimal numbers 
$q,$ $q_1,$ $p,$ $p_1,$ which satisfy the conditions
(\ref{z1})--(\ref{z4}). As we mentioned above, the
numbers $q,$ $q_1$ are different. At that $q_1\ll q$
(the case of Legendre polynomials). As we saw in the previous
sections, we cannot take different numbers 
$p,$ $p_1$ for the case of trigonometric functions. Thus, we have to
choose $q=p$ in (\ref{x2})--(\ref{x5}). This leads
to huge computational costs (see very complex formula (\ref{x5})).
From the other hand, we can choose different numbers $q$
in (\ref{x2})--(\ref{x5}). At that we must exclude
the random variables $\xi_q^{(i)},$ $\mu_q^{(i)}$ from 
(\ref{x2})--(\ref{x5}).

At this situation for the case $i_1\ne i_2,$ $i_2\ne i_3,$ 
$i_1\ne i_3$ we have

\vspace{1mm}
\begin{equation}
\label{zzz3}
\frac{3(T-t)^{2}}{2\pi^2}\Biggl(\frac{\pi^2}{6}-
\sum_{r=1}^{p^{*}}\frac{1}{r^2}\Biggr)\le (T-t)^4,
\end{equation}

\vspace{2mm}
\begin{equation}
\label{zzz4}
(T-t)^3\Biggl(\frac{5}{36}-\frac{1}{2\pi^2}\sum_{r=1}^{p_1^{*}}\frac{1}{r^2}
\Biggl.
\Biggl.-\frac{79}{32\pi^4}\sum_{r=1}^{p_1^{*}}\frac{1}{r^4}-
\frac{1}{4\pi^4}\sum_{\stackrel{r,l=1}{{}_{r\ne l}}}^{p_1^{*}}
\frac{5l^4+4r^4-3r^2l^2}{r^2 l^2 \left(r^2-l^2\right)^2}\Biggr)\le (T-t)^4,
\end{equation}

\vspace{5mm}
\noindent
where the left-hand sides of (\ref{zzz3}), (\ref{zzz4}) correspond
to (\ref{555}), (\ref{x5}) but without $\xi_q^{(i)},$ $\mu_q^{(i)}$.
In Table 4, we can see minimal numbers 
$p^{*},$ $p_1^{*}$, which satisfy the conditions
(\ref{zzz3}), (\ref{zzz4}). 

Moreover,

$$
{\sf M}\left\{\left(J_{(01)T,t}^{(0 i_1)}-
J_{(01)T,t}^{(0 i_1)q}
\right)^2\right\}=
{\sf M}\left\{\left(J_{(10)T,t}^{(i_1 0)}-
J_{(10)T,t}^{(i_1 0)q}
\right)^2\right\}=
$$

\begin{equation}
\label{1001x}
=\frac{(T-t)^{3}}{2\pi^2}\Biggl(\frac{\pi^2}{6}-
\sum_{r=1}^q \frac{1}{r^2}\Biggr)\ne 0,
\end{equation}

\vspace{5mm}
\noindent
where
$J_{(01)T,t}^{(0 i_1)q}$, $J_{(10)T,t}^{(i_1 0)q}$ 
are defined by the formulas (\ref{x3}), (\ref{x4}).

It is not difficult to see that the numbers $q_{\rm trig}$ in Table 1
correspond to minimal numbers $q_{\rm trig}$, which satisfy 
the condition

\vspace{1mm}
$$
\frac{(T-t)^{3}}{2\pi^2}\Biggl(\frac{\pi^2}{6}-
\sum_{r=1}^{q_{\rm trig}} \frac{1}{r^2}\Biggr)\le (T-t)^4.
$$

\vspace{4mm}

From the other hand, the right-hand sides of (\ref{y3a}), (\ref{y4a}) 
include only two random variables.
In this situation we can again talk about
the advantage of Ledendre polynomials.

\begin{table}
\centering
\caption{Numbers $q,$ $q_1$}
\label{tab:3}      
\begin{tabular}{p{1.1cm}p{2.1cm}p{2.1cm}p{2.1cm}p{2.1cm}p{2.1cm}p{2.1cm}}
\hline\noalign{\smallskip}
$T-t$&$0.08222$&$0.05020$&$0.02310$&$0.01956$\\
\noalign{\smallskip}\hline\noalign{\smallskip}
$q$&19&51&235&328\\
\noalign{\smallskip}\hline\noalign{\smallskip}
$q_1$&1&2&5&6\\
\noalign{\smallskip}\hline\noalign{\smallskip}
\end{tabular}
\end{table}

\begin{table}
\centering
\caption{Numbers $p,$ $p_1,$ $p^{*},$ $p_1^{*}$}
\label{tab:4}      
\begin{tabular}{p{1.1cm}p{2.1cm}p{2.1cm}p{2.1cm}p{2.1cm}p{2.1cm}p{2.1cm}}
\hline\noalign{\smallskip}
$T-t$&$0.08222$&$0.05020$&$0.02310$&$0.01956$\\
\noalign{\smallskip}\hline\noalign{\smallskip}
$p$&8&21&96&133\\
\noalign{\smallskip}\hline\noalign{\smallskip}
$p_1$&1&1&3&4\\
\noalign{\smallskip}\hline\noalign{\smallskip}
$p^{*}$&23&61&286&398\\
\noalign{\smallskip}\hline\noalign{\smallskip}
$p_1^{*}$&1&2&4&5\\
\noalign{\smallskip}\hline\noalign{\smallskip}
\end{tabular}
\end{table}

\begin{table}
\centering
\caption{Confirmation of the formula (\ref{zzz4})}
\label{tab:5}      
\begin{tabular}{p{2.1cm}p{1.7cm}p{1.7cm}p{2.1cm}p{2.3cm}p{2.3cm}p{2.3cm}}
\hline\noalign{\smallskip}
$\varepsilon/(T-t)^3$&$0.0629$&$0.0097$&$0.0010$&$1.0129\cdot 10^{-4}$&
$1.0132\cdot 10^{-5}$\\
\noalign{\smallskip}\hline\noalign{\smallskip}
$q$&1&10&100&1000&10000\\
\noalign{\smallskip}\hline\noalign{\smallskip}
\end{tabular}
\vspace{5mm}
\end{table}

In Table 5, we can see the numerical confirmation of 
the formula (\ref{zzz4}) ($\varepsilon$ is a left-hand side 
of the formula (\ref{zzz4})).

\vspace{5mm}

\section{Comparative Analysis of the Efficiency of 
Application of Legendre Polynomials 
and Trigonometric Functions for the Integral
$J_{(011)T,t}^{*(0 i_1 i_2)}$}

\vspace{5mm}

In this section, we compare computational costs for the iterated Stratonovich
stochastic integral $J_{(011)T,t}^{*(0 i_1 i_2)}$ 
$(i_1, i_2=1,\ldots,m)$ within the frames
of the method of generalized multiple Fourier series for 
the systems of Legendre polynomials
and trigomomenric functions.

Using Theorem 3 for the case of trigonometric system of 
functions, we obtain
\cite{3}-\cite{17} (also see \cite{1}-\cite{2})

$$
J_{(011)T,t}^{*(0 i_1 i_2)q}=(T-t)^{2}\Biggl(\frac{1}{6}
\zeta_{0}^{(i_1)}\zeta_{0}^{(i_2)}-\frac{1}{2\sqrt{2}\pi}
\sqrt{\alpha_q}\xi_q^{(i_2)}\zeta_0^{(i_1)}+\Biggr.
$$

$$
+\frac{1}{2\sqrt{2}\pi^2}\sqrt{\beta_q}\Biggl(
\mu_q^{(i_2)}\zeta_0^{(i_1)}-2\mu_q^{(i_1)}\zeta_0^{(i_2)}\Biggr)+
$$

$$
+\frac{1}{2\sqrt{2}}\sum_{r=1}^{q}
\Biggl(-\frac{1}{\pi r}
\zeta_{2r-1}^{(i_2)}
\zeta_{0}^{(i_1)}+
\frac{1}{\pi^2 r^2}\left(
\zeta_{2r}^{(i_2)}
\zeta_{0}^{(i_1)}-
2\zeta_{2r}^{(i_1)}
\zeta_{0}^{(i_2)}\right)\Biggr)-
$$

$$
-
\frac{1}{2\pi^2}\sum_{\stackrel{r,l=1}{{}_{r\ne l}}}^{q}
\frac{1}{r^2-l^2}\Biggl(
\zeta_{2r}^{(i_1)}
\zeta_{2l}^{(i_2)}+
\frac{l}{r}
\zeta_{2r-1}^{(i_1)}
\zeta_{2l-1}^{(i_2)}
\Biggr)+
$$

$$
+
\sum_{r=1}^{q}
\Biggl(\frac{1}{4\pi r}\left(
\zeta_{2r}^{(i_1)}
\zeta_{2r-1}^{(i_2)}-
\zeta_{2r-1}^{(i_1)}
\zeta_{2r}^{(i_2)}\right)+\Biggr.
$$

\begin{equation}
\label{t1}
\Biggl.\Biggl.+
\frac{1}{8\pi^2 r^2}\left(
3\zeta_{2r-1}^{(i_1)}
\zeta_{2r-1}^{(i_2)}+
\zeta_{2r}^{(i_2)}
\zeta_{2r}^{(i_1)}\right)\Biggr)\Biggr).
\end{equation}

\vspace{6mm}

For the case $i_1\ne i_2$ from Theorem 8 we get 
\cite{3}-\cite{17} (also see \cite{1}-\cite{2})

\vspace{1mm}

$$
{\sf M}\left\{\left(J_{(011)T,t}^{*(0 i_1 i_2)}-
J_{(011)T,t}^{*(0 i_1 i_2)q}\right)^2\right\}=
\frac{(T-t)^4}{4}\Biggl(\frac{1}{9}-
\frac{1}{2\pi^2}\sum_{r=1}^q \frac{1}{r^2}-\Biggr.
$$

\vspace{1mm}
\begin{equation}
\label{t2}
\Biggl.-\frac{5}{8\pi^4}\sum_{r=1}^q \frac{1}{r^4}-
\frac{1}{\pi^4}\sum_{\stackrel{k,l=1}{{}_{k\ne l}}}^q
\frac{k^2+l^2}{l^2\left(l^2-k^2\right)^2}\Biggr).
\end{equation}

\vspace{5mm}

Analogues of the formulas (\ref{t1}), (\ref{t2}) for the case of 
Legendre polynomials will look as follows \cite{3}-\cite{17}
(also see \cite{1}-\cite{2})

\vspace{1mm}
$$
J_{(011)T,t}^{*(0 i_1 i_2)q}=\frac{T-t}{2}J_{(11)T,t}^{*(i_1 i_2)q}
+\frac{(T-t)^2}{4}\Biggl(
\frac{1}{\sqrt{3}}\zeta_0^{(i_2)}\zeta_1^{(i_1)}+\Biggr.
$$

\vspace{2mm}
\begin{equation}
\label{t3}
+\Biggl.\sum_{i=0}^{q}\Biggl(
\frac{(i+1)\zeta_{i+2}^{(i_2)}\zeta_{i}^{(i_1)}
-(i+2)\zeta_{i}^{(i_2)}\zeta_{i+2}^{(i_1)}}
{\sqrt{(2i+1)(2i+5)}(2i+3)}+
\frac{\zeta_i^{(i_1)}\zeta_{i}^{(i_2)}}{(2i-1)(2i+3)}\Biggr)\Biggr),
\end{equation}

\vspace{5mm}
\noindent
where

$$
J_{(11)T,t}^{*(i_1 i_2)q}=
\frac{T-t}{2}\Biggl(\zeta_0^{(i_1)}\zeta_0^{(i_2)}+\sum_{i=1}^{q}
\frac{1}{\sqrt{4i^2-1}}\left(
\zeta_{i-1}^{(i_1)}\zeta_{i}^{(i_2)}-
\zeta_i^{(i_1)}\zeta_{i-1}^{(i_2)}\right)\Biggr),
$$

\vspace{7mm}

$$
{\sf M}\left\{\left(J_{(011)T,t}^{*(0 i_1 i_2)}-J_{(011)T,t}^{*(0 i_1 i_2)q}
\right)^2\right\}=
$$

\vspace{2mm}
$$
=
\frac{(T-t)^4}{16}\left(\frac{5}{9}-
2\sum_{i=2}^q\frac{1}{4i^2-1}-
\sum_{i=1}^q
\frac{1}{(2i-1)^2(2i+3)^2}-
\right.
$$

\vspace{2mm}
\begin{equation}
\label{t4}
-\left.\sum_{i=0}^q\frac{(i+2)^2+(i+1)^2}{(2i+1)(2i+5)(2i+3)^2}
\right)\ \ \ (i_1\ne i_2).
\end{equation}

\vspace{6mm}

In Tables 6 and 7, we can see the numerical confirmation of 
the formulas (\ref{t2}) and (\ref{t4}) ($\varepsilon$ is 
the right-hand side of (\ref{t2}), (\ref{t4})).

\begin{table}
\centering
\caption{Confirmation of the formula (\ref{t2})}
\label{tab:6}      
\begin{tabular}{p{2.1cm}p{1.7cm}p{1.7cm}p{2.1cm}p{2.3cm}p{2.3cm}p{2.3cm}}
\hline\noalign{\smallskip}
$4\varepsilon/(T-t)^4$&0.0540&0.0082&$8.4261\cdot 10^{-4}$
&$8.4429\cdot 10^{-5}$&
$8.4435\cdot 10^{-6}$\\
\noalign{\smallskip}\hline\noalign{\smallskip}
$q$&1&10&100&1000&10000\\
\noalign{\smallskip}\hline\noalign{\smallskip}
\end{tabular}
\end{table}

\begin{table}
\centering
\caption{Confirmation of the formula (\ref{t4})}
\label{tab:7}      
\begin{tabular}{p{2.1cm}p{1.7cm}p{1.7cm}p{2.1cm}p{2.3cm}p{2.3cm}p{2.3cm}}
\hline\noalign{\smallskip}
$16\varepsilon/(T-t)^4$&0.3797&0.0581&0.0062&$6.2450\cdot 10^{-4}$&$6.2495\cdot 10^{-5}$\\
\noalign{\smallskip}\hline\noalign{\smallskip}
$q$&1&10&100&1000&10000\\
\noalign{\smallskip}\hline\noalign{\smallskip}
\end{tabular}
\end{table}

Let us compare the complexity of the formulas (\ref{t1}) and (\ref{t3}).
The formula (\ref{t1}) includes the double sum

$$
\frac{1}{2\pi^2}\sum_{\stackrel{r,l=1}{{}_{r\ne l}}}^{q}
\frac{1}{r^2-l^2}\Biggl(
\zeta_{2r}^{(i_1)}
\zeta_{2l}^{(i_2)}+
\frac{l}{r}
\zeta_{2r-1}^{(i_1)}
\zeta_{2l-1}^{(i_2)}
\Biggr).
$$

\vspace{3mm}

Thus, the formula (\ref{t1}) is more complex than the formula (\ref{t3})
even if we take identical numbers $q$ in these formulas.
As we noted above, the number $q$ in (\ref{t1}) must be equal to the
number $q$ from the formula (\ref{555}), so it is much 
larger than the number $q$ 
from the formula (\ref{t3}). As a result, we have 
obvious advantage of the formula (\ref{t3})
in computational costs. As we mentioned above,
if we will not use the random 
variables $\xi_q^{(i)}$ and 
$\mu_q^{(i)},$ then the number $q$ in (\ref{t1}) can be chosen smaller, but
the mean-square error of approximation of the
stochastic integral $J_{(11)T,t}^{(i_1 i_2)}$ will be three times larger
(see (\ref{801})). Moreover, in this case the stochastic integrals 
$J_{(01)T,t}^{(0 i_1)}$, $J_{(10)T,t}^{(i_1 0)},$ 
$J_{(001)T,t}^{(00i_1)}$ (with Gaussian distribution)
will be approximated worse. In this situation we can again talk about
the advantage of Ledendre polynomials.

\vspace{5mm}

\section{Conclusions}

\vspace{5mm}

Summing up the results of previous sections, we 
can come to the following conclusions.

\vspace{2mm}

1.\;We can talk about approximately equal computational costs
for the formulas (\ref{555}) and (\ref{401}). This means that
computational costs for implementing the Milstein scheme (explicit
one-step  
strong Taylor--Ito numerical method with the order $\gamma=1.0$ of accuracy 
for Ito stochastic differential equations \cite{19}) 
for the case of Legendre polynomials and for the case 
of trigonometric functions are approximately the same.

\vspace{2mm}

2.\;If we will not use the random 
variables $\xi_q^{(i)}$ (see (\ref{555})), then 
the mean-square error of approximation of the
stochastic integral $J_{(11)T,t}^{(i_1 i_2)}$ will be three times larger
(see (\ref{801})).
In this situation, we can talk about
the advantage of Ledendre polynomials in
the Milstein method.
Moreover, in this case the stochastic integrals 
$J_{(01)T,t}^{(0 i_1)}$, $J_{(10)T,t}^{(i_1 0)},$ 
$J_{(001)T,t}^{(00i_1)}$ (with Gaussian distribution)
will be approximated worse.

\vspace{2mm}

3.\;If we talk about the explicit one-step strong Taylor--Ito scheme 
with the order $\gamma=1.5$ of accuracy 
for Ito stochastic differential equations, then 
the numbers $q,$ $q_1$ (see (\ref{y3}), (\ref{y4})) 
are different. At that $q_1\ll q$
(the case of Legendre polynomials). 
The number $q$ must be the same in (\ref{x2})--(\ref{x5}) 
(the case of trigonometric functions).
This leads to huge computational costs (see very complex formula (\ref{x5})).
From the other hand, we can take different numbers $q$
in (\ref{x2})--(\ref{x5}). At that we should exclude
the random variables $\xi_q^{(i)},$ $\mu_q^{(i)}$ from 
(\ref{x2})--(\ref{x5}). This leads to another 
problems, which we discussed above (see Conclusion 1).

\vspace{2mm}

4.\;In addition, 
the author supposes that effect described in Conclusion 
3 will be more impressive when 
analyzing more complex sets of iterated Ito and Stratonovich 
stochastic
integrals (when $\gamma=$ $2.0,$ $2.5,$ $3.0,$ $\ldots $; here $\gamma$
has the same meaning as in Conclusion 3). 
This supposition is based on the fact that the polynomial 
system of functions has the significant advantage (compared with 
the trigonometric system) for approximation of iterated stochastic 
integrals for which not all weight functions are equal to 1.

\vspace{5mm}

\section{Further Development of Multiple Fourier--Legendre Series Approach to
the Mean-Square Approximation of Iterated Ito and Stratonovich
Stochastic Integrals of Multiplicities 3 to 5} 

\vspace{5mm}

From Theorems 1, 2 for $k=4$ and $5$ we obtain

\vspace{2mm}
$$
J_{(\lambda_1\lambda_2\lambda_3\lambda_4)T,t}^{(i_1 i_2 i_3 i_4)q}
=\sum_{j_1,j_2,j_3,j_4=0}^{q}
C_{j_4 j_3 j_2 j_1}\Biggl(
\zeta_{j_1}^{(i_1)}\zeta_{j_2}^{(i_2)}\zeta_{j_3}^{(i_3)}\zeta_{j_4}^{(i_4)}
\Biggr.
-
$$
$$
-
{\bf 1}_{\{i_1=i_2\ne 0\}}
{\bf 1}_{\{j_1=j_2\}}
\zeta_{j_3}^{(i_3)}
\zeta_{j_4}^{(i_4)}
-
{\bf 1}_{\{i_1=i_3\ne 0\}}
{\bf 1}_{\{j_1=j_3\}}
\zeta_{j_2}^{(i_2)}
\zeta_{j_4}^{(i_4)}
-
$$
$$
-
{\bf 1}_{\{i_1=i_4\ne 0\}}
{\bf 1}_{\{j_1=j_4\}}
\zeta_{j_2}^{(i_2)}
\zeta_{j_3}^{(i_3)}
-
{\bf 1}_{\{i_2=i_3\ne 0\}}
{\bf 1}_{\{j_2=j_3\}}
\zeta_{j_1}^{(i_1)}
\zeta_{j_4}^{(i_4)}-
$$
$$
-
{\bf 1}_{\{i_2=i_4\ne 0\}}
{\bf 1}_{\{j_2=j_4\}}
\zeta_{j_1}^{(i_1)}
\zeta_{j_3}^{(i_3)}
-
{\bf 1}_{\{i_3=i_4\ne 0\}}
{\bf 1}_{\{j_3=j_4\}}
\zeta_{j_1}^{(i_1)}
\zeta_{j_2}^{(i_2)}
+
$$
$$
+
{\bf 1}_{\{i_1=i_2\ne 0\}}
{\bf 1}_{\{j_1=j_2\}}
{\bf 1}_{\{i_3=i_4\ne 0\}}
{\bf 1}_{\{j_3=j_4\}}+
$$
$$
+
{\bf 1}_{\{i_1=i_3\ne 0\}}
{\bf 1}_{\{j_1=j_3\}}
{\bf 1}_{\{i_2=i_4\ne 0\}}
{\bf 1}_{\{j_2=j_4\}}+
$$
\begin{equation}
\label{leto5003a}
+\Biggl.
{\bf 1}_{\{i_1=i_4\ne 0\}}
{\bf 1}_{\{j_1=j_4\}}
{\bf 1}_{\{i_2=i_3\ne 0\}}
{\bf 1}_{\{j_2=j_3\}}\Biggr),
\end{equation}

\vspace{5mm}

$$
J_{(\lambda_1\lambda_2\lambda_3\lambda_4\lambda_5)T,t}^{(i_1 i_2 i_3 i_4 i_5)q_1}
=\sum_{j_1,j_2,j_3,j_4,j_5=0}^{q_1}
C_{j_5j_4 j_3 j_2 j_1}\Biggl(
\zeta_{j_1}^{(i_1)}\zeta_{j_2}^{(i_2)}\zeta_{j_3}^{(i_3)}\zeta_{j_4}^{(i_4)}
\zeta_{j_5}^{(i_5)}
\Biggr.
-
$$
$$
-
{\bf 1}_{\{i_1=i_2\ne 0\}}
{\bf 1}_{\{j_1=j_2\}}
\zeta_{j_3}^{(i_3)}
\zeta_{j_4}^{(i_4)}
\zeta_{j_5}^{(i_5)}-
{\bf 1}_{\{i_1=i_3\ne 0\}}
{\bf 1}_{\{j_1=j_3\}}
\zeta_{j_2}^{(i_2)}
\zeta_{j_4}^{(i_4)}
\zeta_{j_5}^{(i_5)}-
$$
$$
-
{\bf 1}_{\{i_1=i_4\ne 0\}}
{\bf 1}_{\{j_1=j_4\}}
\zeta_{j_2}^{(i_2)}
\zeta_{j_3}^{(i_3)}
\zeta_{j_5}^{(i_5)}-
{\bf 1}_{\{i_1=i_5\ne 0\}}
{\bf 1}_{\{j_1=j_5\}}
\zeta_{j_2}^{(i_2)}
\zeta_{j_3}^{(i_3)}
\zeta_{j_4}^{(i_4)}-
$$
$$
-{\bf 1}_{\{i_2=i_3\ne 0\}}
{\bf 1}_{\{j_2=j_3\}}
\zeta_{j_1}^{(i_1)}
\zeta_{j_4}^{(i_4)}
\zeta_{j_5}^{(i_5)}-
{\bf 1}_{\{i_2=i_4\ne 0\}}
{\bf 1}_{\{j_2=j_4\}}
\zeta_{j_1}^{(i_1)}
\zeta_{j_3}^{(i_3)}
\zeta_{j_5}^{(i_5)}-
$$
$$
-
{\bf 1}_{\{i_2=i_5\ne 0\}}
{\bf 1}_{\{j_2=j_5\}}
\zeta_{j_1}^{(i_1)}
\zeta_{j_3}^{(i_3)}\zeta_{j_4}^{(i_4)}-
{\bf 1}_{\{i_3=i_4\ne 0\}}
{\bf 1}_{\{j_3=j_4\}}
\zeta_{j_1}^{(i_1)}
\zeta_{j_2}^{(i_2)}
\zeta_{j_5}^{(i_5)}-
$$
$$
-
{\bf 1}_{\{i_3=i_5\ne 0\}}
{\bf 1}_{\{j_3=j_5\}}
\zeta_{j_1}^{(i_1)}
\zeta_{j_2}^{(i_2)}
\zeta_{j_4}^{(i_4)}
-{\bf 1}_{\{i_4=i_5\ne 0\}}
{\bf 1}_{\{j_4=j_5\}}
\zeta_{j_1}^{(i_1)}
\zeta_{j_2}^{(i_2)}
\zeta_{j_3}^{(i_3)}+
$$
$$
+
{\bf 1}_{\{i_1=i_2\ne 0\}}
{\bf 1}_{\{j_1=j_2\}}
{\bf 1}_{\{i_3=i_4\ne 0\}}
{\bf 1}_{\{j_3=j_4\}}\zeta_{j_5}^{(i_5)}+
{\bf 1}_{\{i_1=i_2\ne 0\}}
{\bf 1}_{\{j_1=j_2\}}
{\bf 1}_{\{i_3=i_5\ne 0\}}
{\bf 1}_{\{j_3=j_5\}}\zeta_{j_4}^{(i_4)}+
$$
$$
+
{\bf 1}_{\{i_1=i_2\ne 0\}}
{\bf 1}_{\{j_1=j_2\}}
{\bf 1}_{\{i_4=i_5\ne 0\}}
{\bf 1}_{\{j_4=j_5\}}\zeta_{j_3}^{(i_3)}+
{\bf 1}_{\{i_1=i_3\ne 0\}}
{\bf 1}_{\{j_1=j_3\}}
{\bf 1}_{\{i_2=i_4\ne 0\}}
{\bf 1}_{\{j_2=j_4\}}\zeta_{j_5}^{(i_5)}+
$$
$$
+
{\bf 1}_{\{i_1=i_3\ne 0\}}
{\bf 1}_{\{j_1=j_3\}}
{\bf 1}_{\{i_2=i_5\ne 0\}}
{\bf 1}_{\{j_2=j_5\}}\zeta_{j_4}^{(i_4)}+
{\bf 1}_{\{i_1=i_3\ne 0\}}
{\bf 1}_{\{j_1=j_3\}}
{\bf 1}_{\{i_4=i_5\ne 0\}}
{\bf 1}_{\{j_4=j_5\}}\zeta_{j_2}^{(i_2)}+
$$
$$
+
{\bf 1}_{\{i_1=i_4\ne 0\}}
{\bf 1}_{\{j_1=j_4\}}
{\bf 1}_{\{i_2=i_3\ne 0\}}
{\bf 1}_{\{j_2=j_3\}}\zeta_{j_5}^{(i_5)}+
{\bf 1}_{\{i_1=i_4\ne 0\}}
{\bf 1}_{\{j_1=j_4\}}
{\bf 1}_{\{i_2=i_5\ne 0\}}
{\bf 1}_{\{j_2=j_5\}}\zeta_{j_3}^{(i_3)}+
$$
$$
+
{\bf 1}_{\{i_1=i_4\ne 0\}}
{\bf 1}_{\{j_1=j_4\}}
{\bf 1}_{\{i_3=i_5\ne 0\}}
{\bf 1}_{\{j_3=j_5\}}\zeta_{j_2}^{(i_2)}+
{\bf 1}_{\{i_1=i_5\ne 0\}}
{\bf 1}_{\{j_1=j_5\}}
{\bf 1}_{\{i_2=i_3\ne 0\}}
{\bf 1}_{\{j_2=j_3\}}\zeta_{j_4}^{(i_4)}+
$$
$$
+
{\bf 1}_{\{i_1=i_5\ne 0\}}
{\bf 1}_{\{j_1=j_5\}}
{\bf 1}_{\{i_2=i_4\ne 0\}}
{\bf 1}_{\{j_2=j_4\}}\zeta_{j_3}^{(i_3)}+
{\bf 1}_{\{i_1=i_5\ne 0\}}
{\bf 1}_{\{j_1=j_5\}}
{\bf 1}_{\{i_3=i_4\ne 0\}}
{\bf 1}_{\{j_3=j_4\}}\zeta_{j_2}^{(i_2)}+
$$
$$
+
{\bf 1}_{\{i_2=i_3\ne 0\}}
{\bf 1}_{\{j_2=j_3\}}
{\bf 1}_{\{i_4=i_5\ne 0\}}
{\bf 1}_{\{j_4=j_5\}}\zeta_{j_1}^{(i_1)}+
{\bf 1}_{\{i_2=i_4\ne 0\}}
{\bf 1}_{\{j_2=j_4\}}
{\bf 1}_{\{i_3=i_5\ne 0\}}
{\bf 1}_{\{j_3=j_5\}}\zeta_{j_1}^{(i_1)}+
$$
\begin{equation}
\label{leto5004}
+\Biggl.
{\bf 1}_{\{i_2=i_5\ne 0\}}
{\bf 1}_{\{j_2=j_5\}}
{\bf 1}_{\{i_3=i_4\ne 0\}}
{\bf 1}_{\{j_3=j_4\}}\zeta_{j_1}^{(i_1)}\Biggr),
\end{equation}

\vspace{5mm}
\noindent
where 

\vspace{-2mm}
$$
J_{(\lambda_1\lambda_2\lambda_3\lambda_4)T,t}^{(i_1 i_2 i_3 i_4)},\ \ \
J_{(\lambda_1\lambda_2\lambda_3\lambda_4\lambda_5)
T,t}^{(i_1 i_2 i_3 i_4 i_5)}
$$

\vspace{3mm}
\noindent
are defined by the formula (\ref{ito1});
$q_1<q;$ 
${\bf 1}_A$ is the indicator of the 
set $A;$ $i_1, i_2, i_3, i_4, i_5 =0, 1,\ldots,m,$ and

\vspace{1mm}
$$
C_{j_4j_3j_2j_1}=\int\limits_{t}^{T}\phi_{j_4}(u)
\int\limits_{t}^{u}\phi_{j_3}(z)
\int\limits_{t}^{z}\phi_{j_2}(y)
\int\limits_{t}^{y}
\phi_{j_1}(x)dx dy dz du=
$$

\vspace{1mm}
$$
=\sqrt{(2j_1+1)(2j_2+1)(2j_3+1)(2j_4+1)}\ \frac{(T-t)^{2}}{16}\
\bar
C_{j_4j_3j_2j_1},
$$

\vspace{3mm}
$$
C_{j_5j_4j_3j_2j_1}=\int\limits_{t}^{T}\phi_{j_5}(v)
\int\limits_{t}^{v}\phi_{j_4}(u)
\int\limits_{t}^{u}\phi_{j_3}(z)
\int\limits_{t}^{z}\phi_{j_2}(y)
\int\limits_{t}^{y}
\phi_{j_1}(x)dx dy dz du dv
=
$$

\vspace{1mm}
$$
=\sqrt{(2j_1+1)(2j_2+1)(2j_3+1)(2j_4+1)(2j_5+1)}\ \frac{(T-t)^{5/2}}{32}\ 
\bar
C_{j_5j_4j_3j_2j_1},
$$

\vspace{3mm}
\begin{equation}
\label{204}
\bar C_{j_4j_3j_2j_1}=\int\limits_{-1}^{1}P_{j_4}(u)
\int\limits_{-1}^{u}P_{j_3}(z)
\int\limits_{-1}^{z}P_{j_2}(y)
\int\limits_{-1}^{y}
P_{j_1}(x)dx dy dz du,
\end{equation}

\vspace{3mm}
\begin{equation}
\label{205}
\bar C_{j_5j_4j_3j_2j_1}=
\int\limits_{-1}^{1}P_{j_5}(v)
\int\limits_{-1}^{v}P_{j_4}(u)
\int\limits_{-1}^{u}P_{j_3}(z)
\int\limits_{-1}^{z}P_{j_2}(y)
\int\limits_{-1}^{y}
P_{j_1}(x)dx dy dz du dv,
\end{equation}

\vspace{5mm}

\noindent
where $P_i(x)$ $(i=0, 1, 2,\ldots)$ is the Legendre polynomial,

\vspace{-1mm}
$$
\zeta_{j}^{(i)}=
\int\limits_t^T \phi_{j}(s) d{\bf w}_s^{(i)}
$$

\vspace{3mm}
\noindent
are independent standard Gaussian random variables
for various
$i$ or $j$ (if $i\ne 0$),

\vspace{2mm}
$$
\phi_j(x)=\sqrt{\frac{2j+1}{T-t}}P_j\left(\left(
x-\frac{T+t}{2}\right)\frac{2}{T-t}\right),\ \ \ j=0, 1, 2,\ldots
$$

\vspace{4mm}
\noindent
is a complete orthonormal system
of Legendre polynomials in the space $L_2([t, T]).$

Note that the Fourier--Legendre coefficients $\bar C_{j_4j_3j_2 j_1},$
$\bar C_{j_5j_4j_3j_2 j_1},$ and 
$\bar C_{j_3j_2 j_1}$ (see (\ref{ww}))
can be calculated exactly 
using DERIVE or MAPLE (computer algebra systems). 
Several tables with these coefficients can be found in \cite{3}-\cite{10aaa}, \cite{13},
\cite{15f}.
The database with $270,000$ of exactly calculated Fourier--Legendre 
coefficients is descibed in \cite{Kuz-Kuz}, \cite{Mikh-1}.
Note that the mentioned Fourier--Legendre coefficients
not depend on the integration step $T-t$ of
numerical methods for Ito stochastic differential equations.
So, $T-t$ can be not a constant in this approach.

From (\ref{star00011}) ($0<T-t<1$) we obtain

\begin{equation}
\label{12}
{\sf M}\left\{\left(
J_{(\lambda_1\lambda_2\lambda_3\lambda_4)T,t}^{(i_1 i_2 i_3 i_4)}-
J_{(\lambda_1\lambda_2\lambda_3\lambda_4)T,t}^{(i_1 i_2 i_3 i_4)q}\right)^2\right\}\le
24\left(\frac{(T-t)^{4}}{24}-\sum_{j_1,j_2,j_3,j_4=0}^{q}
C_{j_4j_3j_2j_1}^2\right),
\end{equation}

\begin{equation}
\label{13}
{\sf M}\left\{\left(
J_{(\lambda_1\lambda_2\lambda_3\lambda_4\lambda_5)T,t}^{(i_1 i_2 i_3 i_4 i_5)}-
J_{(\lambda_1\lambda_2\lambda_3\lambda_4\lambda_5)T,t}^{(i_1 i_2 i_3 i_4 i_5)q_1}
\right)^2\right\}\le
120\left(\frac{(T-t)^{5}}{120}-\sum_{j_1,j_2,j_3,j_4,j_5=0}^{q_1}
C_{j_5j_4j_3j_2j_1}^2\right).
\end{equation}

\vspace{5mm}

Note that in practice the numbers $q, q_1$ in (\ref{y4}), 
(\ref{leto5003a}), (\ref{leto5004}) can be selected not large.
For example, for the case of pairwise different
$i_1, i_2, i_3, i_4, i_5=1,\ldots,m$ we obtain

\vspace{1mm}
\begin{equation}
\label{800a}
{\sf M}\left\{\left(
J_{(111)T,t}^{(i_1i_2 i_3)}-
J_{(111)T,t}^{(i_1i_2 i_3)6}\right)^2\right\}=
\frac{(T-t)^{3}}{6}-\sum_{j_1,j_2,j_3=0}^{6}
C_{j_3j_2j_1}^2
\approx
0.01956000(T-t)^3,
\end{equation}

\vspace{2mm}
\begin{equation}
\label{804a}
{\sf M}\left\{\left(
J_{(1111)T,t}^{(i_1i_2 i_3 i_4)}-
J_{(1111)T,t}^{(i_1i_2 i_3 i_4)2}\right)^2\right\}=
\frac{(T-t)^{4}}{24}-\sum_{j_1,j_2,j_3,j_4=0}^{2}
C_{j_4j_3j_2j_1}^2\approx
0.02360840(T-t)^4,
\end{equation}

\vspace{2mm}
\begin{equation}
\label{805a}
{\sf M}\left\{\left(
J_{(11111)T,t}^{(i_1i_2 i_3 i_4 i_5)}-
J_{(11111)T,t}^{(i_1i_2 i_3 i_4 i_5)1}\right)^2\right\}=
\frac{(T-t)^5}{120}-\sum_{j_1,j_2,j_3,j_4,j_5=0}^{1}
C_{j_5j_4j_3j_2j_1}^2\approx
0.00759105(T-t)^5.
\end{equation}

\vspace{7mm}

From Theorems 3--6 we have

\vspace{1mm}
$$
J_{(\lambda_1\lambda_2,\lambda_3)T,t}^{*(i_1i_2 i_3)q}
=\sum_{j_1,j_2,j_3=0}^{q}
C_{j_3j_2j_1}
\zeta_{j_1}^{(i_1)}\zeta_{j_2}^{(i_2)}\zeta_{j_3}^{(i_3)},
$$

\vspace{2mm}
$$
J_{(\lambda_1\lambda_2,\lambda_3\lambda_4)T,t}^{*(i_1i_2 i_3i_4)q_1}
=\sum_{j_1,j_2,j_3,j_4=0}^{q_1}
C_{j_4j_3j_2j_1}
\zeta_{j_1}^{(i_1)}\zeta_{j_2}^{(i_2)}\zeta_{j_3}^{(i_3)}\zeta_{j_4}^{(i_4)},
$$

\vspace{2mm}
$$
J_{(\lambda_1\lambda_2,\lambda_3\lambda_4\lambda_5)T,t}^{*(i_1i_2 i_3i_4i_5)q_2}
=\sum_{j_1,j_2,j_3,j_4,j_5=0}^{q_2}
C_{j_5j_4j_3j_2j_1}
\zeta_{j_1}^{(i_1)}\zeta_{j_2}^{(i_2)}\zeta_{j_3}^{(i_3)}\zeta_{j_4}^{(i_4)}\zeta_{j_5}^{(i_5)},
$$

\vspace{6mm}
\noindent
where

\vspace{-1mm}
$$
J_{(\lambda_1\lambda_2\lambda_3)T,t}^{*(i_1i_2 i_3)},\
J_{(\lambda_1\lambda_2\lambda_3\lambda_4)T,t}^{*(i_1i_2 i_3i_4)},\
J_{(\lambda_1\lambda_2\lambda_3\lambda_4\lambda_5)T,t}^{*(i_1i_2 i_3i_4i_5)}
$$

\vspace{5mm}
\noindent
are defined by the formula (\ref{str1}).

The values

$$
{\sf M}\left\{\left(
J_{(111)T,t}^{*(i_1i_2 i_3)}-
J_{(111)T,t}^{*(i_1i_2 i_3)6}\right)^2\right\},\ \ \
{\sf M}\left\{\left(
J_{(1111)T,t}^{*(i_1i_2 i_3 i_4)}-
J_{(1111)T,t}^{*(i_1i_2 i_3 i_4)2}\right)^2\right\},
$$

\vspace{2mm}
$$
{\sf M}\left\{\left(
J_{(11111)T,t}^{*(i_1i_2 i_3 i_4 i_5)}-
J_{(11111)T,t}^{*(i_1i_2 i_3 i_4 i_5)1}\right)^2\right\}
$$

\vspace{5mm}
\noindent
are equal to the right-hand sides of (\ref{800a})--(\ref{805a})
for the case of pairwise different $i_1, i_2, i_3, i_4, i_5=1,\ldots,m.$

Note that the optimization of the mean-square approximation
procedures for the itertaed Ito stochastic integrals (\ref{ito}) of multiplicities
1 to 5 is carried out in \cite{Mikh-2}, \cite{Mikh-2aaaaa}.

\vspace{5mm}

\section{Theorems 1--7 from Point
of View of the Wong--Zakai Approximation}

\vspace{5mm}

The iterated Ito stochastic integrals and solutions
of Ito SDEs are complex and important functionals
from the independent components ${\bf f}_{s}^{(i)},$
$i=1,\ldots,m$ of the multidimensional
Wiener process ${\bf f}_{s},$ $s\in[0, T].$
Let ${\bf f}_{s}^{(i)p},$ $p\in\mathbb{N}$ 
be some approximation of
${\bf f}_{s}^{(i)},$
$i=1,\ldots,m$.
Suppose that 
${\bf f}_{s}^{(i)p}$
converges to
${\bf f}_{s}^{(i)},$
$i=1,\ldots,m$ if $p\to\infty$ in some sense and has
differentiable sample trajectories.

A natural question arises: if we replace 
${\bf f}_{s}^{(i)}$
by ${\bf f}_{s}^{(i)p},$
$i=1,\ldots,m$ in the functionals
mentioned above, will the resulting
functionals converge to the original
functionals from the components 
${\bf f}_{s}^{(i)},$
$i=1,\ldots,m$ of the multidimentional
Wiener process ${\bf f}_{s}$?
The answere to this question is negative 
in the general case. However, 
in the pioneering works of Wong E. and Zakai M. \cite{W-Z-1},
\cite{W-Z-2},
it was shown that under the special conditions and 
for some types of approximations 
of the Wiener process the answere is affirmative
with one peculiarity: the convergence takes place 
to the iterated Stratonovich stochastic integrals
and solutions of Stratonovich SDEs and not to iterated 
Ito stochastic integrals and solutions
of Ito stochastic differential equations.
The piecewise 
linear approximation 
as well as the regularization by convolution 
\cite{W-Z-1}-\cite{Watanabe} relate the 
mentioned types of approximations
of the Wiener process. The above approximation 
of stochastic integrals and solutions of SDEs 
is often called the Wong--Zakai approximation.

Let ${\bf w}_{\tau},$ $\tau\in[0, T]$ is a random vector with 
an $m+1$ components: ${\bf w}_{\tau}^{(i)}={\bf f}_{\tau}^{(i)}$ 
for $i=1,\ldots,m$ and 
${\bf w}_{\tau}^{(0)}=\tau,$\ 
${\bf f}_{\tau}^{(i)}$ $(i=1,\ldots,m)$
are independent standard Wiener processes.

It is well known that the following representation 
takes place \cite{Lipt}, \cite{7e}

\begin{equation}
\label{um1x}
{\bf w}_{\tau}^{(i)}-{\bf w}_{t}^{(i)}=
\sum_{j=0}^{\infty}\int\limits_t^{\tau}
\phi_j(s)ds\ \zeta_j^{(i)},\ \ \ \zeta_j^{(i)}=
\int\limits_t^T \phi_j(s)d{\bf w}_s^{(i)},
\end{equation}

\vspace{4mm}
\noindent
where $\tau\in[t, T],$ $t\ge 0,$
$\{\phi_j(x)\}_{j=0}^{\infty}$ is an arbitrary complete 
orthonormal system of functions in the space $L_2([t, T]),$ and
$\zeta_j^{(i)}$ are independent standard Gaussian 
random variables for various $i$ or $j.$
Moreover, the series (\ref{um1x}) converges for any $\tau\in [t, T]$
in the mean-square sense.

Let ${\bf w}_{\tau}^{(i)p}-{\bf w}_{t}^{(i)p}$ be 
the mean-square approximation of the process
${\bf w}_{\tau}^{(i)}-{\bf w}_{t}^{(i)},$
which has the following form

\vspace{-3mm}
\begin{equation}
\label{um1xx}
{\bf w}_{\tau}^{(i)p}-{\bf w}_{t}^{(i)p}=
\sum_{j=0}^{p}\int\limits_t^{\tau}
\phi_j(s)ds\ \zeta_j^{(i)}.
\end{equation}

\vspace{3mm}

From (\ref{um1xx}) we obtain

\vspace{-4mm}
\begin{equation}
\label{um1xxx}
d{\bf w}_{\tau}^{(i)p}=
\sum_{j=0}^{p}
\phi_j(\tau)\zeta_j^{(i)} d\tau.
\end{equation}

\vspace{4mm}

Consider the following iterated Riemann--Stieltjes
integral

\begin{equation}
\label{um1xxxx}
\int\limits_t^T
\psi_k(t_k)\ldots \int\limits_t^{t_2}\psi_1(t_1)
d{\bf w}_{t_1}^{(i_1)p_1}\ldots d{\bf w}_{t_k}^{(i_k)p_k},
\end{equation}

\vspace{4mm}
\noindent
where $p_1,\ldots,p_k\in\mathbb{N},$\ \ $i_1,\ldots,i_k=0,1,\ldots,m,$ 

\begin{equation}
\label{um1xxx1}
d{\bf w}_{\tau}^{(i)p}=
\left\{\begin{matrix}
d{\bf f}_{\tau}^{(i)p}\ &\hbox{\rm for}\ \ \ i=1,\ldots,m\cr\cr\cr
d\tau^p\ &\hbox{\rm for}\ \ \ i=0
\end{matrix}
,\right.
\end{equation}

\vspace{4mm}
\noindent
and $d{\bf f}_{\tau}^{(i)p},$ $d\tau^p$ are defined by the relation (\ref{um1xxx}).

Let us substitute (\ref{um1xxx}) into (\ref{um1xxxx})

\begin{equation}
\label{um1xxxx1}
\int\limits_t^T
\psi_k(t_k)\ldots \int\limits_t^{t_2}\psi_1(t_1)
d{\bf w}_{t_1}^{(i_1)p_1}\ldots d{\bf w}_{t_k}^{(i_k)p_k}=
\sum\limits_{j_1=0}^{p_1}\ldots \sum\limits_{j_k=0}^{p_k}
C_{j_k \ldots j_1}\prod\limits_{l=1}^k \zeta_{j_l}^{(i_l)},
\end{equation}

\vspace{4mm}
\noindent
where 
$$
\zeta_j^{(i)}=\int\limits_t^T \phi_j(s)d{\bf w}_s^{(i)}
$$ 

\vspace{2mm}
\noindent
are independent standard Gaussian random variables for various 
$i$ or $j$ (in the case when $i\ne 0$),
${\bf w}_{s}^{(i)}={\bf f}_{s}^{(i)}$ for
$i=1,\ldots,m$ and 
${\bf w}_{s}^{(0)}=s,$

$$
C_{j_k \ldots j_1}=\int\limits_t^T\psi_k(t_k)\phi_{j_k}(t_k)\ldots
\int\limits_t^{t_2}
\psi_1(t_1)\phi_{j_1}(t_1)
dt_1\ldots dt_k
$$

\vspace{4mm}
\noindent
is the Fourier coefficient.

To best of our knowledge \cite{W-Z-1}-\cite{Watanabe}
the approximations of the Wiener process
in the Wong--Zakai approximation must satisfy fairly strong
restrictions
\cite{Watanabe}
(see Definition 7.1, pp.~480--481).
Moreover, approximations of the Wiener process that are
similar to (\ref{um1xx})
were not considered in \cite{W-Z-1}, \cite{W-Z-2}
(also see \cite{Watanabe}, Theorems 7.1, 7.2).
Therefore, the proof of analogs of Theorems 7.1 and 7.2 \cite{Watanabe}
for approximations of the Wiener 
process based on its series expansion (\ref{um1x})
should be carried out separately.
Thus, the mean-square convergence of the right-hand side
of (\ref{um1xxxx1}) to the iterated Stratonovich stochastic integral 
(\ref{str})
does not follow from the results of the papers
\cite{W-Z-1}, \cite{W-Z-2} (also see \cite{Watanabe},
Theorems 7.1, 7.2).

From the other hand, Theorems 1--7 from this 
paper can be considered as the proof of the
Wong--Zakai approximation for the iterated 
Stratonovich stochastic integrals (\ref{str}) of multiplicities 1 to 6
based on the approximation (\ref{um1xx}) of the Wiener process.
At that, the Riemann--Stieltjes integrals (\ref{um1xxxx}) converge
(according to Theorems 1--7)
to the appropriate Stratonovich 
stochastic integrals (\ref{str}). Recall that
$\{\phi_j(x)\}_{j=0}^{\infty}$ (see (\ref{um1x}), (\ref{um1xx}), and
Theorems 3--7)
is a complete 
orthonormal system of Legendre polynomials or 
trigonometric functions 
in the space $L_2([t, T])$.

To illustrate the above reasoning, 
consider two examples for the case $k=2,$
$\psi_1(s),$ $\psi_2(s)\equiv 1;$ $i_1, i_2=1,\ldots,m.$

The first example relates to the piecewise linear approximation
of the multidimensional Wiener process (these approximations 
were considered in \cite{W-Z-1}-\cite{Watanabe}).

Let ${\bf b}_{\Delta}^{(i)}(t),$ $t\in[0, T]$ be the piecewise
linear approximation of the $i$th component ${\bf f}_t^{(i)}$
of the multidimensional standard Wiener process ${\bf f}_t,$
$t\in [0, T]$ with independent components
${\bf f}_t^{(i)},$ $i=1,\ldots,m,$ i.e.

$$
{\bf b}_{\Delta}^{(i)}(t)={\bf f}_{k\Delta}^{(i)}+
\frac{t-k\Delta}{\Delta}\Delta{\bf f}_{k\Delta}^{(i)},
$$

\vspace{3mm}
\noindent
where 

\vspace{-2mm}
$$
\Delta{\bf f}_{k\Delta}^{(i)}={\bf f}_{(k+1)\Delta}^{(i)}-
{\bf f}_{k\Delta}^{(i)},\ \ \
t\in[k\Delta, (k+1)\Delta),\ \ \ k=0, 1,\ldots, N-1.
$$

\vspace{4mm}

Note that w.~p.~1

\vspace{-1mm}
\begin{equation}
\label{pridum}
\frac{d{\bf b}_{\Delta}^{(i)}}{dt}(t)=
\frac{\Delta{\bf f}_{k\Delta}^{(i)}}{\Delta},\ \ \
t\in[k\Delta, (k+1)\Delta),\ \ \ k=0, 1,\ldots, N-1.
\end{equation}

\vspace{4mm}

Consider the following iterated Riemann--Stieltjes
integral

\vspace{1mm}
$$
\int\limits_0^T
\int\limits_0^{s}
d{\bf b}_{\Delta}^{(i_1)}(\tau)d{\bf b}_{\Delta}^{(i_2)}(s),\ \ \ 
i_1,i_2=1,\ldots,m.
$$

\vspace{4mm}

Using (\ref{pridum}) and additive property of the Riemann--Stieltjes integrals, 
we can write w.~p.~1

\vspace{2mm}
$$
\int\limits_0^T
\int\limits_0^{s}
d{\bf b}_{\Delta}^{(i_1)}(\tau)d{\bf b}_{\Delta}^{(i_2)}(s)=
\int\limits_0^T
\int\limits_0^{s}
\frac{d{\bf b}_{\Delta}^{(i_1)}}{d\tau}(\tau)d\tau
\frac{d {\bf b}_{\Delta}^{(i_2)}}{d s}(s)
ds =
$$

\vspace{3mm}
$$
=
\sum\limits_{l=0}^{N-1}\int\limits_{l\Delta}^{(l+1)\Delta}
\left(
\sum\limits_{q=0}^{l-1}\int\limits_{q\Delta}^{(q+1)\Delta}
\frac{\Delta{\bf f}_{q\Delta}^{(i_1)}}{\Delta}d\tau+
\int\limits_{l\Delta}^{s}
\frac{\Delta{\bf f}_{l\Delta}^{(i_1)}}{\Delta}d\tau\right)
\frac{\Delta{\bf f}_{l\Delta}^{(i_2)}}{\Delta}ds=
$$

\vspace{3mm}
$$
=\sum\limits_{l=0}^{N-1}\sum\limits_{q=0}^{l-1}
\Delta{\bf f}_{q\Delta}^{(i_1)}
\Delta{\bf f}_{l\Delta}^{(i_2)}+
\frac{1}{\Delta^2}\sum\limits_{l=0}^{N-1}
\Delta{\bf f}_{l\Delta}^{(i_1)}
\Delta{\bf f}_{l\Delta}^{(i_2)}
\int\limits_{l\Delta}^{(l+1)\Delta}
\int\limits_{l\Delta}^{s}d\tau ds=
$$

\vspace{3mm}
\begin{equation}
\label{oh-ty}
=\sum\limits_{l=0}^{N-1}\sum\limits_{q=0}^{l-1}
\Delta{\bf f}_{q\Delta}^{(i_1)}
\Delta{\bf f}_{l\Delta}^{(i_2)}+
\frac{1}{2}\sum\limits_{l=0}^{N-1}
\Delta{\bf f}_{l\Delta}^{(i_1)}
\Delta{\bf f}_{l\Delta}^{(i_2)}.
\end{equation}

\vspace{6mm}

Using (\ref{oh-ty}) it 
is not difficult to show 
that

\vspace{1mm}
$$
\hbox{\vtop{\offinterlineskip\halign{
\hfil#\hfil\cr
{\rm l.i.m.}\cr
$\stackrel{}{{}_{N\to \infty}}$\cr
}} }
\int\limits_0^T
\int\limits_0^{s}
d{\bf b}_{\Delta}^{(i_1)}(\tau)d{\bf b}_{\Delta}^{(i_2)}(s)=
\int\limits_0^T
\int\limits_0^{s}
d{\bf f}_{\tau}^{(i_1)}d{\bf f}_{s}^{(i_2)}+
\frac{1}{2}{\bf 1}_{\{i_1=i_2\}}\int\limits_0^T ds=
$$

\vspace{3mm}
\begin{equation}
\label{uh-111}
=
\int\limits_0^{*T}
\int\limits_0^{*s}
d{\bf f}_{\tau}^{(i_1)}d{\bf f}_{s}^{(i_2)},
\end{equation}

\vspace{5mm}
\noindent
where $\Delta\to 0$ if $N\to\infty$ ($N\Delta=T$).

Obviously, (\ref{uh-111}) agrees with Theorem 7.1 (see \cite{Watanabe},
p.~486).

The next example relates to the approximation
of the Wiener process based on its series expansion
(\ref{um1x}) for $t=0$, where
$\{\phi_j(x)\}_{j=0}^{\infty}$ 
is a complete 
orthonormal system of Legendre polynomials or 
trigonometric functions 
in the space $L_2([0, T])$.

Consider the following iterated Riemann--Stieltjes
integral

\vspace{-1mm}
\begin{equation}
\label{abcd1}
\int\limits_0^T
\int\limits_0^{s}
d{\bf f}_{\tau}^{(i_1)p}d{\bf f}_{s}^{(i_2)p},\ \ \ 
i_1,i_2=1,\ldots,m,
\end{equation}

\vspace{3mm}
\noindent
where $d{\bf f}_{\tau}^{(i)p}$ is defined by the
relation
(\ref{um1xxx}).

Let us substitute (\ref{um1xxx}) into (\ref{abcd1}) 

\vspace{-1mm}
\begin{equation}
\label{set18}
\int\limits_0^T
\int\limits_0^{s}
d{\bf f}_{\tau}^{(i_1)p}d{\bf f}_{s}^{(i_2)p}=
\sum\limits_{j_1,j_2=0}^p
C_{j_2 j_1} \zeta_{j_1}^{(i_1)}\zeta_{j_2}^{(i_2)},
\end{equation}

\vspace{3mm}
\noindent
where 
$$
C_{j_2 j_1}=
\int\limits_0^T \phi_{j_2}(s)\int\limits_0^s
\phi_{j_1}(\tau)d\tau ds
$$

\vspace{3mm}
\noindent
is the Fourier coefficient; another notations 
are the same as in (\ref{um1xxxx1}).

As we noted above, approximations of the Wiener process that are
similar to (\ref{um1xx})
were not considered in \cite{W-Z-1}, \cite{W-Z-2}
(also see Theorems 7.1, 7.2 in \cite{Watanabe}).
Furthermore, the extension of the results of Theorems 7.1 and 7.2
\cite{Watanabe} to the case under consideration is
not obvious.

On the other hand, we can apply the theory built in Chapters 1 and 2
of the monographs \cite{10a}-\cite{10aaa}. More precisely, 
using 
Theorem 3, we obtain from (\ref{set18}) the desired result

\vspace{-1mm}
\begin{equation}
\label{umen-bl}
\hbox{\vtop{\offinterlineskip\halign{
\hfil#\hfil\cr
{\rm l.i.m.}\cr
$\stackrel{}{{}_{p\to \infty}}$\cr
}} }
\int\limits_0^T
\int\limits_0^{s}
d{\bf f}_{\tau}^{(i_1)p}d{\bf f}_{s}^{(i_2)p}=
\hbox{\vtop{\offinterlineskip\halign{
\hfil#\hfil\cr
{\rm l.i.m.}\cr
$\stackrel{}{{}_{p\to \infty}}$\cr
}} }
\sum\limits_{j_1,j_2=0}^p
C_{j_2 j_1} \zeta_{j_1}^{(i_1)}\zeta_{j_2}^{(i_2)}=
\int\limits_0^{*T}
\int\limits_0^{*s}
d{\bf f}_{\tau}^{(i_1)}d{\bf f}_{s}^{(i_2)}.
\end{equation}

\vspace{5mm}

From the other hand, by Theorems 1, 2
(see (\ref{leto5001})) for the case
$k=2$ we obtain from (\ref{set18}) the following relation

\vspace{-2mm}
$$
\hbox{\vtop{\offinterlineskip\halign{
\hfil#\hfil\cr
{\rm l.i.m.}\cr
$\stackrel{}{{}_{p\to \infty}}$\cr
}} }
\int\limits_0^T
\int\limits_0^{s}
d{\bf f}_{\tau}^{(i_1)p}d{\bf f}_{s}^{(i_2)p}=
\hbox{\vtop{\offinterlineskip\halign{
\hfil#\hfil\cr
{\rm l.i.m.}\cr
$\stackrel{}{{}_{p\to \infty}}$\cr
}} }
\sum\limits_{j_1,j_2=0}^p
C_{j_2 j_1} \zeta_{j_1}^{(i_1)}\zeta_{j_2}^{(i_2)}=
$$

\vspace{2mm}
$$
=
\hbox{\vtop{\offinterlineskip\halign{
\hfil#\hfil\cr
{\rm l.i.m.}\cr
$\stackrel{}{{}_{p\to \infty}}$\cr
}} }
\sum\limits_{j_1,j_2=0}^p
C_{j_2 j_1} \biggl(\zeta_{j_1}^{(i_1)}\zeta_{j_2}^{(i_2)}-
{\bf 1}_{\{i_1=i_2\}}{\bf 1}_{\{j_1=j_2\}}\biggr)+
{\bf 1}_{\{i_1=i_2\}}\sum\limits_{j_1=0}^{\infty}
C_{j_1 j_1}=
$$

\vspace{2mm}
\begin{equation}
\label{umen-blx}
=
\int\limits_0^T
\int\limits_0^{s}
d{\bf f}_{\tau}^{(i_1)}d{\bf f}_{s}^{(i_2)}+
{\bf 1}_{\{i_1=i_2\}}\sum\limits_{j_1=0}^{\infty}
C_{j_1 j_1}.
\end{equation}

\vspace{5mm}

Since
$$
\sum\limits_{j_1=0}^{\infty}
C_{j_1 j_1}=\frac{1}{2}\sum\limits_{j_1=0}^{\infty}
\left(\int\limits_0^T \phi_j(\tau)d\tau\right)^2
=\frac{1}{2}
\left(\int\limits_0^T \phi_0(\tau)d\tau\right)^2=\frac{1}{2}
\int\limits_0^T ds,
$$

\vspace{4mm}
\noindent
then from (\ref{umen-blx}) we obtain (\ref{umen-bl}).

\end{document}